\documentclass[12pt]{article}

\usepackage{a4wide}
\usepackage{amssymb}
\usepackage{amsfonts}
\usepackage{amsmath}

\date{}

\newtheorem{proposition}{Proposition}[section]
\newtheorem{theorem}[proposition]{Theorem}
\newtheorem{lemma}[proposition]{Lemma}

\newtheorem{corollary}[proposition]{Corollary}

\def\fd{ {\rm fd}}
\def\lfd{ {\rm lfd}}
\def\fdim{ {\rm d}}

\def\GK{{\rm  GK}\,}
\def\Kdim{{\rm K.dim }\,}

\def\Hom{{\rm Hom}}
\def\der{\partial }
\def\deri{{\partial }_i}

\def\nFM0{{\nu }_{F,M_0}}
\def\nFN0{{\nu }_{F,N_0}}
\def\nGN0{{\nu }_{G,N_0}}

\def\nF{{\nu }_{F}}

\def\nG{ {\nu }_{G} }

\def\N0{ {\bf N}_0 }

\def\t{\otimes}
\def\g{\gamma}

\def\ra{\rightarrow}

\def\Xpm{X^{\pm }}

\def\l1{{\lambda}_1}

\def\p{{\bf p}}

\def\a{\alpha}
\def\a0{ {\alpha }_0}
\def\a1{ {\alpha }_1}

\def\l{\lambda}
\def\o{\omega}

\def\nFGM0{{\nu }_{F,G,M_0}}

\def\lc{{\rm lc}}
\def\nFN0{{\nu}_{F,N_0}}

\def\sm{{\sigma}^m}

\def\sm1{{\sigma}^{-1}}

\def\smtp1{{\sigma}^{-t+1}}

\def\o{\omega }
\def\S1{S^{-1}}

\def\Xpm1{X^{\pm 1}_1}

\def\sPM1{{\sigma }^{\pm 1}}
\def\sMP1{{\sigma }^{\mp 1 }}

\def\b{\beta}
\def\d{\delta}

\def\di{{\rm d.ind}}

\def\L{\Lambda}
\def\O{\Omega}

\def\G{\Gamma}

\def\CA{{\cal A}}

\def\CD{{\cal D}}

\def\Ytm1{Y^{t-1}}
\def\Yim1{Y^{i-1}}

\def\CL{{\cal L}}
\def\CM{{\cal M}}
\def\CN{{\cal N}}
\def\CS{{\cal S}}
\def\CF{{\cal F}}

\def\ass{{\rm ass}}

\def\CZ{{\cal Z}}

\def\Aut{{\rm Aut}}

\def\bA{\overline{A}}
\def\Der{{\rm Der }}
\def\ad{{\rm ad }}
\def\dim{{\rm dim }}
\def\char{{\rm char }}

\def\bi{\bar{i}}

\def\ker{ {\rm ker } }

\def\gr{ {\rm gr} }

\def\SL2Z{ {\rm SL}_2({\bf Z}) }

\def\CZ{ {\cal Z}}

\def\CL{{\cal L}}

\def\Gp1{ G^{1 , 1 } }
\def\P11{ P^{-1 , 1 } }
\def\Pp1{ P^{1 , 1 } }

\def\lcm{{\rm lcm }}

\def\CV{{\cal V}}

\def\nCLsr{{}^\nu\kern-2pt {\cal L}^{\sigma , \rho  }}
\def\nP{{}^\nu \kern-2pt P}
\def\nL{{}^\nu\kern-2pt L}
\def\nLL{{}^\nu\kern-2pt \Lambda}
\def\nPsr{{}^\nu\kern-2pt P^{\sigma , \rho  }}
\def\nLsr{{}^\nu\kern-2pt L^{\sigma , \rho  }}
\def\nuCL{{}^\nu\kern-2pt  {\cal L}}
\def\nCLsr{{}^\nu\kern-2pt {\cal L}^{\sigma , \rho  }}
\def\nCL1m{{}^\nu\kern-2pt {\cal L}^{-1 , 1  }}
\def\x1nu{x^\frac{1}{\nu}}
\def\xm1nu{x^{-\frac{1}{\nu}}}

\def\CB{{\cal B}}
\def\SL2Z{{\rm SL}_2( \mathbb{Z})}
\def\hSL2Z{\widehat{{\rm SL}}_2( \mathbb{Z})}

\def\bwnm1{ \overline{w}_n^{-1}}

\def\twnm1{ \widetilde{w}_n^{-1}}

\def\im{{\rm im}}

\def\o0{\overline{0}}
\def\o1{\overline{1}}

\def\o{\omega}

\def\CB{{\cal B}}
\def\SL2Z{{\rm SL}_2( \mathbb{Z})}
\def\hSL2Z{\widehat{{\rm SL}}_2( \mathbb{Z})}

\def\bwnm1{ \overline{w}_n^{-1}}

\def\twnm1{ \widetilde{w}_n^{-1}}

\def\b0{ \overline{0}}

\def\bx{ \overline{x}}

\def\bi{ \overline{i}}

\def\Gp{\mathfrak{p}}

\def\Gm{\mathfrak{m}}
\def\Gn{\mathfrak{n}}

\def\lcm{{\rm lcm}}
\def\ann{{\rm ann}}

\def\CD{{\cal D}}
\def\End{{\rm End}}

\def\gr{{\rm gr}}

\def\a{\alpha}
\def\b{\beta}

\def\CD{{\cal D}}

\def\HS{{\rm HS}}
\def\di!{\frac{\der^i}{i!}}
\def\dik!{\frac{\der^k_i}{k!}}

\def\Dim{{\rm Dim}}
\def\ud{{\rm ud}}
\def\lud{{\rm lud}}
\def\bder{{\overline{\der}}}

\def\derb{\partial^{[\beta]}}
\def\e{\varepsilon}
\def\Lepk{\Lambda_{\varepsilon}^{[p^k]}}
\def\bpk{{\bf p}^k}
\def\CC{{\cal C}}

\def\CV{{\cal V}}
\def\CDPn{{\cal D}(P_n)}
\def\HS{{\rm HS}}
\def\HSK{{\rm HS}_K}

\def\hsK{{\rm hs}_K}
\def\bk{{\bf k}}
\def\CT{{\cal T}}
\def\coker{{\rm coker}}
\def\Max{{\rm Max}}
\def\CU{{\cal U}}

\def\Tk{T_{\bf k}}
\def\Tl{T_{\bf l}}
\def\TkGm{T_{k(\mathfrak{m})}}
\def\Zk{Z_{\bf k}}
\def\dj!{\frac{\der^j}{j!}}
\begin{document}

\author{V. V. \  Bavula 
}

\title{ Dimension, multiplicity, holonomic modules, and an analogue of the inequality of Bernstein
for  rings of differential operators in prime characteristic
}

\maketitle

$${\it Dedicated \;  to  \; Joseph \; Bernstein\; on  \; the\; ocassion\; of\; his\;  60th \; birthday}$$

\begin{abstract}
Let $K$ be an {\em arbitrary} field of characteristic $p>0$ and
$\CD (P_n)$ be the ring of differential operators on a polynomial
algebra $P_n$ in $n$ variables. A long anticipated {\em analogue
of the inequality of Bernstein} is proved for the ring $\CD
(P_n)$. In fact, {\em three different proofs} are given of this
inequality (two of which are essentially characteristic free): the
first one is based on the concept of the {\em filter dimension},
the second - on the concept of a {\em set of holonomic subalgebras
with multiplicity}, and the third  works only for {\em finitely
presented} modules and follows from a description of these modules
(obtained in the paper).
 On the way, analogues  of the concepts of
 (Gelfand-Kirillov) {\em dimension, multiplicity, holonomic modules} are
 found in prime characteristic (giving answers to old questions of finding such analogs). An idea is very
 simple - to find  characteristic free generalizations (and
 proofs) which in characteristic zero give known results and  in
 prime characteristic - generalizations. An analogue of the {\em Quillen's Lemma} is proved for
 simple {\em finitely presented} $\CD (P_n)$-modules. Moreover, for each such module $L$,
  ${\rm End}_{\CD (P_n)}(L)$ is a finite {\em separable field} extension of $K$
  and $\dim_K({\rm End}_{\CD (P_n)}(L))$ is equal to the {\em multiplicity} $e(L)$ of $L$.  In contrast to the
 characteristic zero case where the Geland-Kirillov dimension of a
 nonzero finitely generated  $\CD (P_n)$-module $M$ can be {\em any natural} number
 from the interval $[n,2n]$, in the prime characteristic, the
 (new)
 dimension $\Dim (M)$  can be \underline{{\em any real}} number
 from the interval $[n,2n]$. It is proved that {\em every
 holonomic}
 module has {\em finite length} but in contrast to the characteristic
 zero case it is {\em not true} neither that a nonzero finitely generated module of
 dimension $n$ is holonomic nor that a holonomic module is
 {\em finitely presented}. Some of the surprising  results are $(i)$  {\em each simple finitely
 presented} $\CD (P_n)$-module $M$ is {\em holonomic} having
 the multiplicity which is a {\em natural} number (in characteristic zero rather  the
  {\em opposite} is true, i.e. $\GK (M)=2n-1$, as a rule), $(ii)$
  the dimension $\Dim (M)$ of a nonzero {\em finitely presented}
  $\CD (P_n)$-module $M$ can be \underline{{\em any natural}} number
 from the interval $[n,2n]$, $(iii)$ the multiplicity $e(M)$
 {\em exists} for each {\em finitely presented} $\CD (P_n)$-module $M$
 and $e(M)\in \mathbb{Q}$, the multiplicity $e(M)$ is a {\em
  natural number} if $\Dim (M)=n$, and can be {\em arbitrary
  small} rational number if $\Dim (M)>n$.


$$ {\bf Contents}$$
\begin{enumerate}
\item Introduction. \item Filter dimension of algebras and
modules. \item Dimension of (not necessarily finitely generated or
Noetherian) algebras and dimension of their finitely generated
modules. \item  An analogue of the inequality of Bernstein for the
ring of differential operators $\CD (P_n)$ with polynomial
coefficients. \item
 Description of finitely
presented $\CD (P_n)$-modules, multiplicity and (Hilbert) almost
polynomials. \item Classification of simple finitely presented
$\CD (P_n)$-modules.  \item Classification of tiny simple
 (non-finitely presented) $\CD (P_n)$-modules.
\item Multiplicity of each finitely presented $\CD (P_n)$-module
is a natural number. \item Holonomic sets of subalgebras with
multiplicity, every holonomic $\CD (P_n)$-module has finite
length.
\end{enumerate}

\end{abstract}


\section{Introduction}
 Throughout the paper, $K$ is a
field,  $P_n=K[x_1, \ldots , x_n]$  a polynomial algebra in $n$
variables over the field $K$, a module  means a {\em left} module,
$\t =\t_K$, $\GK $ stands for the Gelfand-Kirillov dimension.

 In characteristic zero, the ring $\CD (P_n)$ of differential
operators on $P_n$ (so-called, the {\em Weyl algebra}) has
pleasant properties: it is a simple finitely generated Noetherian
domain of Gelfand-Kirillov dimension $\GK (\CD (P_n))=2n$ equipped
with a standard filtration such that the associated graded algebra
$\gr\, \CD (P_n)$ is an affine commutative algebra. {\em None} of
these properties, {\em except simplicity}, holds for the ring $\CD
(P_n)$ in prime characteristic. Moreover, in prime characteristic
the ring $\CD (P_n)$ has a lot of nilpotent elements and zero
divisors. This has a
 serious implication that the standard approach of studying $\CD
 (P_n)$-modules via reduction to modules over affine commutative
 algebras simply is {\em not available}.

 Key ingredients of the theory of (algebraic) $\CD $-modules in
 characteristic zero are the {Gelfand-Kirillov dimension,
 multiplicity, Hilbert polynomial, the inequality of Bernstein, and
 holonomic modules}. In prime characteristic, straightforward
 generalizations of these either do not exist or give `wrong'
 answers (as in the case of the Gelfand-Kirillov dimension:  $\GK
 (\CD (P_n))=n$ in prime characteristic rather than $2n$ as it
 `should' be and it is in characteristic zero).

 In 70'th and 80'th, for rings of differential operators in prime characteristic
  natural questions were posed (see, for example, questions 1-4 in
 \cite{Smith85LNM}) [some of them are still open] that can be
 summarized as {\em to find generalizations of the mentioned
 concepts and results (that results in `good theory'} expectation
 of which was/is high, see, the remark of ${\rm Bj\ddot{o}rk}$ in
 \cite{Smith85LNM}). One of the question in \cite{Smith85LNM} is
 {\em to give a definition of holonomic module in prime
 characteristic}. In characteristic zero, holonomic modules have
 remarkable {\em homological} properties based on which
   Mebkhout and Narvaez-Macarro \cite{Meb-Nar-MacLNM90} gave a
 definition of holonomic module. Another approach (based on the
 {\em Cartier Lemma}) was taken by Bogvad \cite{BogvadJA95} who defined, so-called,
 {\em filtration holonomic} modules. This one is more close to the
 original idea of holonomicity in characteristic zero.
 Note that the two mentioned concepts of holonomicity in prime
 characteristic appeared {\em before} analogues of the
 Gelfand-Kirillov dimension and the inequality of Bernstein have
 been found.

In the present paper, analogues of the Gelfand-Kirillov dimension,
multiplicity, the inequality of Bernstein, and holonomic modules
are found in prime characteristic based on a simple idea - to find
characteristic free generalizations (and proofs) which in
characteristic zero give known concepts (and proofs) and in prime
characteristic - generalizations.

{\bf Filtrations of standard type and the dimension $\Dim $}. A
part of the success story in studying various finitely generated
(Noetherian) algebras is the  class of finite dimensional
filtrations that are equivalent to standard filtrations (a
standard filtration is determined in the obvious way by a finite
set of algebra generators). In general, for an algebra which {\em
not finitely generated} (like $\CD (P_n)$ in prime characteristic)
there is {\em no obvious} choice of finite dimensional filtrations
but for the algebra $\CD (P_n)$ there is  an obvious one -
filtrations that `correspond' to standard filtrations in
characteristic zero, in the present paper  they are called {\em
filtrations of standard type} and an analogue of the {\em
Bernstein filtration} is called the {\em canonical filtration}
$F=\{ F_i\}_{i\geq 0}$ on $\CD (P_n)$ and
$\dim_K(F_i)={i+2n\choose 2n}=\frac{1}{(2n)!}i^{2n}+\cdots $,
$i\geq 0$. Now, a finitely generated $\CD (P_n)$-module $M=\CD
(P_n)M_0$ $(\dim_K(M_0)<\infty)$ is equipped with the  {\em
filtration of standard type} $\{ M_i:=F_iM_0\}$ and one can define
the  {\em dimension} of $M$: $\Dim (M):=\g (i\mapsto \dim_K(M_i))$
where $\g$ denotes the `growth' of function. In particular, $\Dim
(\CD (P_n))=2n$.

{\bf An analogue of the inequality of Bernstein}.

\begin{theorem}\label{Berin}
({\bf The inequality of Bernstein}, \cite{Ber72}) Let $K$ a field
of characteristic zero. Then $\GK (M)\geq n$ for all nonzero
finitely generated $\CD (P_n)$-modules $M$.
\end{theorem}

An analogue of this inequality exists for an arbitrary simple
finitely generated algebra.

\begin{theorem}\label{FFIaf}
\cite{Bavcafd} Let $A$ be a simple finitely generated algebra.
Then
$$
 \GK (M)\geq \frac{\GK (A)}{\fdim (A)+\max \{ \fdim (A),
1\} }
$$
 for all nonzero finitely generated $A$-modules $M$ where $d(A)$
 is a (left) filter dimension of $A$.
\end{theorem}
In particular, $d(\CD (P_n))=1$ (see \cite{Bavjafd}) and $\GK (\CD
(P_n))=2n$ in characteristic zero, and so $\GK (M)\geq
\frac{2n}{1+1}=n$ (Theorem \ref{Berin}).

In Section \ref{pDimWn}, a generalization of Theorem \ref{FFIaf}
(Theorem \ref{FFI}) is given for a simple (not necessarily
finitely generated or Noetherian) algebra equipped with a finite
dimensional filtration.

\begin{theorem}\label{FFIint}
 Let $A$ be a simple algebra with a finite dimensional filtration $F=\{ A_i\}$.
 Then
$$
 \Dim (M)\geq \frac{\Dim (A)}{\fdim(A)+\max \{ \fdim(A),
1\} }
$$
 for all nonzero finitely generated $A$-modules $M$ where $\fdim
 $ is the filter dimension.
\end{theorem}
Applying this result to the algebra $\CD (P_n)$ in prime
characteristic one obtains an analogue of the inequality of
Bernstein in prime characteristic.

\begin{theorem}\label{BerchpInt}
 Let $K$ a field
of characteristic $p>0$. Then $\Dim (M)\geq n$ for all nonzero
finitely generated $\CD (P_n)$-modules $M$.
\end{theorem}

{\it Proof}. Since $\Dim (\CD (P_n))=2n$ and $d(\CD (P_n))=1$
(Theorem \ref{pfdWeyl=1}), applying Theorem \ref{FFIint} we have
$\Dim (M)\geq \frac{2n}{1+1}=n$. $\Box$

The proof is essentially characteristic free.

In characteristic zero, the Gelfand-Kirillov dimension of a
nonzero finitely generated $\CD (P_n)$-module can be  {\em any
natural} number from the interval $[n, 2n]$.

\begin{theorem}\label{}
 Let $K$ be a field of characteristic $p>0$.
\begin{enumerate}
\item (Theorem \ref{3Jl05}) For each \underline{real} number $d$
from the interval $[n,2n]$ there exists a cyclic $\CD
(P_n)$-module $M$ such that $\Dim (M)=d$. \item (Theorem
\ref{23Jun05}) The dimension $\Dim (N)$ of a nonzero finitely
presented $\CD (P_n)$-module $N$ can be any \underline{natural}
number from the interval $[n,2n]$.
\end{enumerate}
\end{theorem}

{\bf Holonomic modules}. A function $f: \mathbb{N}\ra \mathbb{N}$
has {\em polynomial growth} if there exists a polynomial $p(t)\in
\mathbb{Q}[t]$ such that $f(i)\leq p(i)$ for $i\gg 0$. In
characteristic zero, a nonzero finitely generated $\CD
(P_n)$-module $M$ is {\em holonomic} iff $\GK (M)=n$ iff the
function $i\mapsto \dim_K(M_i)$ has polynomial growth of degree
$n$ (i.e. $ \dim_K(M_i)\leq p(i)$ for $i\gg 0$, and  $\deg_t
(p(t))=n$) for some/any standard filtration $\{ M_i\}$ on $M$.

{\it Definition}. In prime characteristic, a nonzero finitely
generated $\CD (P_n)$-module $M$ is {\em holonomic} iff the
function $i\mapsto \dim_K(M_i)$ has polynomial growth of degree
$n$ (i.e. $ \dim_K(M_i)\leq p(i)$ for $i\gg 0$, and $\deg_t
(p(t))=n$) for some (then any) {\em filtration of standard type}
$\{ M_i\}$ on $M$.

\begin{itemize}
\item (Proposition \ref{Mbknh}) {\em In prime characteristic,
there exists a cyclic non-holonomic non-Noetherian $\CD
(P_n)$-module $M$ with} $\Dim (M)=n$. \item (Theorem \ref{hDPnfl})
{\em In prime characteristic, each holonomic module has finite
length and it doest not exceed its `multiplicity'}.
\end{itemize}
These two results show that even having the  analogue of the
Gelfand-Kirillov dimension and the  analogue of the inequality of
Bernstein the `straightforward' generalization of holonomicity
(namely, $\Dim (M)=n$) simply is {\em not correct}.

{\bf Holonomic sets of subalgebras with multiplicity}. For a {\em
simple} algebra $A$ (which is not necessarily finitely generated
or Noetherian), existence of {\em holonomic set of subalgebras
with multiplicity} is another reason why an analogue of the
inequality of Bernstein holds and why each holonomic $A$-module
has finite length (Theorem \ref{ChBI1}).

\begin{itemize}
\item (Theorems \ref{1c27Jn05} and \ref{27Jn05}). {\em In prime
characteristic, the algebra $\CD (P_n)$ has a holonomic set of
subalgebras with multiplicity $1$ (given explicitly).}
\end{itemize}

{\it Definition}. In prime characteristic, a set $\CC =\{ C_\nu
\}_{\nu \in \CN}$ of subalgebras of the algebra $\CD (P_n)$ is
called a {\em holonomic set of subalgebras with mulitplicity} $e$
if for each nonzero $\CD (P_n)$-module $M$ there exists a nonzero
finite dimensional vector subspace $V$ of $M$ such that
$$ \dim_K(C_{\nu, i}V)\geq \frac{e}{n!}i^n+\cdots, \;\;\; i\gg
0,$$
 for some $\nu \in \CN$ where $\{ C_{\nu , i}:=C_\nu\cap F_i\}$ is the
 induced filtration on the algebra $C_\nu$ from the canonical
 filtration $F=\{ F_i\}$ on the algebra $\CD (P_n)$ and the three dots mean $o(i^n)$, smaller terms.

 {\bf Finitely presented $\CD (P_n)$-modules and multiplicity}.
Briefly, in prime characteristic {\em finitely presented} $\CD
(P_n)$-modules  behave similarly as {\em finitely generated} $\CD
(P_n)$-modules in characteristic zero (Theorem \ref{12Jul05}): for
{\em each finitely presented} $\CD (P_n)$-module $M$, the Poincare
series of it is a {\em rational} function, though its Hilbert
function is not a polynomial but an {\em almost polynomial} degree
of which coincides with the dimension $\Dim (M)$ of $M$ (and it
can be {\em any natural number} from the interval $[n,2n]$, this
gives another proof of an analogue of the inequality of Bernstein
for finitely presented $\CD (P_n)$-modules, Theorem
\ref{23Jun05}), and the multiplicity exits for $M$ (Theorem
\ref{23Jun05}). The differences are  $(i)$ in prime
characteristic,  finitely presented $\CD (P_n)$-modules have
transparent structure and are described by Theorem \ref{23Jun05},
but in characteristic zero the category of finitely generated $\CD
(P_n)$-modules is still a mystery, $(ii)$ for each {\em natural}
number $d$ such that $n<d\leq 2n$, there exists a cyclic finitely
presented $\CD (P_n)$-module $M$ with $\Dim (M)=d$ and with {\em
arbitrary small} multiplicity $e(M)$, Lemma \ref{a22Jul05} (in
characteristic zero, multiplicity is a {\em natural} number),
though the multiplicity of {\em every holonomic finitely
presented} $\CD (P_n)$-module is a {\em natural} number (Theorem
\ref{t12Jul05}), $(iii)$ and what is completely unexpected is that
{\em each simple finitely presented} $\CD (P_n)$-module is {\em
holonomic} (Corollary \ref{1c23Jun05}), and if, in addition, the
field $K$ is {\em algebraically closed} then the multiplicity is
{\em always} $1$ (Corollary \ref{c8Jul05}).

{\bf A classification of simple finitely presented $\CD
(P_n)$-modules and an analogue of Quillen's Lemma}.
  {\em In prime characteristic (see Theorem \ref{8Jul05}),}
  \begin{itemize}
\item {\em A clasification of simple finitely presented $\CD
(P_n)$-modules is obtained.} \item {\em Every simple finitely
presented $\CD (P_n)$-module $M$ is holonomic, and} \item ({\bf An
analogue of Quillen's Lemma}) {\em its endomorphism algebra} ${\rm
End}_{\CD (P_n)}(M)$ {\em is a finite separable field extension of
$K$, and } \item $\dim_K({\rm End}_{\CD (P_n)}(M))=e(M)$, {\em the
multiplcity of $M$, and } \item {\em if, in addition, the field
$K$ is algebraically closed then always $e(M)=1$.}
  \end{itemize}

  {\bf A classification of tiny simple  $\CD
  (P_n)$-modules}. A classification is obtained of the `smallest'
simple $\CD (P_n)$-modules (see Theorems \ref{10Jul05} and
\ref{8Jul05}), they are called {\em tiny} modules. Theorem
\ref{8Jul05} describes the set of tiny {\em finitely presented}
$\CD (P_n)$-modules and Theorem \ref{10Jul05} classifies the set
of tiny {\em non-finitely presented} $\CD (P_n)$-modules. They
turned out to be {\em holonomic} with multiplicities which are
{\em natural} numbers. Briefly, they have the same properties as
simple finitely presented $\CD (P_n)$-modules.

Results of this paper have been generalized for the ring of
differential operators on a {\em smooth irreducible affine
algebraic variety}, \cite{Bav-holmodII}.


\section{Filter dimension of algebras and
modules}\label{pfdim}

The filter dimension is one of the key ingredients in the proof of
an analogue of the inequality of Bernstein in prime
characteristic.

Originally, the filter dimension was defined for any {\em finitely
generated} algebra $A$ and any {\em finitely generated}
$A$-module.  In this section, the concept of filter dimension of
algebras and modules will be extended to a class of not
necessarily finitely generated algebras.

{\bf The concept of growth}. Let $\CF $ be the set of all
functions from the set of natural numbers $ \mathbb{N}=\{ 0, 1,
\ldots \}$ to itself. For each function $f\in \CF $, the
non-negative real number or $\infty $ defined as
$$ \g (f):=\inf \{  r\in \mathbb{R}\, | \, f(i)\leq i^r\; {\rm for
}\; i\gg 0\}$$ is called the  {\bf degree} of $f$. The function
$f$ has {\bf polynomial growth} if $\g (f)<\infty $. Let $f,g,
p\in \CF $, and $p(i)=p^*(i)$ for $i\gg 0$ where $p^*(t)\in
\mathbb{Q}[t]$ (a polynomial algebra with coefficients from the
field  of rational numbers). Then
\begin{eqnarray*}
\g (f+g)\leq \max \{ \g  (f), \g (g)\}, & & \g (fg)\leq \g (f)+ \g
(g),\\
\g (p)=\deg_t(p^*(t)), & & \g (pg)= \g (p)+ \g
(g).\\
\end{eqnarray*}

{\bf The equivalence relation on the class of filtrations}. Let
$A$ be an algebra over an {\em arbitrary} field $K$. Recall that a
{\em filtration} $F=\{ A_i\}_{i\geq 0}$ of the algebra $A$ is an
ascending chain of vector subspaces of $A$:
$$ A_0\subseteq A_1\subseteq  \cdots \subseteq A_i\subseteq \cdots ,
 \;\;\; A=\cup_{i\geq 0}A_i,
\;\; K\subseteq A_0, \;\; A_iA_j\subseteq A_{i+j}, \;\; i,j\geq
0.$$ The filtration $F$ is a {\bf finite dimensional} filtration
(or a {\bf finite} filtration, for short) provided $\dim_K\, F_i
<\infty$ for all $i\geq 0$. Filtrations  $F=\{ A_i\} $ and $G=\{
B_i \}$ on $A$ are called {\em equivalent} $(F\sim G)$ if there
exist natural numbers $a,b,c,d$ such that $a>0$, $c>0$ and
$$A_i\subseteq B_{ai+b}\;\;\; \hbox{and}\;\;\; B_i\subseteq A_{ci+d}\;\;\; \hbox{for}\;
i\gg 0 .$$ The equivalent filtrations  $F=\{ A_i\} $ and $G=\{ B_i
\}$ on $A$ are called {\em strongly  equivalent} if $a=c=1$. A
similar definition exists for filtrations on modules rather than
algebras.

Clearly, this is an equivalence relation on the  class of all
filtrations of the algebra $A$. For a filtration $F$,
 $\widetilde{F}$ denotes the equivalence class of the filtration
 $F$.  If one of the inclusions above  holds, say the first, we write
$F\leq G$.

{\bf The Gelfand-Kirillov dimension}. If  $A=K\langle a_1, \ldots
, a_s\rangle $ is a {\em finitely generated}
 $K$-algebra. The finite dimensional filtration $F=\{ A_i\}$
associated with algebra generators $ a_1, \ldots , a_s$:
$$ A_0:=K\subseteq A_1:=K+\sum_{i=1}^sKa_i\subseteq \cdots \subseteq
A_i:=A_1^i\subseteq \cdots $$ is  called the {\bf standard
filtration} for the algebra $A$. Let $M=AM_0$ be a finitely
generated $A$-module where $M_0$ is a finite dimensional
generating subspace of the $A$-module $M$. The finite dimensional
filtration $\{ M_i:=A_iM_0\}$ is called the {\bf standard
filtration} for the $A$-module $M$. All standard filtrations of an
algebra $A$ (or a finitely generated $A$-module) are equivalent.

{\it Definition}. $\GK (A):=\g (i\mapsto \dim_K(A_i))$ and $ \GK
(M):=\g (i\mapsto \dim_K(M_i))$ are called the {\bf
Gelfand-Kirillov} dimensions of the algebra $A$ and the $A$-module
$M$ respectively.

It is easy to prove that the Gelfand-Kirillov dimension  of the
algebra (resp. the module)  does not depend on the choice of the
standard filtration of the algebra (resp. and the choice of the
generating subspace of the module) see \cite{KL} for details. This
is a direct consequence of the fact that all the standard
filtrations are equivalent.

 The results we are going
to generalize first were proved for {\em finitely generated}
algebras (and their {\em finitely generated} modules)  equipped
with standard filtrations. Here we extend results to arbitrary
filtrations (mainly finite dimensional) on a not necessarily
finitely generated algebras. The results do not depend on a
filtration inside its equivalence class, but, in general, they do
depend on the equivalence class. The choice of the equivalence
class depends on a concrete class of algebras.

Our main motivation is an equivalence class of finite dimensional
  filtrations on a ring of differential operators $\CD (A)$ in
prime characteristic that in characteristic zero coincide with the
class of all the {\em standard} filtrations on the algebra $\CD
(A)$.

{\bf The return functions and the (left) filter dimension}. Let
$A$ be a filtered algebra with a filtration $F=\{ A_i\}$, and let
$M=AM_0$ be a finitely generated $A$-module with a {\em finite
dimensional} generating subspace $M_0$. Then $M=\cup_{i\geq 0}M_i$
is a filtered $A$-module with the filtration $\{ M_i:=A_iM_0\}$
which obviously does depend on the filtration $F$ and a generating
subspace $M_0$. When one fixes the filtration $F$ then distinct
finite dimensional subspaces of the $A$-module $M$ give {\em
equivalent} filtrations on the module $M$.

The next definition appeared in  \cite{Bavcafd} in case of
standard filtrations.

 {\it Definition}. The function $\nFM0 :\mathbb{N}\ra
\mathbb{N}\cup \{ \infty \}$,
$$ \nFM0 (i):=\min \{ j\in\mathbb{N}\cup \{ \infty \}: \;
A_jM_{i,gen}\supseteq M_0\;\; {\rm for \; all}\;\; M_{i,gen} \}$$
is called the {\bf return function} of the $A$-module $M$
associated with the filtration $F=\{ A_i\}$ of the algebra $A$ and
the generating subspace $M_0$ of the $A$-module $M$ where
$M_{i,gen}$ runs through {\em all finite dimensional generating
subspaces} for the $A$-module $M$ such that $M_{i,gen}\subseteq
M_i$.

Suppose, in addition, that the algebra $A$ is a {\em simple}
algebra.   The {\bf return function} $\nu_F \in \CF $ and the {\bf
left return function} $\l_F\in \CF $ for the algebra $A$ with
respect to the filtration $F:= \{ A_i\}$ for the algebra $A$ are
defined by the rules:
\begin{eqnarray*}
\nu_F(i)&:=& \min \{ j\in \mathbb{N}\cup \{ \infty \} \, | \,\,
1\in A_jaA_j\;\, {\rm
for \; all}\; \,  0\neq a\in A_i\},\\
\l_F(i)&:=& \min \{ j\in \mathbb{N}\cup \{ \infty \} \, | \, 1\in
AaA_j\; \; {\rm for \; all}\; \; 0\neq a\in A_i\},
\end{eqnarray*}
where $A_jaA_j$ is the vector subspace of the algebra $A$ spanned
over the field $K$ by the elements $xay$ for all $x,y\in A_j$; and
$AaA_j$ is the left ideal of the algebra $A$ generated by the set
$aA_j$. Similarly, the {\bf unit return function} $\nu_F^u \in \CF
$ and the {\bf left unit return function} $\l_F^u\in \CF $  are
defined (where $U=U(A)$ is the {\em group of units}, i.e.
invertible elements of $A$):
\begin{eqnarray*}
\nu_F^u(i)&:=& \min \{ j\in \mathbb{N}\cup \{ \infty \} \, | \,\,
U(A) \cap A_jaA_j\neq \emptyset \;\, {\rm
for \; all}\; \,  0\neq a\in A_i\},\\
\l_F^u(i) &:=& \min \{ j\in \mathbb{N}\cup \{ \infty \} \, |
\,U(A) \cap AaA_j\neq \emptyset \; \; {\rm for \; all}\; \; 0\neq
a\in A_i\}.
\end{eqnarray*}
Clearly, 
\begin{equation}\label{4rf}
\l^u_F(i)\leq \l_F(i)\leq \nu_F(i)\;\; {\rm and}\;\; \l^u_F(i)\leq
\nu^u_F(i)\leq \nu_F(i)\;\; \; {\rm for \; all}\;\;\; i\geq 0.
\end{equation}

The next result shows that under a mild restriction the four
return functions take only {\em finite} values. In general, there
is no reason to believe that values of the return functions are
always finite, but for central simple algebras equipped with an
arbitrary finite dimensional filtration this is always the case
(see the next lemma). Recall that the centre of a simple algebra
is a field.
\begin{lemma}\label{fdfin}
Let $A$ be a simple  algebra equipped with a finite dimensional
filtration $F=\{ A_i\}$ such that the  centre $Z(A)$ of the
algebra $A$ is an algebraic field extension of $K$. Then the four
return functions take  finite values.
\end{lemma}

{\it Proof}. In a view of (\ref{4rf}), it suffices to prove the
lemma for the return function $\nu_F$, that is $\nu_F(i)<\infty$
for all $i\geq 0$.

 The centre $Z=Z(A)$ of the simple algebra $A$ is a field that
 contains $K$. Let $\{ \o_j\, | \, j\in J\}$ be a $K$-basis for the $K$-vector
space $Z$. Since $\dim_K(A_i)<\infty $, one can find  finitely
many $Z$-linearly independent elements, say $a_1,\ldots , a_s$, of
$A_i$ such that $A_i\subseteq Za_1+\cdots +Za_s$. Next, one can
find a finite subset, say $J'$, of $J$ such that $A_i\subseteq
Va_1+\cdots +Va_s$ where $V=\sum_{j\in J'}K\o_j$. The field $K'$
generated over $K$ by the elements $\o_j$, $j\in J'$, is a finite
field extension of $K$ (i.e. $\dim_K(K')<\infty $) since $Z/K$ is
algebraic, hence $K'\subseteq A_n$ for some $n\geq 0$. Clearly,
$A_i\subseteq K'a_1+\cdots +K'a_s$.

The $A$-bimodule ${}_AA_A$ is simple with ring of endomorphisms
${\rm End}({}_AA_A)\simeq Z$. By the {\em Density} Theorem,
\cite{Pierceb}, 12.2, for each integer $1\leq j \leq s$, there
exist elements of the algebra $A$, say $x_1^j, \ldots , x_m^j,
y_1^j, \ldots , y_m^j$, $m=m(j)$, such that for all $1\leq l\leq
s$
$$ \sum_{k=1}^m x_k^ja_ly_k^j=\delta_{j,l}, \;\; {\rm the \; Kronecker
\; delta}.$$ Let us fix a natural number, say $d=d_i$, such that
$A_d$ contains all the elements $x_k^j$, $y_k^j$, and the field
$K'$. We claim that $\nu_F(i)\leq 2d$. Let $0\neq a\in A_i$. Then
$a=\l_1a_1+\cdots +\l_sa_s$ for some $\l_i\in K'$. There exists
$\l_j\neq 0$. Then $\sum_{k=1}^m\l_j^{-1}x_k^ja_jy^j_k=1$, and $
 \l_j^{-1}x_k^j, y^j_k\in A_{2d}$. This proves the claim and the lemma.
  $\Box $

{\it Remark}.  If the field $K$ is {\em uncountable} then
automatically the centre $Z(A)$ of a simple finitely generated
algebra $A$ is algebraic over $K$  (since $A$ has a countable
$K$-basis and the rational function field $K(x)$ has uncountable
basis over $K$ since elements $\frac{1}{x+\l }$, $\l \in K$, are
$K$-linearly independent).

In what follows we will assume that the four return functions {\bf
do not take infinite value}.

\begin{lemma}      
 Let $A$ be an algebra equipped with two equivalent
 filtrations $F=\{ A_i\}$ and $G =\{ B_i\}.$
\begin{enumerate}
\item  Let $M$ be a finitely generated $A$-module. Then $\g (\nFM0
)=\g (\nGN0 )$ for any finite dimensional generating subspaces
$M_0$ and $N_0$ of the $A$-module $M$. \item  If, in addition, $A$
is a simple algebra  then  $\g (\nF )=\g (\nG )$,  $\g (\l_F)=\g
(\l_G)$, and $\g (\nu_F)=\g (\nu_{F\otimes F^o, K})$ where
$\nu_{F\otimes F^o, K}$ is the return function of the $A\otimes
A^o$-module $A$ and $A^o$ is the opposite algebra to $A$. \item
If, in addition, $A$ is a simple algebra then $\g (\nu_F^u )=\g
(\nu_G^u )$ and $\g (\l_F^u)=\g (\l_G^u)$.
\end{enumerate}
\end{lemma}

{\it Proof}. 1. The module $M$ has two filtrations $\{M_i=A_iM_0
\}$ and $\{N_i=B_iN_0 \} $. Let $\nu =\nFM0 $    and $\mu =\nGN0
$.

First, we consider two special cases, then the general case will
follow easily from these two. Suppose first that $F=G$. Choose a
natural number $s$ such that $M_0\subseteq N_s$ and $N_0\subseteq
M_s$, then $N_i\subseteq M_{i+s}$ and $M_i\subseteq N_{i+s}$ for
all $i\geq 0$. Let $N_{i,gen}$ be any generating subspace for the
$A$-module $M$ such that $N_{i,gen}\subseteq N_i$.  Since
$M_0\subseteq A_{\nu (i+s)}N_{i,gen}$ for all $i\geq 0$ and $N_0
\subseteq A_sM_0,$ we have $N_0\subseteq A_{\nu
(i+s)+s}N_{i,gen}$,  hence $\mu (i)\leq \nu (i+s)+s$ and finally
$\g (\mu )\leq \g (\nu ).$ By symmetry, the opposite inequality is
true and  so $\g (\mu )=\g (\nu ) $.

Now, suppose that $M_0=N_0$.  Since  $F\sim G$ one can choose
natural numbers $a,b,c,d$ such that $a>0$, $c>0$ and
$$A_i\subseteq B_{ai+b}\;\;\; \hbox{and}\;\;\; B_i\subseteq A_{ci+d}\;\;\;
\hbox{for}\;\;\;  i\gg 0 .$$

 Then $N_i=B_iN_0\subseteq A_{ci+d}M_0=M_{ci+d}$ for all $i\geq
 0$, hence $N_0=M_0\subseteq A_{\nu (ci+d)}N_{i,gen}\subseteq B_{a\nu
(ci+d)+b} N_{i,gen}$, therefore $\mu (i)\leq a\nu (ci+d)+b$ for
all $i\geq 0$,  hence $\g (\mu )\leq \g (\nu ) $. By symmetry, we
get
 the  opposite inequality which implies  $\g (\mu )=\g (\nu )$.

In the general case, $\g (\nFM0 )=\g (\nFN0 )=\g (\nGN0 ) $.

2. The algebra $A$ is simple, equivalently, it is  a simple (left)
$A\t A^o$-module where $A^o$  is the {\em opposite} algebra to
$A$. The opposite algebra has the  filtration $F^o=\{ A_i^o\}$.
 The tensor product of algebras $A\t A^o$, so-called,  the {\em
enveloping algebra} of $A$, has the  filtration $F\t F^o= \{
C_n\}$ which is the tensor product of the  filtrations $F$ and
$F^o$, that is, $C_n=\sum \{ A_i\t A_j^o,i+j\leq n \} $. Let ${\nu
}_{ F\t F^o, K }$ be the return function of the $A\t A^o$-module
$A$ associated with the filtration $F\t F^o$ and the generating
subspace $K$. Then
$$\nF (i)\leq {\nu }_{ F\t F^o,K }(i)\leq 2\nF (i) \; \hbox{for all}\;  i\geq 0,$$
 and so
\begin{equation}\label{gnuFK}
\g (\nF )=\g ({\nu }_{F\t F^o,K}),
\end{equation}
 and, by the first statement, we
have  $\g (\nu_F)= \g ({\nu }_{F\t F^o,K})=\g ({\nu }_{G\t
G^o,K})=\g (\nu_G)$, as required. Using a similar argument as in
the proof of the first statement one can proof that $\g (\l_F)=\g
(\l_G)$. We leave this as an exercise.

$3$. Let $F$, $G$, $a,b,c,d$ be as above. Let $U=U(A)$ be the
group of units of the algebra $A$, and let $\l:=\nu_F^u$ and $ \mu
:=\nu^u_G$ (resp. $\l :=\l^u_F$ and $\mu :=\l_G^u$). We prove two
cases simultaneously. Let $x$ be a nonzero element of $A_i$. Then
$0\neq x\in B_{ai+b}$ and
\begin{eqnarray*}
\emptyset &\neq & U\cap B_{\mu (ai+b)}xB_{\mu (ai+b)}\subseteq
U\cap A_{c\mu (ai+b)+d}xA_{c\mu (ai+b)+d},\\
\emptyset &\neq & U\cap AxB_{\mu (ai+b)}\subseteq  U\cap A
xA_{c\mu (ai+b)+d} \;\;\; {\rm respectively}.
\end{eqnarray*}
In both cases, $\g (\l )\leq \g (c\mu (ai+b)+d)\leq \g (\mu )$. By
symmetry, the inverse inequality is also true, and so $\g (\l ) =
\g (\mu )$. $\Box $

{\it Definition}.  $\fd (M)=\g (\nFM0 )$ is the {\bf filter
dimension} of the $A$-module $M$, and $\fd (A):= \fd ({}_{A\t
A^o}A)$ is the {\bf filter dimension} of the algebra $A$. If, in
addition, the algebra $A$ is simple, then $\fd (A)=\g (\nu_F)$,
$\lfd (A):=\g (\l_F)$ is called the {\bf left filter dimension} of
the algebra $A$, $\ud (A)=\g (\nu_F^u)$ is called the {\bf unit
dimension} of $A$, and  $\lud (A):=\g (\l_F^u)$ is called the {\bf
left unit dimension} of the algebra $A$.

By the previous lemma, the definitions make sense provided an
equivalence class of filtrations is fixed. We will always assume
that we have fixed such a class. A particular choice of an
equivalence class of filtrations depends on a class of algebras we
study. For finitely generated algebras such an equivalence class
as a rule is the equivalence class of all standard filtrations,
but for algebras that are not finitely generated there is {\em no
obvious choice} of an equivalence class of filtrations.

For standard filtrations the concept of (left) filter dimension
first appeared in  \cite{Bavcafd}.

By (\ref{4rf}),
\begin{equation}\label{4rf1}
\lud (A)\leq \lfd (A)\leq \fd (A)\;\;\; {\rm and}\;\;\; \lud
(A)\leq\ud (A)\leq \fd (A).
\end{equation}


\section{Dimension of (not necessarily finitely generated or Noetherian)
algebras and dimension of their finitely generated
modules}\label{pDimWn}

 Theorem \ref{FFI} is the main result of this section, it is a kind of the inequality
 of Bernstein but for an arbitrary simple algebra (not necessarily finitely generated)
  equipped with a finite dimensional filtration.  In this
section, let $A$ be an algebra over an {\em arbitrary} field $K$
with a {\em finite dimensional} filtration $F=\{ A_i\}$. Let
$M=AM_0$ be a finitely generated  $A$-module with a finite
dimensional generating subspace $M_0$. Then $M$ has a finite
dimensional filtration $\{ M_i:=A_iM_0\}$. Suppose that $G=\{
B_i\}$ is a finite dimensional filtration on $A$ equivalent to the
filtration $F$ and let $N_0$ be another finite dimensional
generating subspace for the $A$-module $M$. Then the $A$-module
$M$ has a second finite dimensional filtration $\{ N_i:=
A_iN_0\}$. It follows easily that $\g (\dim_K \, A_i)=\g (\dim_K
\, B_i)$ and $\g (\dim_K \, M_i)=\g (\dim_K \, N_i)$.

{\em Definition}. The {\bf dimension} $\Dim \, A$ of the algebra
$A$ and the {\bf dimension} $\Dim \, M$ of the finitely generated
$A$-module $M$ are the numbers $\g (\dim_K \, A_i)$ and $\g
(\dim_K \, M_i)$ respectively.

 So, the dimension $\Dim \, A$ of the
algebra $A$ is an invariant of the algebra $A$  {\em and} the
equivalence class of the filtration $F$. The same is true about
the dimension $\Dim \, M$ of the $A$-module $M$.

If $A$ is a finitely generated algebra and $\{ A_i\}$ is a {\em
standard} filtration then $\Dim (A)=\GK (A)$ and $\Dim (M)=\GK
(M)$.

 In this paper,  $\fdim (A)$ stands for any of the dimensions $\fd (A)$,
 $\lfd (A)$, $\ud (A)$ or $\lud (A)$ of an  algebra $A$ (i.e. $\fdim =\fd , \lfd, \ud, \lud $).

 The four dimensions  appear
  naturally when one tries to find a
{\em lower} bound for the holonomic number (Theorem \ref{FFI}).

 The next theorem is a generalization of the {\bf
 inequality of Bernstein} (Theorem \ref{Berin}) to the class of
simple algebras. This result was first appeared in \cite{Bavcafd,
bie98} in the case of {\em simple finitely generated} algebras
with respect to the class of {\em standard} filtrations and for
$d=\fd, \lfd$.

\begin{theorem}\label{FFI}
 Let $A$ be a simple algebra with a finite dimensional filtration $F=\{ A_i\}$.
 Then
$$
 \Dim (M)\geq \frac{\Dim (A)}{\fdim(A)+\max \{ \fdim(A),
1\} }
$$
 for all nonzero finitely generated $A$-modules $M$ where $\fdim
 =\fd , \lfd, \ud ,  \lud $.
\end{theorem}

{\it Proof}. In a view of (\ref{4rf1}), it suffices to prove the
theorem for $d=\lud$. Let $\l =\l_F^u$ be the left unit return
function associated with the finite dimensional filtration $F$ of
the algebra $A$ and let $0\neq a \in A_i$.  It follows from the
inclusion
$$AaM_{\l (i)}=AaA_{\l (i)}M_0\supseteq (U(A)\cap AaA_{\l (i)})M_0\neq 0$$
that the linear map
$$A_i\ra  \Hom_K (M_{\l (i)},M_{\l (i)+i}), a\mapsto
(m\mapsto  am),$$ is injective, and  so  dim $A_i\,\le $ dim
$M_{\l (i)}$ dim $M_{\l  (i)+i}$.
 Using the above elementary properties of the degree  (see also
\cite{MR}, 8.1.7), we have
\begin{eqnarray*}
\Dim (A)& =& \gamma (\dim \, A_i)\le \gamma (\dim \, M_{\l
(i)})+\gamma
(\dim \, M_{\l (i)+i})\\
&\le &  \gamma (\dim \, M_i )\gamma (\l )+\gamma (\dim \, M_i
)\max \{
\gamma (\l ),1\}\\
&=& \Dim (M)(\lud A+\max \{\lud A, 1\})\\
&\leq & \Dim (M)(\lud A+\max \{\lud A, 1\}). \;\; \Box
\end{eqnarray*}

The {\em inequality of Bernstein} says that $\GK (M)\geq n$ {\em
for any nonzero finitely generated module $M$ over a ring of
differential operators $\CD (X)$ on a smooth irreducible affine
algebraic variety $X$ of dimension $n=\dim \, X$ over a field of
characteristic zero}. Since $\GK (\CD (X))=2n$ and $ \fd (\CD
(X))=\lfd (\CD (X))=1$ (\cite{Bavjafd}), by Theorem \ref{FFI}, we
have a `short' proof of the inequality of Bernstein:
$$ \GK (M)\geq \frac{2n}{1+1}=n.$$

{\it Definition}. $h_A:=\inf \{ \Dim (M)\, | \,  M $  is a nonzero
finitely  generated $A$-module$\}$ is called the {\bf holonomic
number} for the algebra $A$ (with respect to the equivalence class
$\widetilde{F}$ of the finite dimensional filtration $F$).

The result above gives a {\em lower bound} for the holonomic
number of the  simple  algebra $A$:
$$ h_A\geq  \frac{\GK (A)}{\fdim(A)+\max \{ \fdim(A),
1\} }.
$$
\begin{theorem}\label{IIIin}
Let $A$ and $\widetilde{F}$ be as above. Then
$$\Dim (M)\leq \Dim(A)\,\fd (M)$$
for any simple $A$-module $M$.
\end{theorem}

{\it Proof}. Let $\nu =\nu_{F, Km} $ be the return function of the
module  $M$ associated with the  finite dimensional filtration
$F=\{A_i\}$ of the algebra $A$ and a fixed nonzero element  $m\in
M$. Let $\pi :M\ra K$ be a non-zero linear map satisfying $\pi
(m)=1$. Then, for any $i\geq 0$ and any $0\neq u\in M_i:=A_im$:
$1=\pi (m)\in \pi (A_{\nu (i)}u)$, and so the linear map
$$M_i\ra \Hom_K (A_{\nu (i)},K),\,\,u\mapsto (a\mapsto \pi (au)),$$
is an {\em injective} map hence dim $M_i\leq $ dim $A_{\nu (i)}$
and finally  $\Dim (M)\leq \Dim (A)\, \fd (M)$.  $\Box$

\begin{corollary}
Let $A$ be a simple algebra with $\Dim (A)>0$. Then
 $$\fd (A)\geq \frac{1}{2}.$$
\end{corollary}

{\it Proof}.  Clearly,
 $\Dim (A\t A^o)\leq \Dim (A)+\Dim (A^o)=2\Dim (A)$.
Applying Theorem  \ref{IIIin} to  the  simple $A\t A^o$-module
$M=A$ we finish the proof
$$ \Dim (A)=\Dim ({}_{A\t A^o}A)\leq \Dim (A\t A^o)\fd
({}_{A\t A^o}A)\leq 2\Dim (A) \fd (A)$$ hence $\fd (A)\geq
\frac{1}{2}$. $\Box $

\begin{corollary}
Let $A$ be a simple algebra with $\Dim (A)>0$. Then
 $$\fd (M)\geq \frac{1}{\fd (A)+\max \{ \fd (A), 1\} }$$
  for all simple $A$-modules $M$.
\end{corollary}

{\it Proof}. Applying Theorem \ref{FFI} and Theorem \ref{IIIin},
we have the result
$$ \fd (M)\geq \frac{\Dim (M)}{\Dim (A)}\geq
\frac{\Dim (A)}{\Dim (A) (\fd (A)+\max \{ \fd (A), 1\}
)}=\frac{1}{\fd (A)+\max \{ \fd (A), 1\} }. \;\; \Box
$$

In general, it is difficult to find the exact value for the filter
dimension but for the ring of differential operators $\CD (P_n)$
with polynomial coefficients $P_n=K[x_1, \ldots , x_n]$ over  a
field $K$ of characteristic $p>0$ it is easy and one can find it
directly (Theorem \ref{pfdWeyl=1}).


\section{An analogue of the inequality of Bernstein
for the  ring of differential operators $\CD (P_n)$
 with polynomial coefficients}\label{hBIDPn}

In this section,  $K$ is an {\em arbitrary} field of
characteristic $p>0$, $P_n:=K[x_1, \ldots , x_n]$ is a polynomial
algebra, $\CD =\CD (P_n)$ is the ring of differential operators on
$P_n$. In this section, the concepts of {\em filtration of
standard type} and of {\em holonomic module} are introduced, it is
proved that the filter dimension of the ring $\CD (P_n)$ is $1$
(Theorem \ref{pfdWeyl=1}) and an analogue of the inequality of
Bernstein is established (Theorem \ref{BIDPnp}) in prime
characteristic. We start with recalling some facts and properties
of higher derivations (Hasse-Schmidt derivations) which will be
used freely in the paper.

{\bf Higher derivations}. Let us recall basic facts about higher
derivations. For more detail the reader is referred to \cite{Ma},
Sec. 27.

A sequence $\d =(1:={\rm id}_A, \d_1, \d_2, \ldots )$ of
$K$-linear maps from a commutative $K$-algebra $A$ to itself
(where ${\rm id}_A$ is the identity map on $A$) is called a {\em
higher derivation} (or a {\em Hasse-Schmidt derivation}) over $K$
from $A$ to $A$ if, for each $k\geq 0$, 
\begin{equation}\label{dkxyhi}
\d_k(xy)=\sum_{i+j=k}\d_i(x)\d_j(y)\; \;\; {\rm for\; all}\;\;
x,y\in A.
\end{equation}
Clearly, $\d_1\in \Der_K(A)$. These conditions are equivalent to
saying that the map $e:A\ra A[[t]]$, $x\mapsto \sum_{i\geq
0}\d_i(x)t^i$, is a $K$-algebra homomorphism where $A[[t]]$ is a
ring of power series with coefficients from $A$, or equivalently,
that the map
$$e:A[[t]]\ra A[[t]], \;\; t\mapsto t, \;\; x\mapsto \sum_{i\geq
0}\d_i(x)t^i \;\;\; (x\in A),$$ is a $K[[t]]$-algebra
homomorphism. Clearly, $e$ is a $K[[t]]$-algebra automorphism of
$A[[t]]$, and vice versa (any automorphism $e\in
\Aut_{K[[t]]}(A[[t]])$ of the type $e(a)=a+\sum_{i\geq
1}\d_i(a)t^i$ yields a higher derivation $(\d_i)$ where $a\in A$).

The set $\HS_K(A)$ of all higher $K$-derivations from $A$ to $A$
is a subgroup of the group $\Aut_K(A[[t]])$ of all $K$-algebra
automorphisms of $A[[t]]$. It follows immediately that a higher
derivation has a {\em unique} extension to a localization $\S1 A$
of the algebra $A$ at a multiplicative subset $S$ of $A$.

Let $\CD (A)$ be the {\em ring of differential operators} on the
algebra $A$ and let $\{ \CD (A)_i\}_{i\geq 0}$ be its {\em order
filtration}. Recall that $\CD (A)=\cup_{i\geq 0}\CD (A)_i\subseteq
{\rm End}_K(A)$, $\CD (A)_0:= {\rm End}_A(A)\simeq A$, and
$$ \CD (A)_i:= \{ f\in {\rm End}_K(A): \, fx-xf\in \CD
(A)_{i-1}\;\; {\rm for \; all}\;\; x\in A\}, \;\; i\geq 1.$$

Let  $\d = (\d_i)\in \HS_K(A)$. By (\ref{dkxyhi}),
\begin{equation}\label{diCDAi}
\d_i\in \CD (A)_i, \;\; i\geq 0,
\end{equation}
since $\d_ix-x\d_i= \sum_{j=0}^{i-1}\d_{i-j}(x)\d_j$ for all $x\in
A$ and the result follows by induction on $i$. For each $i\geq 1$
and $x\in A$,
$$\sum_{j\geq 0}\d_j(x^{p^i})t^j=e(x^{p^i})=e(x)^{p^i}=\sum_{k\geq
0}\d_k(x)^{p^i}t^{kp^i},$$ and so
$\d_{kp^i}(x^{p^i})=\d_k(x)^{p^i}$ for all $i,k\geq 0$ and $x\in
A$; and $\d_j(x^{p^i})=0$ for all $j$ such that $p^i\not| j$
($p^i$ does not divide $j$). In particular,
$\d_l(KA^{p^{i+1}})=0$ for all $i\geq 1$ and $0<l<p^{i+1}$.

{\bf The higher derivations $(1, \frac{\der}{1!},
\frac{\der^2}{2!}, \ldots )\in \HS_K(K[x])$ where $\der
=\frac{d}{dx}$}. Given a polynomial algebra $K[x]$ in a single
variable $x$ over $K$, the $K$-algebra homomorphism $K[x]\ra
K[x][[t]]$, $f(x)\mapsto f(x+t)=\sum_{i\geq 0}\di! (f)t^i$,
determines the higher derivation
 $(1, \frac{\der}{1!},
\frac{\der^2}{2!}, \ldots )\in \HS_K(K[x])$. If ${\rm char} (K)=0$
then $\dj!$ means $(j!)^{-1} \der^j$, but if ${\rm char} (K)=p>0$
then 
\begin{equation}\label{dxj=bin}
\dj! (x^i)={i\choose j}x^{i-j}
\end{equation}
where ${i\choose j}$ is the binomial in characteristic $p$:
${i\choose j}=0$ if $j>i$, and for $j\leq i$, let $j=\sum j_kp^k$,
$0\leq j_k<p$,  and $i=\sum i_kp^k$, $0\leq i_k<p$. Then
\begin{equation}\label{binijchp}
{i\choose j}=\prod_k {i_k\choose j_k}
\end{equation}
where ${i_k\choose j_k}=0$ if $j_k>i_k$, and ${i_k\choose
j_k}=\frac{i_k!}{j_k!(i_k-j_k)!}$ if $j_k\leq i_k$. The formulas
(\ref{dxj=bin}) and (\ref{binijchp}) are obvious when one looks at
the following product:
$$
(x+t)^i=\prod_k(x^{p^k}+t^{p^k})^{i_k}=\prod_k\sum_{l_k=0}^{i_k}{i_k\choose
l_k}x^{(i_k-l_k)p^k}t^{l_kp^k}= \sum \prod_k {i_k\choose
i_k-l_k}x^{\sum_\nu (i_\nu -l_\nu )p^\nu }t^{\sum_\mu l_\mu p^\mu
}$$ where the sum runs through all $l_0, l_1, \ldots $ that
satisfy $0\leq l_0\leq i_0, 0\leq l_1\leq i_1, \ldots $. The
binomials in characteristic $p>0$ has a remarkable property - the
{\em translation invariance} (with respect to the $p$-{\em adic
scale}): 
\begin{equation}\label{trinvbi}
{p^ki\choose p^kj}={i\choose j}, \;\;\; k\geq 0.
\end{equation}
This follows directly from (\ref{binijchp}). By (\ref{binijchp}),
${i\choose j}\neq 0$ iff $i_k\geq j_k$ for all $k$. It follows
that ${pi\choose (p-1)i}=0$ for all $i\geq 0$.

{\it Remark}. Though $\der^p=0$ but $\frac{\der^p}{p!}\neq 0$
since $\frac{\der^p}{p!}(x^p)=1$ and $\frac{\der^p}{p!}$ is {\em
not} a derivation as
$\frac{\der^p}{p!}(x)x^{p-1}+x\frac{\der^p}{p!}(x^{p-1})=0$
(recall that if $\d $ a derivation then so is $\d^p$).

A higher derivation $\d =(\d_i)\in \HS_K(A)$ is called {\em
iterative} if $\d_i\d_j={i+j\choose i}\d_{i+j}$ for all $i,j\geq
0$. Then a  direct computation shows that 
\begin{equation}\label{dip=0}
\d_{i}^p=0\;\;\; {\rm  for \; all\;} \;\;  i\geq 1,
\end{equation}
$\d_{i}^p=\d_{i}\cdots \d_{i}= {2i\choose i}{3i\choose 2i}\cdots
{pi\choose (p-1)i}\d_{pi}=0\d_{pi}=0$. For $i=1$, we have
$\d_1^p=0$. The higher derivation $(\di! )\in \HS_K(K[x])$ is
iterative as follows directly from the definition of $(\di! )$.

Given $\d \in \Der_K(A)$, then $\d^p\in \Der_K(A)$ and, for any
$a\in A$, $(a\d )^p=a^p\d^p+(a\d)^{p-1}(a)\d$ (the {\em
Hochschild's formula}, \cite{Ma}, 25.5). In the algebra $\CD
(P_n)$, for each $i=1, \ldots , n$, $\der_i^p=0$, and therefore,
\begin{equation}\label{xipd}
(x_i\der_i)^p=x_i\der_i.
\end{equation}

{\bf The higher derivations $(1, \frac{\der_i}{1!},
\frac{\der^2_i}{2!}, \ldots )\in \HS_K(P_n)$, $i=1,\ldots , n$}.
The $K$-algebra homomorphism $P_n\ra P_n[[t]]$, $$f(x_1, \ldots ,
x_n)\mapsto f(x_1, \ldots , x_{i-1}, x_i+t, x_{i+1}, \ldots ,
x_n)=\sum_{i\geq 0}\dik! (f)t^k,$$ gives the higher derivation
 $(1, \frac{\der_i}{1!},
\frac{\der^2_i}{2!}, \ldots )\in \HS_K(P_n)$. If ${\rm char}
(K)=0$ then $\dik!$ means $(k!)^{-1} \der^k_i$, but if ${\rm char}
(K)=p>0$ then repeating the proof of (\ref{dxj=bin}), we see that
\begin{equation}\label{nbinijchp}
\dik! (x_j^l)=\d_{ij}{l\choose k}x_j^{l-k}
\end{equation}
for all $l\geq k\geq 1$  and $1\leq i,j\leq n$ where $\d_{ij}$ is
the Kronecker delta; and $\dik!
 (x_j^l)=0$ if $k>l$.

For an ideal $I$ of the polynomial algebra $P_n$, $I[[t]]$ is an
ideal of the algebra $P_n[[t]]$, and the factor algebra
$P_n[[t]]/I[[t]]\simeq P_n/I[[t]]$. The set $\HSK (P_n,I):=\{ e\in
\HSK (P_n)\, | \, e(I[[t]])=I[[t]]\}$ is a subgroup of the group
$\HSK (P_n)$. Note that 
\begin{equation}\label{eIt=It}
e(I[[t]])=I[[t]] \Leftrightarrow  e(I[[t]])\subseteq I[[t]].
\end{equation}
The implication $(\Rightarrow )$ is obvious. The reverse
implication follows immediately from the fact that, for any
$K$-algebra $A$ and any  higher derivation $\d =(1, \d_1, \ldots
)\in \HSK (A)$, the inverse automorphism to the automorphism
$e(a):=\sum \d_i(a)t^i$ has the form 
\begin{equation}\label{e-1a}
e^{-1}(a)=a+\d_1'(a)t+\cdots +\d_i'(a)t^i+\cdots
\end{equation}
where $\d_i'=\sum \pm \o_{ij}$ is a finite sum where $\o_{ij}$ is
a product of certain $\d_k$'th. Note that $e(I)\subseteq I[[t]]$
iff $\d_i (I)\subseteq I$ for all $i\geq 1$ where $e(p)=\sum
\d_i(p)t^i$. Then the inclusion $e(I)\subseteq I[[t]]$ implies the
inclusions $e^{-1}(I)\subseteq e^{-1}(I[[t]])\subseteq I[[t]]$.
Therefore, $e^{\pm 1}(I[[t]])\subseteq I[[t]]$, and so
$e(I[[t]])=I[[t]]$. Therefore, $\HSK (P_n,I):=\{ e\in \HSK (P_n)\,
| \, e(I[[t]])\subseteq I[[t]]\}$. Then it follows that the set
$\hsK (P_n,I):=\{ e\in \HSK (P_n)\, | \, e(P_n)\subseteq
P_n+I[[t]]\}$ is a normal subgroup of $\HSK (P_n,I)$. The kernel
of the canonical homomorphism of groups 
\begin{equation}\label{hHSI}
\HSK (P_n,I)\ra \HSK (P_n/I),\;\; e\mapsto(p+I\mapsto
e(p)+I[[t]]),
\end{equation}
is equal to $\hsK (P_n,I)$, and so the map in the proposition is a
group monomorphism.
\begin{proposition}\label{HSPnfI}
The map $$\HSK (P_n,I)/\hsK (P_n,I)\ra \HSK (P_n/I),\;\;\; e\cdot
\hsK (P_n,I)\mapsto (p+I\mapsto e(p)+I[[t]]),$$ is an isomorphism
of groups.
\end{proposition}

{\it Proof}. It remains to show that the map is surjective. Given
$\overline{e}\in \HSK (P_n/I)$. For each $i=1,\ldots , n$, $
\overline{e} (x_i+I)=x_i+I+\sum_{j\geq 1}(p_{ij}+I)t^j$ for some
$p_{ij}\in P_n$. The automorphism $\overline{e}$ can be extended
to an element $e\in \HSK (P_n)$ setting $e(x_i)= x_i+\sum_{j\geq
1}p_{ij}t^j$ such that the element $\overline{e}$ is the image of
the element $e$ under the homomorphism (\ref{hHSI}). This proves
the surjectivity.  $\Box $

Suppose, for a moment, that $\char (K)=0$. Then the ring of
differential operators $\CD (P_n)$ is, so-called, the {\bf Weyl
algebra} $A_n=K\langle x_1, \ldots , x_n, \der_1, \ldots ,
\der_n\rangle$, has the {\em standard} filtration $\{
A_{n,i}=\bigoplus_{|\alpha | +| \beta | \leq i} Kx^\alpha
\der^\beta \}_{i\geq 0}$ associated with the set of canonical
generators $x_j $, $\der_j:=\frac{\der}{\der x_j}$, where $\alpha
=(\alpha_1, \ldots ,\alpha_n)$, $\beta = (\beta_1, \ldots ,
\beta_n)\in \mathbb{N}^n$, $x^\alpha :=x_1^{\alpha_1}\cdots
x_n^{\alpha_n}$, $\der^\beta :=\der_1^{\beta_1}\cdots
\der_n^{\beta_n}$, $|\alpha | :=\alpha_1+\cdots +\alpha_n$.

The polynomial algebra $P_n$ as the left $A_n$-module has the
 {\em standard} filtration $\{ P_{n,i}:=A_{n,i}K=\bigoplus_{|\alpha |\leq
i} Kx^\alpha\}$. For each $i\geq 0$,
$$ \dim \, A_{n,i}={i+2n\choose 2n}=\frac{(i+2n)(i+2n-1)\cdots
(i+1)}{(2n)!}\;\;\; {\rm and}\;\;\; \dim \, P_{n,i}={i+n\choose
n},$$ and so $\GK (A_n)=2n$ and $\GK (P_n)=n$. The associated
graded algebra
$$\gr (A_n):= \oplus_{i\geq 0}A_{n,i}/A_{n,i-1}=K[\bx_1, \ldots , \bx_n, \bder_1, \ldots , \bder_n]$$ is isomorphic as a graded
algebra to a polynomial algebra in $2n$ variables with usual
grading.

If $\char (K)=p>0$ then the ring of differential operators  $\CD
:=\CD (P_n)$ on $P_n$ is an algebra  generated by $x_1, \ldots ,
x_n$ and commuting {\em higher derivations} $\frac{\der_i^k}{k!}$,
$i=1, \ldots , n$ and $k\geq 1$ that satisfy the following {\em
defining} relations: 
\begin{equation}\label{DPndef}
[x_i,x_j]=[\frac{\der_i^k}{k!},\frac{\der_j^l}{l!}]=0,\;\;\;
\frac{\der_i^k}{k!}\frac{\der_i^l}{l!}={k+l\choose
k}\frac{\der_i^{k+l}}{(k+l)!}, \;\;\; [\frac{\der_i^k}{k!},
x_j]=\d_{ij}\frac{\der_i^{k-1}}{(k-1)!}, \end{equation}
 for all
$i,j=1, \ldots , n$ and $k,l\geq 1$ where $\d_{ij}$ is the
Kronecker delta and $\frac{\der_i^0}{0!}:=1$. We will use also the
following notation: $\deri^{[k]} :=\frac{\der_i^k}{k!}$.

{\bf The involution $*$}. The $K$-linear map $*:\CD \ra \CD $,
$x_i\mapsto x_i$, $\deri^{[j]}\mapsto (-1)^j\deri^{[j]}$, $i=1,
\ldots , n$, $j\geq 1$, is an {\em involution} of the algebra $\CD
$ ($a^{**}=a$ and $(ab)^*=b^*a^*$). So, the algebra $\CD$ is a
symmetric object, its `left' and `right' properties are the
`same'. In particular, the categories of left and right
$\CD$-modules are `identical'.

{\bf The $\CD (P_n)$-module $P_n$}. The polynomial algebra $P_n$
is a (left) $\End_K(P_n)$-module, $\CD (P_n)$ is a subalgebra of
$\End_K(P_n)$, and so $P_n$  is a (left) $\CD (P_n)$-module. The
$\CD (P_n)$-module $P_n$ is canonically isomorphic to the factor
module $\CD (P_n)/\sum_{0\neq \beta \in \mathbb{N}^n}\CD
(P_n)\frac{\der^\beta}{\beta !}$.

The algebra $\CD $ is not finitely generated and not a (left or
right) Noetherian. It follows from the relations that $F=\{ F_i:=
\bigoplus_{|\alpha | +|\beta | \leq i} Kx^\alpha \frac{\der^\beta
}{\beta!}\}$ is a finite dimensional filtration for the algebra
$\CD$ where $\frac{\der^\beta }{\beta!}:= \frac{\der_1^{\beta_1}
}{\beta_1!}\cdots \frac{\der_n^{\beta_n} }{\beta_n!}$.

{\it Definition}. The filtration $F=\{ F_i\}$
 is called the {\bf canonical filtration} on $\CD (P_n)$. If $M=
 \CD M_0$ ($\dim_K(M_0)<\infty$) is a finitely generated
 $\CD$-module then the finite dimensional filtration $\{ M_i:=
 F_iM_0\}_{i\geq 0}$ is called the {\bf canonical filtration} of
 $M$.

 In characteristic zero, this filtration coincide with the standard
filtration $\{ A_{n,i}\}$. The canonical filtration is
$*$-invariant: $F^*=F$, i.e. $F^*_i=F_i$ for all $i\geq 0$.

{\it Definition}. A (finite dimensional) filtration $\{
F_i'\}_{i\geq 0}$ on the algebra $\CD $ which is {\em equivalent}
to the {\em canonical}  filtration $F$ is called a {\bf filtration
of standard type} of $\CD$. If $M=\CD M_0$ $(\dim_K(M_0)<\infty$)
is a finitely generated $\CD $-module then
 the finite dimensional filtration  $\{ F_i'M_0\}$ is called a {\bf filtration
of standard type} of $M$.

Filtrations of standard type in prime characteristic are correct
generalizations of standard filtrations in zero characteristic.
Each canonical filtration is a filtration of standard type.

The polynomial algebra $P_n$ as a left $\CD $-module has the
filtration of standard type $\{ F_iK=P_{n,i}\}$. Since $\dim \,
F_i=\dim \, A_{n,i}$ and $\dim \, F_iK=\dim \, P_{n,i}$,
$$ \Dim (\CD (P_n))= 2n \;\;\; {\rm and}\;\;\; \Dim (P_n)=\GK
(P_n)=n.$$ Note that the Gelfand-Kirillov dimension $\GK (\CD
(P_n))=n$, not $2n$. The associated graded algebra $\gr \, \CD
:=\oplus_{i\geq 0} F_i/F_{i-1} $ ($F_{-1}:=0$) is a commutative
algebra which is not a finitely generated algebra, the nil-radical
$\Gn $ of the algebra $\gr \, \CD $ is equal to $\sum_{\alpha,
\beta \in \mathbb{N}^n, |\beta |>0} K\bx^\alpha
\frac{\bder^\beta}{\beta!}$ and $\gr (\CD )/\Gn \simeq K[\bx_1,
\ldots , \bx_n]$ is a polynomial algebra.

\begin{theorem}\label{pfdWeyl=1}
Let $\CD (P_n)$ be the ring of differential operators  with
polynomial coefficients $P_n=K[x_1, \ldots , x_n]$ over  a field
$K$ of characteristic $p>0$. Then $\fdim (\CD (P_n))=1$ where
$\fdim =\fd, \lfd , \ud , \lud$.
\end{theorem}

{\it Proof}. Let $\nu =\nu_F$ be the return function of the
algebra $\CD =\CD (P_n)$ associated with the canonical filtration
$F=\{ F_i\}$ on $\CD $. Let us prove that $\nu (i)\leq i$ for all
$i\geq 0$. We use induction on $i$. The case $i=0$ is obvious as
$F_0=K$. Suppose that $i>0$ and the statement is true for all
$i'<i$. Let $a\in F_i\backslash F_{i-1}$. Then $a=\sum a_{\alpha
\beta} x^\alpha \frac{\der^\beta}{\beta !}$ with $|\alpha | +|
\beta | \leq i$ and $a_{\alpha \beta}\in K$. If there exists a
coefficient $a_{\alpha \beta}\neq 0$ for some $\beta \neq 0$, i.e.
$\beta_j \neq 0$  for some $j$, then applying the inner derivation
$\ad \, x_j$ of the algebra $\CD $ to the element $a$ we have a
{\em nonzero} element $x_ja-ax_j\in F_{i-1}$, then induction gives
the result.

Now, we have to consider the case where $a_{\alpha \beta}=0$ for
all $\beta\neq 0$, that is $a\in P_{n,i}\backslash P_{n,i-1}$.
Then there exists a variable, say $x_j$, such that $\deg_{x_j}
(a)>0$ (the degree in $x_j$) and a unique integer $k\geq 0$ such
that $p^k\leq \deg_{x_j}(a)<p^{k+1}$. Then applying the inner
derivation $\ad \, \der_j^{[p^k]}$ of the algebra $\CD $ to the
element $a$ we have a {\em nonzero} element $\der_j^{[p^k]}
a-a\der_j^{[p^k]}\in F_{i-p^k}$, and again induction finishes the
proof of the fact that $\nu (i)\leq i$ for all $i\geq 0$. It
follows that $1\geq \fd (\CD )\geq \fdim (\CD )$ (see
(\ref{4rf})).

The $\CD $-module $P_n$ has dimension $\Dim (P_n)=n$. By Theorem
\ref{FFI},
\begin{eqnarray*}
2n&=&\Dim (\CD )\leq  \Dim (P_n) (\fdim (\CD )+\max \{ \fdim (\CD ), 1\})\\
&\leq &  n(\fdim (\CD )+\max \{ \fdim (\CD ), 1\}),
\end{eqnarray*}
and so  $\fdim (\CD )\geq 1$. Then $\fdim (\CD )=1$, as required.
$\Box $

\begin{theorem}\label{BIDPnp}
({\bf An analogue of the inequality of Bernstein}) Let $M$ be a
nonzero finitely generated $\CD (P_n)$-module where $K$ is a field
of characteristic $p>0$. Then $\Dim (M)\geq n $.
\end{theorem}

{\it Proof}. By  Theorems \ref{FFI} and \ref{pfdWeyl=1},
$$ \Dim (M)\geq \frac{\Dim (\CD (P_n))}{1+1}=\frac{2n}{n}=n.\;\;
\; \Box$$

So, for any nonzero finitely generated $\CD (P_n)$-module $M$:
$n\leq \Dim (M)\leq 2n$. Any intermediate {\em natural} number
occurs: for $n=1$, $\Dim (P_1)=1$ and $\Dim (\CD (P_1))=2$. For
arbitrary $n$, $\CD (P_n)=\CD (P_1)\otimes \cdots \otimes \CD
(P_1)$ ($n$ times). Clearly, $\Dim (P_1^{\otimes s}\otimes \CD
(P_1)^{\otimes (n-s)})=s+2(n-s)=2n-s$. When $s$ runs through $0,
1, \ldots n$, the number $2n-s$ runs through $n, n+1, \ldots ,
2n$.

We say that a function $f:\mathbb{N}\ra \mathbb{N}$ has {\bf
polynomial growth} if there exists a polynomial $p(t)\in
\mathbb{Q}$ such that $f(i)\leq p(i)$ for $i\gg 0$. If a function
has polynomial growth so does any function which is equivalent to
it. We say that a filtration $\{ V_i\}$ has {\em polynomial
growth} if the function $\dim_K\, V_i$ has.

{\it Definition}. A finitely generated $\CD (P_n)$-module $M$ is
called a {\bf holonomic} module if there is a filtration of
standard type on $M$ that has polynomial growth.

Since all filtrations of standard type are equivalent, a finitely
generated $\CD (P_n)$-module $M$ is holonomic iff all filtrations
of standard type on $M$ has polynomial growth. It follows from the
definition that the class of holonomic $\CD (P_n)$-modules is
closed under sub- and factor modules, and under finite direct
sums.


\section{Description of finitely presented $\CD (P_n)$-modules, multiplicity
and (Hilbert) almost polynomials}\label{MACDP}

In this section, we show that in prime characteristic {\em
finitely presented} $\CD (P_n)$-modules  behave similarly as {\em
finitely generated} $\CD (P_n)$-modules in characteristic zero:
for {\em each finitely presented} $\CD (P_n)$-module $M$, the
Poincare series of it is a {\em rational} function, though its
Hilbert function is not a polynomial but an {\em almost
polynomial} and the degree of it coincides with the dimension
$\Dim (M)$ of $M$ (and if $M\neq 0$ then the dimension $\Dim (M)$
can be {\em any natural number} from the interval $[n,2n]$, this
gives another proof of an analogue of the inequality of Bernstein
for finitely presented $\CD (P_n)$-modules, Theorem
\ref{23Jun05}), and the multiplicity exits for $M$ (Theorem
\ref{23Jun05}). The differences are  $(i)$ in prime
characteristic,  finitely presented $\CD (P_n)$-modules have
transparent structure and are described by Theorem \ref{23Jun05},
but in characteristic zero the category of finitely generated $\CD
(P_n)$-modules is  far from being well-understood, $(ii)$ for each
{\em natural} number $d$ such that $n<d\leq 2n$, there exists a
cyclic finitely presented $\CD (P_n)$-module $M$ with $\Dim (M)=d$
and with {\em arbitrary small} multiplicity, Lemma \ref{a22Jul05}
(in characteristic zero, multiplicity is a {\em natural} number),
though the multiplicity of {\em  every holonomic finitely
presented} $\CD (P_n)$-module is a {\em natural} number (Theorem
\ref{t12Jul05}), $(iii)$ and what is completely unexpected is that
{\em each simple finitely presented} $\CD (P_n)$-module is {\em
holonomic} (Corollary \ref{1c23Jun05}), and if, in addition, the
field $K$ is {\em algebraically closed} that the multiplicity is
{\em always} $1$ (Corollary \ref{c8Jul05}).

Let $K$ be an {\em arbitrary} field.

{\bf Quasi and almost polynomials}. A function
$f:\mathbb{N}\ra\mathbb{N}$ is called a {\bf quasi-polynomial}
with a {\em period} $k$ if there exist $k$ polynomials $p_s(t)\in
\mathbb{Q}[t]$, $s\in \mathbb{Z}/k\mathbb{Z}$, such that
$$ f(i)=p_{\bi }(i)\;\;\; {\rm for \; all}\;\;\; i\gg 0,$$
where $\bi :=i+k\mathbb{Z}\in\mathbb{Z}/k\mathbb{Z}$. We say that
the  quasi-polynomial $f$ has {\em coefficients from a set}
$S\subseteq \mathbb{Q}$ if all the polynomials $p_i$ belong to
$S$.

 A quasi-polynomial $f=(p_0, \ldots ,
p_{k-1})$ is called an {\bf almost polynomial} if
 all the polynomials $p_i$ have the {\em same} degree $\deg (f)$ and the
 {\em same} leading coefficient $\lc (f)$ which are called respectively
  the {\em degree} and the
{\em leading coefficient} of $f$. $e(f):=\deg (f)!\lc (f)$ is
called the {\em multiplicity} of $f$. Then $f(i)=\frac{e(f)}{
d!}i^d+\cdots $, $i\gg 0$ where $d=\deg (f)$, and the three dots
mean `smaller' terms.

A function $f:\mathbb{N}\ra\mathbb{N}$ is called a {\bf somewhat
polynomial} if there are two polynomials $p,q\in \mathbb{Q}[t]$ of
the {\em same} degree $d$ such that $p(i)\leq f(i)\leq q(i)$ for
all $i\gg 0$. Then $d$ is called the {\em degree} of $f$.

{\bf Somewhat commutative algebras}. A $K$-algebra $R$ is called a
{\em somewhat commutative} algebra if it has a {\em finite
dimensional} filtration $R=\cup_{i\geq 0} R_i$ such that  $1\in
R_0$ and the associated graded algebra $\gr \, R$ is an {\em
affine commutative}
 algebra. Then the algebra $R$ is a Noetherian finitely
generated algebra. Let us choose homogeneous  $R_0$-algebra
generators of the
 $R_0$-algebra $\gr\, R:= \oplus_{i\geq 0} R_i/R_{i-1}$, say $y_1, \ldots , y_s$ of graded degrees $1\leq
k_1,\ldots , k_s$ respectively (that is $y_i\in
R_{k_i}/R_{k_i-1}$). A filtration $\G =\{ \G_i \}_{i\geq 0}$ of an
$R$-module $M$ is called a {\em good} filtration if the associated
graded $\gr\, R$-module $\gr_\G (M):=\oplus_{i\geq 0}
\G_i/\G_{i-1}$ is finitely generated. An $R$-module $M$ has a good
filtration iff it is finitely generated, and if $\{ \G_i\}$ and
$\{ \O_i\}$ are two good filtrations on $M$, then there exists a
natural number $t$ such that $\G_i\subseteq \O_{i+t}$ and
$\O_i\subseteq \G_{i+t}$ for all $i\geq 0$. If an $R$-module $M$
is finitely generated and $M_0$ is a finite dimensional generating
subspace of $M$, then the standard filtration $\{ R_iM_0\}$ is
good. The first two statements of
 the following lemma are well-known by specialists (see their proofs in
\cite{BavIzv}, Theorem 3.2 and Proposition 3.3).

\begin{lemma}\label{izvres}
Let $R=\cup_{i\geq 0}R_i$ be a somewhat commutative algebra,
$k=\lcm (k_1,\ldots ,k_s)$, $M$ be a finitely generated $R$-module
of Gelfand-Kirillov dimension $d=\GK (M)$ with good filtration $\G
=\{ \G_i\}$. Then
\begin{enumerate}
\item $\dim_K(\G_i )=\frac{e(M)}{d!}i^d+\cdots $ is an almost
polynomial of period $k$ with coefficients from
$\frac{1}{k^dd!}\mathbb{Z}$ where $e(M)\in\frac{1}{k^d}\mathbb{N}$
is called the multiplicity of $M$. The multiplicity does not
depend on the  choice of the good filtration $\G$. \item The
Poincare series of  $M$, $P_M(\o ):=\sum_{i\geq
0}\dim_K(\G_i)\o^i\in \mathbb{Q}(\o)$, is a rational function of
the form $\frac{f(\o )}{\prod_{i=1}^s(1-\o^{k_i})}$ where $f(\o
)\in \mathbb{Q}[\o ]$. The $P_M(\o )$  has the pole  of order
$d+1$ at $\o =1$, and $e(M)=e_{P_M}$ where $e_{P_M}:=(1-\o
)^{d+1}P_M(\o )|_{\o =1}$ is called the multiplicity of $P_M$.
\item If the elements $y_1, \ldots , y_t$ are nilpotent then the
two previous statements hold replacing the number $k$ by $\lcm
(k_{t+1}, \ldots k_s)$. \item In particular, if all non-nilpotent
generators of the algebra $\gr \, R$ have degree $1$ then $P_M(\o
)=\frac{f(\o )}{(1-\o )^{d+1}}$ for some polynomial $f(\o )\in
\mathbb{Q}[\o ]$ such that $e(M)=f(1)\in \mathbb{N}$, and
$\dim_K(\G_i)= \frac{e(M)}{d!}i^d+\cdots $ for $i\gg 0$ is a
polynomial of degree $d$ with coefficients from
$\frac{1}{d!}\mathbb{Z}$.
\end{enumerate}
\end{lemma}

{\it Proof}. $3$. Repeat the original proof taking into account
that the algebra $R_0\langle y_1, \ldots , y_t\rangle$ is finite
dimensional.

$4$. This statement is obvious.  $\Box $

\begin{corollary}\label{1izvres}
Let $P,Q\in\mathbb{Q}(\o)$ be rational functions having the pole
at $\o =1$ of order $n$ and $m$ respectively. Let $e_P>0$ and
$e_Q>0$ be the multiplicities of $P$ and $Q$ respectively. Then
$n+m-1$ and $e_Pe_Q$ are the order of the pole at $\o =1$ and the
multiplicity of the rational function $(1-\o )PQ$ respectively.
\end{corollary}

{\it Proof}. The first statement is trivial, then the multiplicity
of the rational function $(1-\o )PQ$ is equal to $(1-\o
)^{n+m-1}(1-\o )PQ|_{\o =1}=(1-\o )^nP(1-\o )^mQ|_{\o =1}=e_Pe_Q$.
$\Box $

Till the end of the section  $K$ is an {\em arbitrary} field of
{\em characteristic} $p>0$.

{\bf The algebras} $\L_\e$. For each $\e =(\e_1, \ldots , \e_n)\in
\{\pm 1\}^n$, consider the commutative subalgebra $\L_\e
:=\L_{\e_1}\t \cdots \t\L_{\e_n}$ of the ring of differential
operators $\CD (P_n)=\CD (P_1)\t \cdots \t \CD (P_1)$ where
$\L_{\e_i}:=K[x_i]$, if $\e_i=1$, and $\L_{\e_i}:=K[\deri^*]=K[
\deri^{[1]}, \deri^{[2]}, \ldots ]$, if $\e_i=-1$. The algebra
$\L_\e$  is a tensor product of commutative algebras $\L_{\e_i}$.
Each of the tensor multiples is a naturally $\mathbb{N}$-graded
algebra: $K[x_i]=\oplus_{j\geq 0}Kx_i^j$ and
$K[\deri^*]=\oplus_{j\geq 0}K\deri^{[j]}$. So, the algebra $\L_\e$
is a naturally $\mathbb{N}^n$-graded algebra with respect to the
tensor product of the $\mathbb{N}$-gradings:
$$\L_\e=\bigoplus_{\alpha \in
\mathbb{N}^n}Kl^\alpha , \;\;\; l^\alpha l^\beta ={\alpha +\beta
\choose \beta}_\e l^{\alpha +\beta}\;\;\;{\rm for \; all}\;\;
\alpha , \beta \in \mathbb{N}^n,
$$
where $l^\alpha :=l^\alpha_\e :=l_1^{\alpha_1}\cdots
l_n^{\alpha_n}$, $l_i^{\alpha_i}:=x_i^{\alpha_i}$, if $\e_i=1$,
and $l_i^{\alpha_i}:=\der_i^{[\alpha_i]}$, if $\e_i=-1$. ${\alpha
+\beta \choose \beta}_\e  := \prod_{i=1}^n {\alpha_i +\beta_i
\choose \beta_i}_{\e_i}$ where ${i\choose j}_{-1}:={i\choose j}$
and ${i\choose j}_{1}:=1$. The $\e$-{\em binomials} are {\em
translation invariant}: 
\begin{equation}\label{ebtr}
{\alpha p^k\choose \beta p^k}_\e ={\alpha \choose \beta }_\e
\;\;\; {\rm for \; all}\;\; k\geq 0\;\; {\rm and}\;\;\; \alpha ,
\beta \in \mathbb{N}^n.
\end{equation}
For each $k\geq 0$, $\Lepk := \bigoplus_{\alpha \in
\mathbb{N}^n}Kl^{\alpha p^k}$ is a subalgebra of $\L_\e$. The
translation invariance of the $\e$-binomials implies that the
$K$-linear map
$$\L_\e \ra \Lepk , \;\; l^\alpha \mapsto l^{\alpha p^k}, \;\;\; \alpha \in \mathbb{N}^n, $$
is a $K$-{\em algebra isomorphism}.  There exists the descending
chain of subalgebras of $\L_\e$:
$$ \L_\e:=\L_{\e}^{[p^0]}\supset \L_{\e}^{[p]}\supset \cdots
\supset \Lepk \supset \cdots , \;\;\ \cap_{k\geq 0}\Lepk =K.$$
 For each $k\geq 0$, let $\L_{\e, [p^k]}:=\bigoplus_{\alpha <\bpk } Kl^\alpha = \L_{\e_1,
[p^k]}\t \cdots \t \L_{\e_n, [p^k]}$ where $\bpk := (p^k, \ldots ,
p^k)$, and $\alpha <\bpk$ means that $0\leq \alpha_1<p^k, \ldots
,0\leq \alpha_n<p^k$.  $\dim_K(\L_{\e, [p^k]})=p^{nk}$. The vector
space $\L_{\e, [p^k]}$ is an algebra iff $\e =(-1, \ldots, -1)$,
and, in this case, $\L_{(-1, \ldots , -1), [p^k]}=\L_{-1, [p^k]}\t
\cdots \t \L_{-1, [p^k]}$ is the tensor product of commutative
finite dimensional algebras where each tensor multiple, say
$i$'th, $$\L_{-1, p^k}= K\langle \deri^{[1]} \rangle \t K\langle
\deri^{[p]}\rangle  \t K\langle \deri^{[p^2]}  \rangle \t \cdots
\t K\langle \deri^{[p^{k-1}]}\rangle \simeq (K[t]/(t^p))^{\t k} $$
is the tensor product of commutative local finite dimensional
algebras since $K\langle \deri^{[p^{s}]}\rangle \simeq K[t]/(t^p)$
as $(\deri^{[p^{s}]})^p=0$, $s\geq 1$. Clearly,
\begin{eqnarray*}
\L_{\e, [p^0]}:&=&K\subset \L_{\e, [p]}\subset \cdots \subset
\L_{\e, [p^k]}\subset \cdots , \;\; \;\; \L_\e =\cup_{k\geq
0}\L_{\e,
[p^k]},\\
\L_\e &=& \L_{\e, [p^k]} \Lepk =\L_{\e, [p^k]}\t \Lepk = \Lepk \t
\L_{\e, [p^k]},\;\; \;\; k\geq 0,
\end{eqnarray*}
and $\L_{\e, [p^k]}\L_{\e, [p^l]}\subseteq \L_{\e, [p^{\max \{ k,
l\} }]}$ for $\e =(-1, \ldots , -1)$ and  all  $k,l\geq 0$.

The subalgebra $\L := \L_{-1}=K[\der^{[1]}, \der^{[2]}, \ldots ]$
of the algebra $\CD (K[x])$ is {\em not} a finitely generated
algebra (since $\L =\cup_{k\geq 0}\L_{[p^k]}$ is the union of its
proper subalgebras), it is {\em not} a domain (its nil-radical
$\Gn (\L )$ is equal to $\L_+:=\oplus_{j\geq 1}K\der^{[j]}$), it
is {\em not} a Noetherian algebra as
$$ \L_{[p],+}\t\L^{[p]}\subset
\L_{[p^2],+}\t\L^{[p^2]}\subset\cdots \subset
\L_{[p^k],+}\t\L^{[p^k]}\subset\cdots $$ is a {\em strictly
ascending} chain of ideals of the algebra $\L$ where
$\L_{[p^k],+}:= \oplus_{j=1}^{p^k-1}K\der^{[j]}$, and
$$ \L /(\L_{[p^k],+}\t\L^{[p^k]})\simeq (\L /\L_{[p^k],+})\t\L^{[p^k]}\simeq
K\t \L^{[p^k]}\simeq  \L^{[p^k]}\simeq \L, \;\; k\geq 0.$$ In
spite of the fact that the algebra $\L_\e$, $\e =(-1, \ldots ,
-1)$,  is `zero dimensional' ($\L_\e /\Gn (\L_\e )=\L_\e  / \L_{\e
, +}=K$), it has the rich non-trivial category of modules which in
turn the ring of differential operators $\CD (P_n)$ inherits as a
subcategory (via inducing).

The algebra $\L_\e$ is {\em Noetherian} iff $\e =(1, \ldots , 1)$
(in this case, it is $P_n$, a finitely generated Noetherian
domain). If $\e \neq (1, \ldots , 1)$ then the algebra $\L_\e$ is
not finitely generated, not Noetherian, and not a domain.

{\bf The subalgebras $\CD(P_n)^{[p^k]}$ and $\L_{[p^k]}\t P_n$}.
 For each $\e =(\e_1, \ldots , \e_n)\in \{ \pm 1\}^n$,
\begin{equation}\label{LeLme}
\CD (P_n)=\L_\e\t\L_{-\e}= (\bigoplus_{\alpha \in
\mathbb{N}^n}Kl^\alpha_\e )\t (\bigoplus_{\beta \in
\mathbb{N}^n}Kl^\beta_{-\e} )= \bigoplus_{\alpha , \beta \in
\mathbb{N}^n}Kl^\alpha_\e\t l^\beta_{-\e}
\end{equation}
 as $K$-modules. For
$k\geq 0$, the vector space $\CD (P_n)_\e^{[p^k]}:=\Lepk \t
\L_{-\e }^{[p^k]}= \L_{-\e }^{[p^k]}\t \Lepk$ is a subalgebra of
$\CD (P_n)$ canonically isomorphic to the $K$-algebra $\CD (P_n)$
via the $K$-algebra isomorphism:
$$ \CD (P_n)\ra\CD (P_n)_\e^{[p^k]}, \;\;\; l^\alpha_\e\t
l^\beta_{-\e}\mapsto l^{\alpha p^k} _\e\t l^{\beta p^k}_{-\e},
\;\;\; \alpha , \beta \in \mathbb{N}^n.$$ This follows directly
form the translation invariance of the $\e$-binomials and from the
defining relations (\ref{DPndef}) for the $K$-algebra $\CD (P_n)$
since the $K$-algebra $\CD (P_n)_\e^{[p^k]}$ is generated by the
elements $x_1^{p^k}, \ldots , x_n^{p^k}, \der_1^{[ip^k]}, \ldots ,
\der_n^{[ip^k]}$, $i\geq 1$,  that satisfy the relations
(\ref{DPndef}), these relations are defining because of the
decomposition (\ref{LeLme}). It is obvious that $\CD
(P_n)_\e^{[p^k]}=\CD (P_n)_{\e'}^{[p^k]}$ for all $\e$ and $\e'$
but the decompositions (\ref{LeLme})  are all distinct and they
will be used later in constructing various modules. So, we drop
the subscript $\e$.
 There exists a descending chain of
{\em isomorphic} subalgebras of $\CD (P_n)$: 
\begin{equation}\label{CDchd}
\CD (P_n):= \CD (P_n)^{[p^0]}\supset \CD (P_n)^{[p]}\supset \CD
(P_n)^{[p^2]}\supset \cdots \supset \CD (P_n)^{[p^k]}\supset
\cdots ,
\end{equation}
 and $\cap_{k\geq 0}\CD (P_n)^{[p^k]}=K$. So, the $K$-linear map
\begin{equation}\label{CDchd1}
\CD (P_n)\ra \CD (P_n)^{[p^k]}, \;\; x_i\mapsto x_i^{p^k}, \;\;
\der_i^{[j]}\mapsto \der_i^{[jp^k]},
\end{equation}
is a $K$-algebra {\em isomorphism}.

For each $k\geq 0$, the vector space $\CD (P_n)_{\e , [p^k]}:=
\L_{\e , [p^k]}\t \L_{-\e , [p^k]}= \L_{-\e , [p^k]}\t \L_{\e ,
[p^k]}$ has dimension $p^{2nk}$ over $K$. Again, it does not
depend on $\e$, and so we drop the subscript $\e$, $\CD
(P_n)_{[p^k]}:=\CD (P_n)_{\e , [p^k]}=\bigoplus_{\alpha <\bpk,
\beta <\bpk} Kx^\alpha \der^{[\beta]}=\bigoplus_{\alpha <\bpk,
\beta <\bpk} K \der^{[\beta]}x^\alpha$.
   Clearly,
$$ K=:\CD (P_n)_{[p^0]}\subset \CD (P_n)_{ [p]}\subset \CD (P_n)_{
[p^2]}\subset \cdots , \;\; \CD (P_n)=\cup_{k\geq 0}\CD (P_n)_{
[p^k]},$$
$$ \CD (P_n)=\CD (P_n)_{ [p^k]}\t \CD (P_n)^{[p^k]}=
\CD (P_n)^{ [p^k]}\t \CD (P_n)_{ [p^k]}, \;\;\; k\geq 0,$$ and $
\CD (P_n)_{ [p^k]}\t \CD (P_n)_{ [p^l]}\subseteq \CD
(P_n)_{[p^{k+l}]}$, for all $k,l\geq 0$.

{\bf The case $\e =(-1, \ldots , -1)$ and the subalgebras
$\L_{[p^k]}\t P_n$}. In this case, we write $\L :=\L_{(-1, \ldots
, -1)}$, $\L_{[p^k]}:= \L_{(-1, \ldots , -1), [p^k]}$, and
$\L^{[p^k]}:= \L_{(-1, \ldots , -1)}^{ [p^k]}$. Then $\CD (P_n)=\L
\t P_n=P_n\t \L$ and $\CD (P_n)=\cup_{k\geq 0}\L_{[p^k]}\t P_n$ is
a union of subalgebras:
$$ P_n:=\L_{[p^0]}\t P_n\subset \L_{[p]}\t
P_n\subset \L_{[p^2]}\t P_n\subset \cdots , \;\;\; \L_{[p^k]}\t
P_n=P_n\t \L_{[p^k]},$$
$$\CD (P_n)=\L^{[p^k]}\t (\L_{[p^k]}\t P_n)= (\L_{[p^k]}\t P_n)\t
\L^{[p^k]}, \;\;\; k\geq 0.$$

For each $k\geq 0$, the algebra $\L_{[p^k]}\t P_n$ is a free left
and right $P_n$-module of rank $p^{nk}$, it is a {\em finitely
generated Noetherian} algebra having the centre $Z_k:=K[x_1^{p^k},
\ldots , x_n^{p^k}]$. {\em The algebra $\L_{[p^k]}\t P_n$ is a
free $Z_k$-module of rank} $p^{2nk}$ since $\L_{[p^k]}\t
P_n=\L_{[p^k]}\t (\oplus_{\alpha <\bpk }Kx^\alpha )\t Z_k$. On the
algebra $\L_{[p^k]}\t P_n$, and one can  consider the {\em induced
filtration} from the canonical filtration $F=\{ F_i\}$ on the
algebra $\CD (P_n)$: 
\begin{equation}\label{FCTk}
 \CT_k=\{ \CT_{k,i}:=\L_{[p^k]}\t P_n\cap F_i=\bigoplus_{ \beta
<\bpk , |\alpha | +|\beta |\leq i} Kx^\alpha
\der^{[\beta]}=\bigoplus_{ \beta <\bpk , |\alpha | +|\beta |\leq
i} K \der^{[\beta]}x^\alpha \}.
\end{equation}
 The associated graded algebra $\gr
  (\L_{[p^k]}\t P_n)$ is naturally isomorphic (as a graded algebra)
to the {\em tensor product of the commutative algebras}
$\L_{[p^k]}\t P_n$ equipped with the tensor product of the induced
filtrations (from the canonical filtration on $\CD (P_n)$). In
particular, $\gr ( \L_{[p^k]}\t P_n)$ is an {\em affine
commutative} algebra with the nil-radical $\L_{[p^k], +}\t P_n$
(where  $\L_{[p^k], +}:= \oplus_{0\neq \beta \in \bpk}
K\der^{[\beta ]}$) which is a completely prime ideal (a prime
ideal is a completely prime if the factor ring modulo the ideal is
a domain)  since
$$\gr ( \L_{[p^k]}\t P_n)/ \L_{[p^k], +}\t P_n \simeq (\L_{[p^k]}/ \L_{[p^k], +})\t P_n
\simeq K\t P_n\simeq P_n.$$ The algebra  $\gr (  \L_{[p^k]}\t
P_n)=\oplus_{i\geq 0} G_i$  is positively graded  (with only
finitely many nonzero components) where
$$G_i=\bigoplus_{ \beta <\bpk ,
|\alpha | +|\beta |= i}  Kx^\alpha \der^{[\beta]}.$$ For a ring
$R$ and a natural number $n\geq 1$, $M_n(R)$ is the ring of
$n\times n$ matrices with entries from $R$.

\begin{lemma}\label{Tkpr}
Let $K$ be a field of characteristic $p>0$, and
$T_k:=T_{k,n}:=\L_{[p^k]}\t P_n$, $k\geq 0$. Then
\begin{enumerate}
\item  The algebra $T_k$ is a somewhat commutative algebra with
respect to the finite dimensional filtration $\CT_k
=\{\CT_{k,i}\}$ having the centre $Z_k=K[x_1^{p^k}, \ldots
,x_n^{p^k}]$ and $\GK (T_{k,n})=n$. In particular, $T_k$ is a
finitely generated Noetherian algebra, and $T_{k,n}=T_{k,1}^{\t
n}$. \item The Poincare series of $T_k$, $P_{T_k}=\sum_{i\geq
0}\dim_K(\CT_{k,i})\o^i=\frac{(1+\o +\o^2+\cdots +
\o^{p^k-1})^n}{(1-\o )^{n+1}}$ and the multiplicity
$e(T_k)=p^{kn}$. \item The Hilbert function of $T_k$ is, in fact,
a polynomial $\dim_K(\CT_{k,i})=\frac{p^{kn}}{n!}i^n+\cdots $,
$i\gg 0$. \item Let $\CZ_k=K(x_1^{p^k}, \ldots ,x_n^{p^k})$ be the
field of fractions of $Z_k$. Then $T_k':=\CZ_k\t_{Z_k}T_k\simeq
M_{p^{kn}}(\CZ_k)$, the matrix algebra. \item The algebra $T_k$ is
a prime algebra of uniform dimension $p^{kn}$, and the
localization $\CS^{-1}T_k$ of $T_k$ at the set $\CS$ of all the
non-zero divisors is isomorphic to the matrix algebra
$M_{p^{kn}}(\CZ_k)$. \item The algebra $T_k$ is preserved by the
involution $*$, $T_k^*=T_k$, and so the algebra $T_k$ is
self-dual. \item The algebra $T_k$ is  faithfully flat over its
centre. \item The left and right Krull dimension of the algebra
$T_k$ is $n$. \item The left and right global  dimension of the
algebra $T_k$ is $n$ but the global dimension of the associated
graded algebra $\gr (T_k)$  is $\infty$ if $k\geq 1$.

\end{enumerate}
\end{lemma}

{\it Proof}. $1$. $P_n$ is a subalgebra of $T_k$, and so $n=\GK
(P_n)\leq \GK (T_k)$. $T_k$ is a finitely generated $Z_k$-module,
and so $\GK (T_k)\leq \GK (Z_k)=n$. Therefore, $\GK (T_k)=n$.

$2$ and $3$. These statements are obvious (see Lemma \ref{izvres}
and Corollary \ref{1izvres}).

$4$. The $\CZ_k$-algebra $T_k'=\bigoplus_{\alpha , \beta
<\bpk}\CZ_k x^\alpha \der^{[\beta ]}$ has dimension $p^{2nk}$ over
the field $\CZ_k$. Consider the $T_k$-module $U:=T_k/(P_n\t
\L_{[p^k], +})\simeq P_n\t(\L_{[p^k]}/\L_{[p^k], +})\simeq P_n\t
K\simeq P_n \overline{1}$ where $\overline{1}$ is the canonical
generator of $U$. The $T_k'$-module
$U':=\CZ_k\t_{Z_k}U=\bigoplus_{\alpha  <\bpk}\CZ_k x^\alpha
\overline{1}$ is simple (use the action of $\der^{[\beta ]}$ on
$x^\alpha $), $\dim_{\CZ_k}(U')=p^{nk}=\sqrt{\dim_{\CZ_k}(T_k')}$,
and $\End_{T_k'}(U')\simeq \cap_{0<\beta <\bpk}\ann \,
\der^{[\beta]}\simeq \CZ_k$. Therefore, $T_k'\simeq
M_{p^{kn}}(\CZ_k)$.

$5$. Since $\CZ_k\backslash \{ 0\} \subseteq \CS$, it follows from
statement $4$ that $\CS^{-1}T_k\simeq T_k'\simeq
M_{p^{kn}}(\CZ_k)$, which implies that $T_k$ is a prime algebra of
uniform dimension $p^{kn}$.

$6$ and $7$. These statements are obvious.

$8$. By statement $6$, the left and right Krull dimension of $T_k$
are equal. By statement $7$, $\Kdim (T_k)\geq \Kdim (Z_k)=n$. The
algebra $T_k$ is a finitely generated $Z_k$-module, hence $\Kdim
(T_k)\leq \Kdim (Z_k)=n$, and so $\Kdim (T_k)=n$.

$9$. Straightforward.
 $\Box $

Since the canonical generators of the commutative
$\mathbb{N}$-graded algebra $\gr (\L_{[p^k]}\t P_n)$ that are not
nilpotent all have graded degree $1$ the next result follows from
Lemma \ref{izvres}.

\begin{lemma}\label{Tfgap}
Let $M$ be a finitely generated $\L_{[p^k]}\t P_n$-module of
Gelfand-Kirillov dimension $d$ equipped with a standard filtration
$\{ M_i:= \CT_{k,i}M_0\}$ where $M_0$ is a finite dimensional
generating space for $M$. Then $\dim_K(M_i
)=\frac{e(M)}{d!}i^d+\cdots $, $i\gg 0$,  is a  polynomial of
 degree $d$ with coefficients from $\frac{1}{d!}\mathbb{Z}$
where $e(M)\in\mathbb{N}$ is called the multiplicity of $M$. The
multiplicity does not depend on the choice of a good filtration
$\G$. The degree $d$ can be any natural number from the interval
$[0,n]$ (see Lemma \ref{Tkpr}).
\end{lemma}

Let $\CD :=\CD (P_n)$ and $T_k:=\L_{[p^k]}\t P_n$. Consider free
finitely generated (left) $\CD$-modules $\CD^\mu$ and $\CD^\nu$
where $\mu , \nu \geq 1$. The set $\Hom_\CD (\CD^\mu , \CD^\nu )$
of all the $\CD$-module homomorphisms from $\CD^\mu$ to $\CD^\nu$
 can be identified with the set of all $\mu \times \nu$ matrices
$M_{\mu , \nu}(\CD )$ with coefficients from $\CD $. On this
occasion, it is convenient to write homomorphisms on the {\em
right}. Then $M_{\mu , \nu}(\CD )=M_{\mu , \nu}(\cup_{k\geq
0}T_k)=\cup_{k\geq 0}M_{\mu , \nu}(T_k )$ is the union of matrix
algebras. Let $M$ be a {\em finitely presented} $\CD$-module, that
is $M=\coker (A)$ where $\CD^\mu \stackrel{A}{\longrightarrow}
\CD^\nu$, $v\mapsto vA$,  $v=(v_1, \ldots , v_\mu)$,  and $A\in
M_{\mu , \nu}(\CD )$. Then $A\in M_{\mu , \nu}(T_k)$ for some $k$,
and $M':= \coker (T_k^\mu \stackrel{A}{\longrightarrow} T_k^\nu )$
is a finitely presented $T_k$-module. Applying the {\em exact}
functor $\CD \t_{T_k}-$ to the exact sequence of $T_k$-modules
$T_k^\mu \stackrel{A}{\longrightarrow} T_k^\nu \ra M'\ra 0$ one
obtains the exact sequence of $\CD$-modules $\CD_k^\mu
\stackrel{A}{\longrightarrow} \CD_k^\nu \ra \CD \t_{T_k}M'\ra 0$.
 Therefore,

\begin{equation}\label{Mfpr}
M\simeq \CD\t_{T_k}M',
\end{equation}
and so {\em each finitely presented $\CD$-module is isomorphic to
an induced module from a finitely generated $T_k$-module}. The
next result describes finitely presented $\CD (P_n)$-modules and
gives as a result  an analogue of the inequality of Bernstein for
them.

\begin{theorem}\label{23Jun05}
Let $M$ be a nonzero finitely presented $\CD (P_n)$-module. Then
$M\simeq \CD\t_{T_k}M'$ for a finitely generated $T_k$-module
$M'$. Let $\{ M_i'\}$ be a standard filtration for the
$T_k$-module $M'$ from
 Lemma \ref{Tfgap} and $\dim_K(M_i'
)=\frac{e(M')}{d!}i^d+\cdots $ for $i\gg 0$ where $d=\GK (M')$.
Let $\{ M_i:=F_iM_0'\}$ be the filtration of standard type on the
$\CD (P_n)$-module $M$. Then
\begin{enumerate}
\item $\dim_K(M_i )=\frac{e(M')}{p^{kn}(n+d)!}i^{n+d}+\cdots $ is
an almost polynomial of period $p^k$ with coefficients from
$\frac{1}{p^{k(n+d)}(n+d)!}\mathbb{Z}$,  and
$e(M)=\frac{e(M')}{p^{kn}}\in \frac{1}{p^{kn}}\mathbb{N}$. \item
 The dimension $\Dim (M)=n+d\geq n$ is equal to $t-1$ where $t$ is  the order of the
 pole of the Poincare series $P_M(\o )=\sum_{i\geq
 0}\dim_K(M_i)\o^i$ at the point $\o =1$, and the multiplicity
 $e(M)=(1-\o )^{\Dim (M)+1}P_M(\o )|_{\o =1}$. The dimension $\Dim (M)$
  of $M$ can be any natural number from the interval $[n,2n]$.
\end{enumerate}
\end{theorem}

{\it Proof}. The subalgebra $\L^{[p^k]}$ of $\CD (P_n)$ has the
induced filtration $\{\L^{[p^k]}_i:=\L^{[p^k]}\cap
F_i=\bigoplus_{p^k|\beta |\leq i}K\der^{[p^k\beta]}\}$. Therefore,
$$ P:= \sum_{i\geq 0}\dim_K(\L^{[p^k]}_i)\o^i=\frac{1}{(1-\o
)(1-\o^{p^k})^n}\;\;\; {\rm and}\;\;\; e_P:= (1-\o)^{n+1}P|_{\o
=1}=\frac{1}{p^{kn}}.$$ It follows from the equality
$M=\L^{[p^k]}\t M'$ that $M_i=\sum_{j+k\leq i}\L^{[p^k]}_j\t
M_k'$. Therefore, $R:=\sum_{i\geq 0}\dim_K(M_i)\o^i=(1-\o )PQ$
where $Q:=\sum_{i\geq 0}\dim_K(M_i')\o^i$. By Corollary
\ref{1izvres}, $e(M)=e_R=e_Pe_Q=\frac{1}{p^{kn}}e(M')$ and $\Dim
(M)=n+d\geq n$, and so
 $\dim_K(M_i
)=\frac{e(M')}{p^{kn}(n+d)!}i^{n+d}+\cdots $, by Lemma
\ref{izvres}. The rest is obvious (Lemma \ref{izvres}). $\Box$

\begin{lemma}\label{a22Jul05}
For each $s=0, 1, \ldots, n-1$, $\CD (P_n)=\CD (P_1)\t \CD
(P_1)^{\t s}\t \CD (P_1)^{\t (n-s-1)}$. For each $k\in
\mathbb{N}$, consider the cyclic finitely presented $\CD
(P_n)$-module $M(k,s):= M(k)\t \CD (P_1)^{\t s}\t \L_{-1}^{\t
(n-s-1)}$ where $M(k):= \CD (P_1)\t_{T_k}T_k/T_k\L_{[p^k],+}$ is
the $\CD (P_1)$-module. Then $\Dim \, M(k,s)=n+1+s$ and
$e(M(k,s))=\frac{1}{p^k}$. So, the multiplicity of a non-holonomic
finitely presented $\CD (P_n)$-module can be arbitrary small (for
each possible dimension $n,\ldots , 2n$).
\end{lemma}

{\it Remark}. By contrast, the multiplicity of each {\em holonomic
finitely presented} $\CD (P_n)$-module is {\em natural} number
(Theorem \ref{t12Jul05}).

 {\it Proof}. The $T_k$-module
$N:=T_k/T_k\L_{[p^k],+}=P_1\overline{1} \simeq {}_{P_1}P_1$ has
the standard filtration $\{ N_i:=
\CT_{k,i}\overline{1}=\bigoplus_{j=0}^iKx_1^j\overline{1}\}$.
Therefore, $\dim_K(N_i)=i+1$, and so $\GK (N)=1$ and $e(N)=1$. By
Theorem \ref{23Jun05}, $\Dim (M(k))=2$ and
$e(M(k))=\frac{1}{p^k}$. The $\CD (P_1)^{\t s}\t \CD (P_1)^{\t
(n-s-1)}$-module $\CD (P_1)^{\t s}\t \L_{-1}^{\t (n-s-1)}$ has
dimension $2s+n-s-1=n+s-1$ and multiplicity $1$. Using Corollary
\ref{1izvres}, we have $\Dim (M)=2+n+s-1=n+1+s$ and
$e(M)=\frac{1}{p^k}\cdot 1=\frac{1}{p^k}$.  $\Box $

\begin{corollary}\label{c12Jul05}
Each short exact sequence $0\ra N\ra M\ra L\ra 0$ of finitely
presented $\CD (P_n)$-modules is obtained from a short exact
sequence $0\ra N'\ra M'\ra L'\ra 0$ of finitely presented
$T_k$-modules for some $k\geq 0$ by tensoring on $\CD
(P_n)\t_{T_k}-$.
\end{corollary}

{\it Proof}. Let $\CD :=\CD (P_n)$. The $\CD $-modules $N$, $M$,
and $L$ are finitely presented, so one can fix a commutative
diagram

$$\begin{array}{ccccccccc}
  0 & \ra  & \CD^{m_1} &  \ra  & \CD^{m_1+m_2} & \ra  & \CD^{m_2} & \ra & 0 \\
 &   & \downarrow &  & \downarrow &  & \downarrow &  &  \\
  0 & \ra  & \CD^{n_1} & \ra  & \CD^{n_1+n_2} & \ra  & \CD^{n_2} & \ra &
  0\\
 &   & \downarrow &  & \downarrow &  & \downarrow &  &  \\
0 & \ra  & N & \ra  & M & \ra  & L & \ra & 0 \\
&   & \downarrow &  & \downarrow &  & \downarrow &  &  \\
&   & 0 &  & 0 &  & 0&  &  \\
\end{array}
$$
with exact rows and columns. One can find a (large) $k$ such that
all the matrices that correspond to the (six) maps between $D^*$'s
 have coefficients from the algebra $T_k$. The diagram above is
 obtained from the following commutative diagram (with exact rows
 and columns) of $T_k$-modules (with the same matrices = maps)

$$\begin{array}{ccccccccc}
  0 & \ra  & T_k^{m_1} & \ra  & T_k^{m_1+m_2} & \ra  & T_k^{m_2} & \ra & 0 \\
 &   & \downarrow &  & \downarrow &  & \downarrow &  &  \\
  0 & \ra  & T_k^{n_1} & \ra  & T_k^{n_1+n_2} & \ra  & T_k^{n_2} & \ra &
  0\\
 &   & \downarrow &  & \downarrow &  & \downarrow &  &  \\
0 & \ra  & N' & \ra  & M' & \ra  & L' & \ra & 0 \\
&   & \downarrow &  & \downarrow &  & \downarrow &  &  \\
&   & 0 &  & 0 &  & 0&  &  \\
\end{array}
$$
by tensoring on $\CD \t_{T_k}-$.  $\Box $

\begin{corollary}\label{1c23Jun05}
Each simple finitely presented $\CD (P_n)$-module is holonomic.
\end{corollary}

{\it Proof}. Let $M$ be a simple finitely presented $\CD
(P_n)$-module. By Theorem \ref{23Jun05}, $M\simeq \CD
(P_n)\t_{T_k}M'$ for a finitely generated $T_k$-module $M'$ which
must be simple. The algebra $T_k$ is a somewhat commutative
algebra which is finitely generated module over its centre $Z$
which is an affine algebra. Therefore, by the Quillen's Lemma
(\cite{MR}, 9.7.3), every element of $\End_{T_k}(M')$ is
algebraic, this implies that each simple $T_k$-module is finite
dimensional over the field $K$, and so $d=\GK (M')=0$. By Theorem
\ref{23Jun05}, $M$ is a holonomic $\CD (P_n)$-module.  $\Box $

\begin{theorem}\label{cth23Jun05}
Let $M$ be a nonzero finitely presented $\CD (P_n)$-module. The
following statements are equivalent.
\begin{enumerate}
\item $M$ is a holonomic $\CD (P_n)$-module.\item $\Dim (M)=n$.
\item $\Dim (M)<n+1$.
\end{enumerate}
\end{theorem}

{\it Proof}. The implications $1\Rightarrow 2\Rightarrow 3$ are
obvious, and the implication $3\Rightarrow 1$ follows from Theorem
\ref{23Jun05}. $\Box $

{\it Remarks}. 1. If a finitely generated $\CD (P_n)$-module $M$
is not finitely presented that Theorem \ref{cth23Jun05} is not
true. There exists a cyclic non-holonomic $\CD (P_n)$-module $M$
with $\Dim (M)=n$ (Proposition \ref{Mbknh}), and there are plenty
of cyclic $\CD (P_n)$-modules having dimension $d$ such that
$n<d<n+1$ (Theorem \ref{3Jl05}).

2. In characteristic zero, the multiplicity of a holonomic $\CD
(P_n)$-module is a natural number, so it can't be arbitrary small.
This is the reason why each holonomic $\CD (P_n)$-module has
finite length. Though the same  is true in prime characteristic
(Theorem \ref{hDPnfl}), Theorem \ref{23Jun05} does {\em not} give
a {\em uniform} lower bound for multiplicity of holonomic finitely
presented $\CD (P_n)$-modules, so one {\em cannot repeat} the
arguments of the characteristic zero case even for finitely
presented modules. Note that there are plenty of holonomic modules
that are not finitely presented.

\begin{theorem}\label{12Jul05}
Let $K$ be a field of characteristic $p>0$, and $0\ra N\ra M\ra
L\ra 0$ be a short exact sequence of finitely presented $\CD
(P_n)$-modules. Then
\begin{enumerate}
\item there exist finite dimensional filtrations $\{ N_i\}$, $\{
M_i\}$, and $\{ L_i\}$ on the modules $N$, $M$, and $L$
respectively such that the last two are filtrations of standard
type and the first one is strongly equivalent to a filtration of
standard type on $N$ and such  that
$\dim_K(M_i)=\dim_K(N_i)+\dim_K(L_i)$, $i\geq 0$. \item
$\Dim_K(M)=\max \{ \Dim(N), \Dim(L)\}$. \item Precisely one of the
following statements is true
\begin{enumerate}
\item $\Dim (N)<\Dim (M)=\Dim (L)$ and $e(M)=e(L)$, \item $\Dim
(L)<\Dim (M)=\Dim (N)$ and $e(M)=e(N)$, \item $\Dim (N)=\Dim
(M)=\Dim (L)$ and $e(M)=e(N)+e(L)$.
\end{enumerate}
\end{enumerate}
\end{theorem}

{\it Proof}. $1$. By Corollary \ref{c12Jul05}, the  short exact
sequence $0\ra N\ra M\ra L\ra 0$ is obtained from a short exact
sequence of finitely generated $T_k$-modules  $0\ra N'\ra M'\ra
L'\ra 0$ by tensoring on $\CD \t_{T_k}-$. The algebra $T_k$ is
somewhat commutative with respect to the induced filtration $\CT$
from the canonical filtration $F=\{ F_i\}$ on the algebra $\CD
=\CD (P_n)$. Let $\{M_i':= \CT_{k,i}M_0\}$ be a standard
filtration on the $T_k$-module $M'$ and $\{ L_i':=\CT_{k,i}L_0\}$
be its image on $L'$ which is a standard filtration on $L'$. It is
a  well-known fact that the induced filtration $\{ N_i':=N'\cap
M_i'\}$ is good, and each good filtration is strongly equivalent
to a standard filtration. Then
$\dim_K(M_i')=\dim_K(N_i')+\dim_K(L_i')$, $i\geq 0$. Since $\CD
=\L^{[p^k]}\t T_k$ and the subalgebra $\L^{[p^k]}$ of $\CD$ has
the induced filtration $\{ \L^{[p^k]}_i:=\L^{[p^k]}\cap
F_i=\bigoplus_{p^k|\beta | \leq i} K\der^{[p^k\beta ]}\}$, it
follows that $\{ M_i:=F_iM_0=\bigoplus_{p^k|\beta | +j\leq i}
\der^{[p^k\beta ]}\t M_j'\}$ and $\{
L_i:=F_iL_0=\bigoplus_{p^k|\beta | +j\leq i} \der^{[p^k\beta ]}\t
L_j'\}$ are  filtrations of standard type on $M$ and $L$
respectively, and that $\{ N_i:=\bigoplus_{p^k|\beta | +j\leq i}
\der^{[p^k\beta ]}\t N_j'\}$ is a finite dimensional filtration on
$N$ that is strongly equivalent to a filtration of standard type
on $N$, and that $\dim_K(M_i)=\dim_K(N_i)+\dim_K(L_i)$,  $i\geq
0$. This proves statement 1.

$2$ and $3$. These statements follow from statement 1 and Theorem
\ref{23Jun05}. $\Box $


\section{Classification of simple finitely presented
$\CD (P_n)$-modules}\label{CLSFPDPn}

In this section, $K$ is an {\em arbitrary} field of characteristic
$p>0$.

In this section, a classification of {\em simple finitely
presented} $\CD (P_n)$-modules is obtained (Theorem \ref{8Jul05})
which looks particularly nice for algebraically closed fields
(Corollary \ref{c8Jul05}). It will be proved that {\em every
simple finitely presented} $\CD (P_n)$-module $M$ is {\em
holonomic}, the endomorphism algebra $\End_{\CD (P_n)}(M)$ is a
{\em finite separable field} over $K$, and the multiplicity $e(M)$
is equal to $\dim_K(\End_{\CD (P_n)}(M))$, and so it is a {\em
natural number} (Theorem \ref{8Jul05}). Plenty of holonomic $\CD
(P_n)$-modules will be considered. Some of the results of this
section are used as an inductive step in proving an analogue of
the inequality of Bernstein in Section \ref{HSSA}.

For an algebra $A$, $\widehat{A}$ denotes the set of all the
isoclasses of simple $A$-modules, and $[M]$ denotes the isoclass
of a simple $A$-module $M$.

Let  $\e =(\e_1, \ldots , \e_n)\in \{ \pm 1\}^n$, and  let $s$ be
the number of positive coordinates of $\e$. The algebra $\L_\e$ is
isomorphic to the tensor product $P_s\t \L (t)$ of the polynomial
algebra $P_s$ and $\L (t):=\L_{-1}^{\t t}$ where $t=n-s$.

The nil-radical $\Gn (\L_\e)$ of the algebra $\L_\e $ is $P_s\t \L
(t)_+$ since
 $P_s\t \L_+(t)$ belongs to the nil-radical of the
algebra $\L_\e$, and $\L_\e /(P_s\t \L_+(t))\simeq P_s\t (\L(t)/\L
(t)_+)\simeq P_s$.

\begin{lemma}\label{Lesim}
\begin{enumerate}
\item Let $\L =\L_{(-1, \ldots , -1)}$. Then $K:=\L/\L_+$ is the
only (up to isomorphism) simple $\L $-module. \item Let $\L
(t):=\L_{-1}^{\t t}$ and $\L_\e \simeq P_s\t \L (t)$ for some
$s\geq 1$ such that $ s+t=n$. Then the map $\widehat{P_s}\ra
\widehat{\L_\e}$, $[L]\mapsto [L=L\t\L (t) /\L_+(t)]$, is a
bijection.
\end{enumerate}
\end{lemma}

{\it Proof}. $1$. Note that  $\L_+$ is the nil-radical of the
algebra $\L$ and $\L/\L_+=K$. Then
$\widehat{\L}=\widehat{\L/\L_+}=\widehat{K}$, and so $K=\L/\L_+$
is the only simple $\L$-module (up to isomorphism).

$2$. Similarly,  $ \L_\e /\Gn (\L_\e)\simeq P_s$. Therefore,
$\widehat{\L_\e}=\widehat{P_s}$, and the result follows.  $\Box $

Given a ring $A$ and its element $a$, let $L_a(b)=ab$ and
$R_a(b)=ba$. Then the maps (from $A$ to itself), $L_a$, $R_a$, and
$\ad \, a=L_a-R_a$ commute. Therefore, $R_a^k=(L_a- \ad \,
a)^k=\sum_{j=0}^k{k\choose j} L_a^{k-j}(-\ad \, a)^j$, $k\geq 0$.
Applying this identity in the case where $a=x_i\in A=\CD (P_n)$,
we see that
$$ \derb x_i^k=\sum_{j=0}^{\beta_i}{k\choose j} x_i^{k-j}\der^{[\beta -je_i]}, \;\;
 \beta \in \mathbb{N}^n, \;\; k\geq 0,$$
 where $e_1:=(1, 0, \ldots , 0), \ldots , e_n:=(0,\ldots , 0, 1)$
and then, for any polynomial $f\in P_n$, 
\begin{equation}\label{cdbp}
[\derb , f]=\sum_{i=1}^n \frac{\der f}{\der x_i}\der^{[\beta
-e_i]}+\cdots =\sum_{i=1}^n \der^{[\beta -e_i]}\frac{\der f}{\der
x_i}+\cdots ,
\end{equation}
where the three dots denote an element of $\CD (P_n)_{|\beta |-2}$
where $ \{ \CD (P_n)_i\}_{i\geq 0}$ is the {\em order} filtration
on $\CD (P_n)$.

For an algebra $R$,  $R^{op}$ (or $R^o$)  stands for the {\em
opposite algebra} ($R=R^{op}$ as vector spaces but multiplication
in $R^{op}$ is given by the rule $a\cdot b=ba$).

\begin{proposition}\label{Lisim}
Let $K$ be a  field of  characteristic $p>0$, $L$ be a simple
$\L_\e$-module. Then
\begin{enumerate}
\item the induced $\CD (P_n)$-module $\CD
(P_n)\t_{\L_\e}L=\bigoplus_{\alpha \in \mathbb{N}^n} l^\alpha\t L$
is a holonomic $\CD (P_n)$-module with $Kl^\alpha \t L\simeq L$ as
$K$-modules where $\L_{-\e }= \bigoplus_{\alpha \in \mathbb{N}^n}
Kl^\alpha$. If, in addition, the field $K$ is perfect then the
$\CD (P_n)$-module $\CD (P_n)\t_{\L_\e}L$ is simple.  \item Let
$F=\{ F_i\}$ be the canonical filtration on $\CD (P_n)$ and $\{
F_iL\}$ be the filtration of standard type on the $\CD
(P_n)$-module $\CD (P_n)\t_{\L_\e}L$. Then
$\dim_K(F_iL)=\dim_K(L){i+n\choose n}$ for all $i\geq 0$. \item
If, in addition, the field $K$ is perfect then the endomorphism
algebra ${\rm End}_{\CD (P_n)}(\CD (P_n)\t_{\L_\e}L)$ is a finite
field extension over $K$ isomorphic to $K$ if $\e =(-1, \ldots ,
-1)$, and to $L'$ if $\e \neq (-1, \ldots , -1)$ where in this
case $\L_\e \simeq P_s\t \L (n-s)$, $s\geq 1$, $L=L'\t K$ (Lemma
\ref{Lesim}). \item The $\CD (P_n)$-module $\CD (P_n)\t_{\L_\e}L$
is finitely presented iff $\e = (1, \ldots , 1)$.
\end{enumerate}
\end{proposition}

{\it Proof}. It follows from the decomposition $\CD :=\CD
(P_n)=\L_{-\e }\t \L_\e=\bigoplus_{\alpha \in
\mathbb{N}^n}l^\alpha \t \L_\e$ that $M:=\CD \t_{\L_\e}L\simeq
\bigoplus_{\alpha \in \mathbb{N}^n}l^\alpha \t L$ and $Kl^\alpha
\t L\simeq L$ as $K$-modules. Then it becomes obvious that, for
$i\geq 0$,
$$ \dim_K(F_iL)=\dim_K((\L_{-\e}\cap F_i)\t
L)=\dim_K(L)\, \dim_K(\L_{-\e}\cap F_i)=\dim_K(L){i+n\choose n}.$$
This proves statement 2 and the fact that $M$ is a holonomic
$\CD$-module.

 It follows from Lemma \ref{Lesim} that the $\CD $-module
$M$ is finitely presented iff so is the $\L_\e$-module $L$ iff $\e
=(1, \ldots , 1)$. This proves statement 4.

Let us prove simplicity of $M$ in the case when the field $K$ is
perfect.  If $\e =(-1, \ldots , -1)$ then, by Lemma \ref{Lesim},
$M=P_n$ with natural action of the ring $\CD $ of differential
operators on it, and so $P_n$ is a simple $\CD $-module  with
$\End_\CD (P_n)=\cap_{\beta \in \mathbb{N}^n}\ker_{P_n}(\derb
)=K$.

It remains to consider the case when $\e \neq (-1, \ldots , -1)$.
In this case (up to order), $\L_\e =P_s\t\L (t)$ for some $s\geq
1$, $t=n-s$. By Lemma \ref{Lesim}, $L=L'\t K$ for some finite
field $L'=P_s/\Gm$ over $K$ where  $\Gm $ is a maximal ideal of
the polynomial algebra $P_s$. Now, $\CD (P_n)=\CD (P_s)\t \CD
(P_t)$, $\CD (P_s)=\L (s)\t P_s$, and $\CD (P_t)=P_t\t \L (t)$.
The $\CD (P_n)$-module $M$ is the tensor product $M_s\t M_t$ of
the $\CD (P_s)$-module $M_s:= \CD (P_s)\t_{P_s}L'\simeq \L (s)\t
L'$ and the $\CD (P_t)$-module $M_t:=\CD (P_t)\t_{\L (t)}\L (t)/\L
(t)_+=\CD (P_t)\t_{\L (t)}K\simeq P_t$. Moreover, $M\simeq \L
(s)\t L'\t P_t$. Since $M_t$ is a simple $\CD (P_t)$-module with
$\End_K(M_t)=K$, to prove the fact that $M$ is a simple $\CD
(P_n)$-module it suffices to show that $M_s$ is a simple $\CD
(P_s)$-module. For each $i=1, \ldots , s$, the kernel of the
$K$-algebra homomorphism $K[x_i]\ra P_s\ra L'=P_s/\Gm$ is
generated by an irreducible polynomial, say $p_i$. By the
assumption, $K$ is a {\em perfect} field, and so the polynomials
$p_i$ and $p_i':=\frac{d p_i}{dx_i}\neq 0$ are co-prime.
Therefore, the multiplication by $p_i'$ yields an invertible
$K$-linear map from the field $L'$ to itself. Let $u$ be a nonzero
element of $M_s$. We have to show that $\CD (P_s)u=M_s$. We use
induction on the degree $d$ of the element $u=\sum_{| \beta |
=d}\derb \t l_\beta +\sum_{| \beta' | <d}\der^{[\beta']} \t
l_{\beta'}$ where $l_\beta\in L'$ (not all are equal to zero) and
$l_{\beta'}\in L'$ where $\beta , \beta'\in \mathbb{N}^s$. The
first sum is called the {\em leading term} of the element $u$. The
case $d=0$ is obvious. So, let $d>0$. There exists $\beta$ in the
leading term of $u$ such that its $i$'th coordinate is a nonzero
one and $l_\beta \neq 0$. By (\ref{cdbp}), the element
$$ p_iu=-\sum_{|\beta | =d}\der^{[\beta -e_i]}\t p_i'l_\beta
+\cdots \neq 0,$$ has degree $<d$. Now, by induction, $\CD
(P_s)p_iu=M_s$, and so $\CD (P_s)u=M_s$, as required. This
finishes the proof of the first statement. It follows that
$\cap_{i=1}^s\ker (p_i)=L'$ in $M_s$ where $p_i:M_s\ra M_s$, $
v\mapsto p_iv$, which  implies that $\ann_{M_s}(\Gm )=L'$, but
$\ann_{M_s}(\Gm )\simeq \End_{\CD (P_s)}(M_s)^{op}$ (here we write
endomorphisms on the same side as scalars, i.e. on the left). Now,
$$\End_\CD (M)^{op}\simeq \ann_M(\Gm)\cap \ann_M(\L (t)_+)=(K\t L'\t P_t)\cap
(\L (s)\t L'\t K)=K\t L'\t K\simeq L'.$$ This proves statement 3
in the case $\e \neq (-1, \ldots , -1)$. The case $\e = (-1,
\ldots , -1)$ has been considered already. The proof of the
proposition is complete. $\Box $

For any algebraic field $L$ over $K$, let $L^{sep}$ be the {\em
maximal separable} subfield of $L$ over $K$, $L^{sep}$ is
generated by all the separable subfields of $L$ over $K$.

If the field $K$ is {\em not} necessarily perfect then the induced
module from Proposition \ref{Lisim} is not a simple module but
rather {\em semi-simple} and its endomorphism algebra is not a
field but rather a {\em direct product of matrix algebras} with
coefficients from
 {\em separable} fields (Lemmas \ref{1Lisim}, \ref{2Lisim}, and Corollary
 \ref{3Lisim}). To prove these facts, first, we consider the
simpler case when $n=1$. A simple $P_1$-module $L$ is, in fact, a
field $L=K[x]/(g)$ where $P_1=K[x]$ and $g(x):=f(x^{p^k})$ is an
irreducible polynomial such that $f(t)\in K[t]$ is an irreducible
{\em separable} polynomial ($\frac{df}{dt}\neq 0$) and $k\geq 0$.
Then $L':= K[x^{p^k}]/(g)\simeq K[t]/(f(t))$ is a finite separable
field extension of $K$, $[L':K]=\deg_t (f(t))$. Clearly,
$L'=L^{sep}$.

\begin{lemma}\label{1Lisim}
Let $K$, $L$, and $g$ be as above.
\begin{enumerate}
\item The factor  algebra $\bA :=\L_{[p^k]}\t P_1/(g)$ of the
subalgebra  $A:=T_k:=\L_{[p^k]}\t P_1$ of $\CD (P_1)$ at the
central element $g$  is isomorphic to the matrix algebra
$M_{p^k}(L')$ of rank $[L:L^{sep}]=p^k$ with coefficients from the
field $L^{sep}:=L':= K[x^{p^k}]/(g)$.
 \item $\End_{\CD (P_1)}(\CD (P_1)\t_{P_1}L)^{op}\simeq \bA$. \item The
 $\CD (P_1)$-module $\CD (P_1)\t_{P_1}L$ is a semi-simple module
 isomorphic to a direct sum of $p^k$ copies of the simple $\CD
 (P_1)$-module $U:=\CD (P_1)\t_A A/A(g, \L_{[p^k], +})=\L^{[p^k]}\t L$,  and
 $\End_{\CD (P_1)}(U)^{op}\simeq L'$.
 \item The map from the set of left ideals of the algebra
 $\bA$ to the set of
  $\CD (P_1)$-submodules of the induced module $\CD
  (P_1)\t_{P_1}L$ given by the rule $\overline{V}\mapsto \CD (P_1)\t_A
  \overline{V}$  is a bijection with
  inverse $N\mapsto N\cap \bA$.
  \item The induced $\CD (P_1)$-module $\CD (P_1)\t_{P_1}L$ is
  simple iff the polynomial $g$ is separable over $K$ (i.e. when
  $k=0$).
\end{enumerate}
\end{lemma}

{\it Remarks}. $1$. This lemma will be used as an inductive step
in Theorem \ref{27Jn05} which is a key result behind the fact that
every holonomic module has finite length (Theorem \ref{hDPnfl}).

$2$. The opposite algebra appears in statement $2$ simply because
we write endomorphisms on the {\em same side} as scalars. The
isomorphism in statement 2 is in fact an identity if one
identifies the opposite algebra of the endomorphism algebra with
the {\em idealizer} of the corresponding left ideal that defines
the cyclic module.

{\it Proof}. Let $\CD = \CD (P_1)$, $P=P_1$, and $g=f(x^{p^k})$.
Recall that $\CD =\L\t P= \L^{[p^k]} \t\L_{[p^k]}\t P=\L^{[p^k]}
\t A$ where $A$ is  a subalgebra of $\CD$, and $K[x^{p^k}]$ is the
centre of the algebra $A$. The induced $\CD$-module
$$\CD\t_PL\simeq \CD /\CD g\simeq \L^{[p^k]}\t \bA=\bigoplus_{i\geq
0}\der^{[ip^k]}\t \bA.$$

It follows from the decomposition $\bA =\L_{[p^k]}\t
P/(g)=\oplus_{0\leq i,j<p^k}\der^{[i]}x^jL'=\oplus_{0\leq
i<p^k}\der^{[i]}L$ that the algebra $\bA$ is a simple algebra with
the centre $L'$ (use $\ad\, x$ and the fact that $L$ is a field),
and $\dim_{L'}( \bA ) =p^{2k}$. In order to prove that the algebra
$\bA$ is isomorphic to the matrix algebra $M_{p^k}(L')$ it
suffices to find a simple $\bA$-module $U'$ such that
$\dim_{L'}(U')=p^k$ and $\End_{\bA}(U')\simeq L'$. One can easily
verify that the module 
\begin{equation}\label{UpAL}
U':= A/A(g, \L_{[p^k],+})\simeq \bigoplus_{0\leq i<p^k}
P\der^{[i]}/(Pg\oplus (\bigoplus_{1\leq i<p^k} P\der^{[i]}))\simeq
P/Pg\simeq L
\end{equation}
satisfy the two conditions above. This proves statement 1.

One can verify (using (\ref{CDchd1}), (\ref{cdbp}), and
separability of $f(t)$) that the $\CD$-module $\CD \t_A
U'=\L^{[p^k]}\t U'$ is a simple module. Now, the $\CD$-module $\CD
\t_P L\simeq \L^{[p^k]}\t \bA \simeq \L^{[p^k]}\t (U')^{p^k}\simeq
(\L^{[p^k]}\t U')^{p^k}$ is a direct sum of $p^k$ copies of the
simple $\CD$-module $\CD \t_A U'$. All the isomorphisms are
natural. Since the set of elements of $\CD \t_P L=\L^{[p^k]}\t
\bA$ that are  annihilated by the left ideal $A(g, \L_{[p^k], +})$
of the algebra $A$  is equal to $\bA$ and the $\CD $-module $\CD
\t_PL$ is semi-simple, statements 2 and 4 follow at once.
Statement 5 is obvious. $\Box $

For each $i=1, \ldots , n$, let $L_i:=K[x_i]/(g_i)$ be a simple
$K[x_i]$-module where $g_i(x_i)=f_i(x_i^{p^{k_i}})$ is an {\em
irreducible} polynomial such that $f_i(t)\in K[t]$ is an
 {\em irreducible separable} polynomial, $k_i\geq 0$, and
$L_i':=K[x_i^{p^{k_i}}]/(g_i)$ is a finite separable field
extension of $K$, $[L_i':K]=\deg_t(f_i)$. Clearly,
$L_i'=L_i^{sep}$.

Consider the $\CD (P_n)$-module $M:=\bigotimes_{i=1}^n M_i$  which
is the tensor product of the induced $\CD (K[x_i])$-modules
$M_i:=\CD (K[x_i])\t_{K[x_i]}L_i$. We keep the notation of Lemma
\ref{1Lisim} adding the subscript $i$ in proper places when
considering the module $M_i$. Clearly,
\begin{eqnarray*}
M&\simeq & \CD (P_n)\t_{P_n}P_n/(g_1, \ldots , g_n)= \L\t
P_n/(g_1,\ldots , g_n)=\L \t (\bigotimes_{i=1}^n
L_i)\\
&=&\bigoplus_{\alpha \in \mathbb{N}^n}\der^{[\alpha ]}\t
(\bigotimes_{i=1}^nL_i)\simeq \CD (P_n)/\CD (P_n)(g_1,\ldots
,g_n).
\end{eqnarray*}
 $\{ \CM_i:=F_i \cdot P_n/(g_1, \ldots , g_n)\}$ is the  filtration of standard type on the
$\CD (P_n)$-module $M$. Then
$$\dim_K(\CM_i)=\prod_{j=1}^n[L_j:K] {i+n\choose
n}=p^k\prod_{i=1}^n \deg_t(f_i(t)) {i+n\choose n}, \;\; i\geq 0,$$
where $k:=k_1+\cdots +k_n$. So, $M$ is a {\em holonomic cyclic
finitely presented} $\CD (P_n)$-module.

By Lemma \ref{1Lisim}, $\End_{\CD (K[x_i])}(M_i)^{op}\simeq
\bA_i\simeq M_{p^{k_i}}(L_i')$. It follows that
$$
 \End_{\CD (P_n)}(M)^{op}\simeq
\bigotimes_{i=1}^n\bA_i \simeq
\bigotimes_{i=1}^nM_{p^{k_i}}(L_i')\simeq
M_{p^k}(\bigotimes_{i=1}^nL_i')\\
\simeq  M_{p^k}(\prod_{\nu =1}^\mu \G_\nu)\simeq \prod_{\nu
=1}^\mu M_{p^k}(\G_\nu ). $$ The tensor product of separable
fields $\t_{i=1}^nL_i'$ is a {\em semi-simple commutative}
algebra, it is a direct product $\prod_{\nu =1}^\mu\G_\nu$ of {\em
finite separable fields} $\G_\nu$ over $K$. The algebra
$\bA:=\otimes_{i=1}^n\bA_i$ is a semi-simple finite dimensional
algebra. Let $V_\nu$, $\nu =1,\ldots , \mu$, be a complete set of
(pairwise non-isomorphic) simple  $\bA$-modules. Then
$\dim_K(V_\nu )=p^k[\G_\nu :K]$ and $\End_A(V_\nu )^{op}\simeq
\G_\nu$.

It follows from the equality $\CD
(P_n)=(\bigotimes_{i=1}^n\L_i^{[p^{k_i}]})\t A$ where
$A:=\otimes_{i=1}^nA_i$, $A_i:=\L_{i,[p^{k_i}]}\t K[x_i]$, that
the $\CD (P_n)$-module $M\simeq \CD (P_n)\t_A\bA\simeq
(\bigoplus_{\nu =1}^\mu \CD (P_n)\t_A V_\nu)^{p^k}$ is a direct
sum of {\em simple} $\CD (P_n)$-modules $U_\nu :=\CD
(P_n)\t_AV_\nu$, and each of them occurs with the {\em same}
multiplicity $p^k$. Summarizing, we have the following lemma which
 is a direct consequence of Lemma \ref{1Lisim}.

 \begin{lemma}\label{2Lisim}
 Let $K$ be an arbitrary field of characteristic $p>0$, the $\CD
 (P_n)$-module $M=\otimes_{i=1}^nM_i$ be the tensor product of
 modules from Lemma \ref{1Lisim}. Then
 \begin{enumerate}
\item The algebra $\bA:=\bigotimes_{i=1}^n\bA_i\simeq \prod_{\nu
=1}^\mu M_{p^k}(\G_\nu)$ where $k:=k_1+\cdots +k_n$ and $\G_\nu$
are finite separable field extensions of $K$. \item $\End_{\CD
(P_n)}(M)^{op}\simeq \bA$. \item The $\CD (P_n)$-module $M$ is a
semi-simple holonomic cyclic finitely presented  module isomorphic
to the direct sum $\oplus_{\nu =1}^\mu U_\nu^{p^k}$ where $U_\nu
:=\CD (P_n)\t_AV_\nu$ is a simple holonomic finitely presented
$\CD (P_n)$-module, and $\End_{\CD (P_n)}(U_\nu )^{op}\simeq
\G_\nu$ is a finite separable field extension of $K$. \item On the
simple $\CD (P_n)$-module
$U_\nu=(\bigotimes_{i=1}^n\L_i^{p^{k_i}})\t V_\nu$ consider the
filtration of standard type $\{ U_{\nu , i}:=F_i\cdot 1\t V_\nu
=\bigoplus_{i_1p^{k_1}+\cdots +i_np^{k_n}\leq
i}\der_1^{[i_1p^{k_1}]}\cdots \der_n^{[i_np^{k_n}]}\t V_\nu\}$.
Then
\begin{enumerate}
\item the Poincare series $P_{U_\nu}=\frac{\dim_K(V_\nu)}{(1-\o
)\prod_{i=1}^n(1-\o^{p^{k_i}})}=\frac{p^k[\G_\nu :K]}{(1-\o
)\prod_{i=1}^n(1-\o^{p^{k_i}})}$ where $k:=k_1+\cdots +k_n$, \item
the multiplicity $e(U_\nu )=[\G_\nu :K]=\dim_K(\End_{\CD
(P_n)}(U_\nu ))$, \item $\dim_K(U_{\nu , i})=\frac{e(U_\nu
)}{n!}i^n+\cdots $, $i\gg 0$, is an almost polynomial with period
$p^{\max \{k_1, \ldots , k_n \}}$.
\end{enumerate}
 \item The map from the
set of left ideals of the algebra $\bA$ to the set of $\CD
(P_n)$-submodules of $M$ given by the rule $\overline{V} \mapsto
\CD (P_n)\t_A \overline{V}$ is a bijection with inverse $N\mapsto
N\cap \bA$. \item The $\CD (P_n)$-module $M$ is simple iff all the
polynomials $g_1, \ldots , g_n$ are separable (i.e. $k_1=\cdots
=k_n=0$) and the tensor product of fields $\otimes_{i=1}^n
L_i'=\otimes_{i=1}^nL_i$ is a field.
 \end{enumerate}
\end{lemma}

\begin{corollary}\label{3Lisim}
Let $K$ be an arbitrary field of characteristic $p>0$, $L$ be a
simple $\L_\e$-module. Then the induced $\CD (P_n)$-module $\CD
(P_n)\t_{\L_\e}L$ is a semi-simple holonomic $\CD (P_n)$-module of
finite length and $\End_{\CD (P_n)}(\CD (P_n)\t_{\L_\e}L)\simeq
\prod_{\nu =1}^\mu M_{n_\nu }(\G_\nu )$ where $\G_\nu$ are finite
separable field extensions of $K$,  $n_\mu\geq 0$, $M_0(\G_\nu
):=K$ (see the proof).
\end{corollary}

{\it Proof}. We keep the notation of Lemma \ref{Lisim} and its
proof. The case when $\e =(-1, \ldots ,-1)$ has been considered
already in the proof of Lemma \ref{Lisim} (in this case, $M=P_n$
and $\End_{\CD (P_n)}(P_n)\simeq K$).

So, we may  assume that $\e \neq (-1, \ldots ,-1)$. In this case,
$\L_\e =P_s\t \L (t)$ for some $s\geq 1$, $t=n-s$, and the $\CD
(P_n)$-module $M=M_s\t M_t$ (see the proof of Lemma \ref{Lisim})
where the $\CD (P_t)$-module $M_t$ is equal to $P_t$ and the $\CD
(P_s)$-module $M_s$ is an epimorphic image of a $\CD (P_s)$-module
$M=\otimes_{i=1}^sM_i$ from Lemma \ref{2Lisim}. Since $P_t$ is a
simple $\CD (P_t)$-module with $\End_{\CD (P_t)}(P_t)=K$, every
$\CD (P_s)\t \CD (P_t)$-submodule of $M_s\t M_t$ is equal to $N\t
M_t$ for some $\CD (P_s)$-submodule $N$ of $M_s$. By Lemma
\ref{2Lisim}, $M\simeq \oplus_{\nu =1}^\mu U_\nu^{n_\nu}\t M_t$
and $\End_{\CD (P_n)}(M)\simeq \prod_{\nu =1}^\mu M_{n_\nu
}(\G_\nu )$ for some $n_\nu \geq 0$ such that $0\leq n_\nu (\leq
p^k)$ where $M_0(\G_\nu ):=K$. $\Box $

Let $\Max (P_n)$ be be the set of all the maximal ideals of the
polynomial algebra $P_n$. Let $\Gm\in \Max (P_n)$, we are going to
determine the structure of the induced $\CD (P_n)$-module $\CD
(P_n)\t_{P_n}P_n/\Gm$ (Lemma \ref{1l8Jul05}, this lemma is central
in proving Theorem \ref{8Jul05}). Note the $P_n/\Gm$ is a finite
field over $K$. For each $i=1, \ldots , n$, there exists a unique
monic irreducible polynomial $g_i\in K[x_i]$ such that
$(g_i)=K[x_i]\cap \Gm$, then $g_i(x_i)=f_i(x_i^{p^{k_i}})$ where
$f_i(t)\in K[t]$ is a monic irreducible separable polynomial for
some $k_i\geq 0$. Note that $g_i$, $f_i$, and $k_i$ are uniquely
determined by the ideal $\Gm$. Let $k(\Gm )=(k_1, \ldots , k_n)$,
$g(\Gm )=(g_1, \ldots , g_n)$, and $f(\Gm )=(f_1, \ldots , f_n)$.
Let $L_i:=K[x_i]/(g_i)$, $L_i':=L_i^{sep}=K[x_i^{p^{k_i}}]/(g_i)$,
$\bigotimes_{i=1}^nL_i'\simeq \prod_{\nu =1}^\mu \G_\nu$ where
$\G_\nu$ are finite separable fields over $K$, let $1=\sum_{\nu
=1}^\mu e_\nu$ be the corresponding sum of primitive orthogonal
idempotents. Let  $A(\Gm ):=\bigotimes_{i=1}^n A(\Gm )_i$ where
$A(\Gm )_i:=\L_{i, [p^{k_i}]}\t K[x_i]$, $\bA (\Gm )_i:=A(\Gm
)_i/A(\Gm )_ig_i$,
$$ \bA (\Gm ):=\bigotimes_{i=1}^n \bA (\Gm )_i\simeq
M_{p^k}(\bigotimes_{i=1}^nL_i')\simeq \prod_{\nu =1}^\mu
M_{p^k}(\G_\nu ), \;\;\;\; \L (\Gm ):= \bigotimes_{i=1}^n
\L_i^{[p^{k_i}]}, $$ where $k:=k_1+\cdots +k_n$.

Let us consider the map $\prod_{\nu =1}^\mu \G_\nu \simeq
\bigotimes_{i=1}^nL_i'\ra P_n/\Gm$ that is the composition of the
inclusion $\bigotimes_{i=1}^nL_i'\ra \bigotimes_{i=1}^nL_i$ and
the natural algebra epimorphism $\bigotimes_{i=1}^nL_i\ra
P_n/\Gm$. Then there exists a {\em unique} $\nu$ such that the map
$\G_\nu \ra P_n/\Gm$ $(e_\nu \mapsto 1)$ is a $K$-algebra
monomorphism. We denote such a unique field $\G_\nu $ by $\G
(\Gm)$. It is obvious that 
\begin{equation}\label{GmPnms}
\G (\Gm )=(P_n/\Gm )^{sep}
\end{equation}
since $e_\mu \mapsto 0$, if $\mu \neq \nu$,
$\bigotimes_{i=1}^nL_i\ra P_n/\Gm$ is an epimorphism, and the
$p^j$'th ($j\gg 1$) power of each element of
$\bigotimes_{i=1}^nL_i$ belongs to $\bigotimes_{i=1}^nL_i'$. The
module
$$ U(\Gm ):=\CD (P_n)\t_{A(\Gm )} V(\Gm )$$
is a   {\em simple holonomic finitely presented} $\CD
(P_n)$-module $U_\nu$ from Lemma \ref{2Lisim} that corresponds to
the field $\G (\Gm )=\G_\nu$ where $V(\Gm ):=V_\nu$.

\begin{lemma}\label{1l8Jul05}
Keep the notation as above. For each maximal ideal $\Gm$ of the
polynomial algebra $P_n$, the induced $\CD (P_n)$-module $\CD
(P_n)\t_{P_n}P_n/\Gm$ is isomorphic to $\frac{[P_n/\Gm
:K]}{[(P_n/\Gm )^{sep} :K]}$ copies of the simple holonomic
finitely presented $\CD (P_n)$-module $U(\Gm )$. In particular,
the $\CD (P_n)$-module $\CD (P_n)\t_{P_n}P_n/\Gm$ is simple iff
the field $P_n/\Gm$ is separable.
\end{lemma}

{\it Proof}. Applying $\CD (P_n)\t_{A(\Gm )}-$ to the natural
epimorphism of $A(\Gm )$-modules $\bA (\Gm )\ra A(\Gm
)\t_{P_n}P_n/\Gm$, we have the natural epimorphism of $\CD
(P_n)$-modules
$$ \CD (P_n)\t_{A(\Gm )}\bA (\Gm )\simeq \CD (P_n)\t_{A(\Gm
)}\prod_{\nu =1}^\mu M_{p^k}(\G_\nu )\ra \CD (P_n)\t_{A(\Gm
)}A(\Gm )\t_{P_n}P_n/\Gm\simeq \CD (P_n)\t_{P_n} P_n/\Gm.$$ Since
$\G_\mu \ra 0$, if $\G_\mu \neq \G_\nu$, we have the natural
epimorphism of $\CD (P_n)$-modules
$$ \CD (P_n)\t_{A(\Gm )}M_{p^k}(\G (\Gm ))\simeq U(\Gm )^{p^k}\ra
\CD (P_n)\t_{P_n}P_n/\Gm .$$ Therefore,  $\CD (P_n)\t_{P_n}P_n/\Gm
\simeq U(\Gm )^s$ for some $s\geq 1$. On the module $\CD
(P_n)\t_{P_n}P_n/\Gm$ consider the filtration of standard type $\{
F_i 1\t P_n/\Gm =\bigoplus_{|\beta | \leq i}\der^{[\beta ]}\t
P_n/\Gm \}$. Then $$\dim_K(F_i 1\t P_n/\Gm )=[P_n/\Gm
:K]{i+n\choose n}=\frac{[P_n/\Gm :K]}{n!}i^n+\cdots , \;\; i\gg
0.$$ By Lemma \ref{2Lisim}, $\dim_K(U(\Gm )_i^s)=\frac{s[\G(\Gm
):K]}{n!}i^n+\cdots $, $i\gg 0$. Since the multiplicity does not
depend on a filtration of standard type, we must have $s[\G (\Gm
):K]=[P_n/\Gm :K]$. This finishes the proof of the  lemma (see
(\ref{GmPnms})). $\Box $

Let $\widehat{\CD (P_n)}({\rm fin. \; pres. })$ be the set of
isoclasses of simple finitely presented $\CD (P_n)$-modules.
Theorem \ref{8Jul05} {\em classifies} these modules and shows that
{\em every} simple finitely presented $\CD (P_n)$-module is
 {\em holonomic}.

\begin{theorem}\label{8Jul05}
Let $K$ be a field of characteristic $p>0$. Then
\begin{enumerate}
\item The map $\Max (P_n)\ra \widehat{\CD (P_n)}({\rm fin. \;
pres. })$, $\Gm \mapsto [U(\Gm ):=\CD (P_n)\t_{A(\Gm )}V(\Gm )]$,
is a bijection with inverse $[M]\mapsto \ass_{P_n}(M)$ (the set of
all associated primes for the $P_n$-module $M$). In particular,
$\ass_{P_n}(U(\Gm ))=\{ \Gm \}$.
 \item Each
simple finitely presented $\CD (P_n)$-module $M$ is a holonomic.
\item (An analogue of Quillen's Lemma). $\End_{\CD (P_n)}(U(\Gm
))\simeq (P_n/\Gm )^{sep}$.\item On the simple $\CD (P_n)$-module
$U (\Gm )=\L (\Gm ) \t V(\Gm)$ consider the filtration of standard
type $\{ U (\Gm )_i:=F_i1\t V(\Gm ) =\bigoplus_{i_1p^{k_1}+\cdots
+i_np^{k_n}\leq i}\der_1^{[i_1p^{k_1}]}\cdots
\der_n^{[i_np^{k_n}]}\t V(\Gm )\}$. Then
\begin{enumerate}
\item the Poincare series $P_{U(\Gm )}=\frac{p^k[(P_n/\Gm
)^{sep}:K]}{(1-\o )\prod_{i=1}^n(1-\o^{p^{k_i}})}$, $k:=k_1+\cdots
+k_n$,  \item the multiplicity $e(U(\Gm ))=[(P_n/\Gm
)^{sep}:K]=\dim_K(\End_{\CD (P_n)}(U(\Gm )))$ is a natural number,
\item $\dim_K(U(\Gm )_i)=\frac{[(P_n/\Gm )^{sep}:K]
)}{n!}i^n+\cdots $, $i\gg 0$, is an almost polynomial with period
$p^{\max \{k_1, \ldots , k_n \}}$.
\end{enumerate}
\end{enumerate}
\end{theorem}

{\it Remark}. $\dim_K(U(\Gm )_i)$ is {\em not} a polynomial (for
$i\gg 0$) iff $\max \{k_1, \ldots , k_n \}>1$.

 {\it Proof}. $1$. Let $M$ be a simple finitely presented
$\CD (P_n)$-module. By Corollary \ref{1c23Jun05} and its proof,
$M\simeq  \CD (P_n)\t_{T_k}M'$ is a holonomic $\CD (P_n)$-module
where $M'$ is a simple finite dimensional $T_k$-module. Then $M'$
is a finite dimensional $P_n$-module as $P_n\subseteq T_k$. Then
the $P_n$-module $M'$ contains a simple $P_n$-module isomorphic to
 $P_n/\Gm$ where $\Gm$ is a maximal ideal of the algebra $P_n$.
Then $M$ is an epimorphic image of the $\CD (P_n)$-module $N:=\CD
(P_n)\t_{P_n}P_n/\Gm$, and so $M\simeq U(\Gm )$, by Lemma
\ref{1l8Jul05}. Note that $N=\cup_{i\geq 1}\ann (\Gm^i)$, and so
$\{ \Gm \} =\ass_{P_n}(N)=\ass_{P_n}(U(\Gm )^s)=\ass_{P_n}(U(\Gm
))$. Therefore, the map $\Gm \mapsto U(\Gm )$ is a bijection with
inverse $M\mapsto \ass_{P_n}(M)$. Statements 2--4 follow from
statement 1 and Lemma \ref{2Lisim}. $\Box $

\begin{corollary}\label{c8Jul05}
Let $K$ be an algebraically closed field of characteristic $p>0$.
Then
\begin{enumerate}
\item The map $\Max (P_n)=K^n\ra \widehat{\CD (P_n)}({\rm fin. \;
pres. })$, $\Gm \mapsto [U(\Gm ):=\CD (P_n)\t_{P_n}P_n/\Gm ]$, is
a bijection with inverse $[M]\mapsto \ass_{P_n}(M)$. In
particular, $\ass_{P_n}(U(\Gm ))=\{ \Gm \}$.\item  $\End_{\CD
(P_n)}(U(\Gm ))\simeq K$.\item On the simple $\CD (P_n)$-module $U
(\Gm )=\L \t P_n/\Gm =\L \overline{1}$ consider the filtration of
standard type $\{ U (\Gm )_i:=F_i\overline{1} =\bigoplus_{|\beta
|\leq i}K\der^{[\beta]}\overline{1}\}$. Then
\begin{enumerate}
\item the Poincare series $P_{U (\Gm )}=\frac{1}{(1-\o )^{n+1}}$,
\item the multiplicity $e(U (\Gm ))=1$, \item $\dim_K(U (\Gm
)_i)={ i+n\choose n}$ is a  polynomial.
\end{enumerate}
\end{enumerate}
\end{corollary}


\section{Classification of tiny simple (non-finitely presented) $\CD
(P_n)$-modules}\label{CTINY}

In this section, $K$ is an {\em arbitrary} field of characteristic
$p>0$.

In this section, we complete a classification of the `smallest'
simple $\CD (P_n)$-modules (see Theorems \ref{10Jul05} and
\ref{8Jul05}), they are called tiny modules. Theorem \ref{8Jul05}
describes the set of tiny {\em finitely presented} $\CD
(P_n)$-modules and Theorem \ref{10Jul05} classifies the set of
tiny {\em non-finitely presented} $\CD (P_n)$-modules. They turned
out to be {\em holonomic} with multiplicities which are natural
numbers. Briefly, they have the same properties as simple finitely
presented $\CD (P_n)$-modules.

Let $\e \in \{\pm 1\}^n$. A $\L_\e$-module $M$ is called a {\em
locally finite} if $\dim_K(\L_\e m)<\infty$ for each element $m\in
M$. We denote by $\CL_\e$ the category of all $\CD (P_n)$-modules
that are locally finite as $\L_\e$-modules. The category $\CL_\e$
 is a full subcategory of the category $\CD (P_n)$-modules (it is
closed under taking sub/factor modules and  direct sums but not
under infinite  direct products).

{\it Definition}. A simple $\CD (P_n)$-module from $\CL_\e $ is
called a {\em tiny} module. The name is inspired by Theorem
\ref{27Jn05} (which roughly speaking says that `typically'
$\dim_K(\L_\e m)=\infty$).

Our aim is to describe the set $\widehat{\CD (P_n)}(\CL_\e )$ of
all the isoclasses of simple $\CD (P_n)$-modules that are locally
finite over $\L_\e$ (Theorem \ref{10Jul05} and Corollary
\ref{c10Jul05}).

For each $\Gm \in \Max (\L_\e )$, the $\CD (P_n)$-module $\CD
(P_n)\t_{\L_\e}\L_\e /\Gm =\L_{-\e}\t \L_\e /\Gm
=\bigoplus_{\alpha \in \mathbb{N}^n}l_{-\e}^\alpha \t \L_\e /\Gm $
is a holonomic $\CD (P_n)$-module as the filtration of standard
type $\{ F_i1\t \L_\e /\Gm =\bigoplus_{|\alpha |\leq
i}l_{-\e}^\alpha \t \L_\e /\Gm \}$ on it  has polynomial growth
$$ \dim_K( F_i1\t \L_\e /\Gm )= \dim_K(\L_\e /\Gm ){i+n\choose
n}, \;\; i\geq 0.$$ For each $j\geq 1$,  $\ass_{\L_\e}(\CD
(P_n)\t_{\L_\e}\L_\e /\Gm^j)=\{ \Gm\}$ and $\CD
(P_n)\t_{\L_\e}\L_\e /\Gm^j\in \CL_\e$. If $\e \neq ( 1, \ldots ,
1)$ then, for each $j\geq 2$, the cyclic $\CD (P_n)$-module $\CD
(P_n)\t_{\L_\e}\L_\e /\Gm^j$ is not Noetherian as $\dim_K(\Gm /
\Gm^2)=\infty$.  It follows that each module $M\in \CL_\e$ is an
epimorphic image of a direct sum of induced modules of the type
$\CD (P_n)\t_{\L_\e}\L_\e /\Gm^j$, and that
$$ M=\bigoplus_{\Gm \in \Max (\L_\e )}M^\Gm$$
is a direct sum of uniquely determined $\CD (P_n)$-submodules
$M^\Gm :=\cup_{i\geq 1}\ann_M (\Gm^i)$ with $\ass_{\L_\e} (M^\Gm
)=\{ \Gm \}$. Therefore, each {\em simple} module from the
category $\CL_\e$ is an epimorphic image of a module of type $\CD
(P_n)\t_{\L_\e}\L_\e /\Gm$.

{\it Example}. For $\e =(1, \ldots , 1)$, i.e. $\L_\e =P_n$, we
have already got the description of the set $\widehat{\CD
(P_n)}(\CL_\e )=\widehat{\CD (P_n)}({\rm fin. \; pres.})$ (Theorem
\ref{8Jul05}).

\begin{theorem}\label{10Jul05}
Let $K$ be a field of characteristic $p>0$ and $\L_\e =\L (t)\t
P_s$ where $t\geq 1$ and $s:=n-t$ (i.e. $\e \neq (1, \ldots ,
1)$). Then $\CD (P_n)=\CD (P_t)\t \CD (P_s)$ and
\begin{enumerate}
\item The map $\Max (P_s)\ra \widehat{\CD (P_n)}(\CL_\e )$, $\Gm
\mapsto \CU(\Gm ):=P_t\t U(\Gm )$, is a bijection with inverse
$M\mapsto \ass_{P_s}(M)$. In particular, $\ass_{P_s}(\CU (\Gm
))=\{ \Gm \}$. \item The map $\Max (\L_\e)\ra \widehat{\CD
(P_n)}(\CL_\e )$, $\L (t)_+\t P_s+\L (t)\t \Gm \mapsto \CU(\Gm )$,
is a bijection with inverse $M\mapsto \ass_{\L_\e}(M)$. In
particular, $\ass_{\L_\e}(\CU (\Gm ))=\{ \L (t)_+\t P_s+\L (t)\t
\Gm\}$.
 \item Each
simple  $\CD (P_n)$-module from $\widehat{\CD (P_n)}(\CL_\e )$ is
a holonomic but not finitely presented. \item (An analogue of
Quillen's Lemma). $\End_{\CD (P_n)}(\CU (\Gm ))\simeq\End_{\CD
(P_t)}(P_t)\t \End_{\CD (P_s)}(U(\Gm ))\simeq K\t (P_s/\Gm
)^{sep}\simeq (P_s/\Gm )^{sep}$.\item On the simple $\CD
(P_n)$-module $\CU (\Gm )=P_t\t \L (\Gm ) \t V(\Gm)$ consider the
filtration of standard type $\{ \CU (\Gm )_i:=F_i1\t1\t V(\Gm )
=\bigoplus_{\alpha \in \mathbb{N}^t, |\alpha | +i_1p^{k_1}+\cdots
+i_sp^{k_s}\leq i}x^\alpha \der_1^{[i_1p^{k_1}]}\cdots
\der_s^{[i_sp^{k_s}]}\t V(\Gm )\}$. Then
\begin{enumerate}
\item the Poincare series $P_{\CU (\Gm )}=(1-\o ) P_{P_t}P_{U(\Gm
)}=\frac{p^k[(P_s/\Gm )^{sep}:K]}{(1-\o
)^{t+1}\prod_{i=1}^s(1-\o^{p^{k_i}})}$ where $k=k_1+\cdots +k_s$,
\item the multiplicity $e(\CU (\Gm ))=e(P_t)e(U(\Gm ))=[(P_s/\Gm
)^{sep}:K] $ and $e(\CU (\Gm ))=\dim_K(\End_{\CD (P_n)}(U(\Gm
)))$, \item $\dim_K(\CU (\Gm )_i)=\frac{[(P_s/\Gm )^{sep}:K]
)}{n!}i^n+\cdots $, $i\gg 0$, is an almost polynomial with period
$p^{\max \{k_1, \ldots , k_s \}}$.
\end{enumerate}
\end{enumerate}
\end{theorem}

{\it Remark}. $\dim_K(\CU (\Gm )_i)$ is {\em not} a polynomial
(for $i\gg 0$) iff $\max \{k_1, \ldots , k_s \}>1$.

 {\it Proof}. Note that the map $\Max (P_s)\ra \Max (\L_\e )$,
 $\Gm \mapsto \L (t)_+\t P_s+\L (t)\t \Gm $, is a bijection. It
 follows that  $P_t=\CD (P_t)/\CD (P_t)\L (t)_+$ is  a simple
 (non-finitely presented) $\CD (P_t)$-module with $\End_{\CD
 (P_t)}(P_t)=K$, and that any simple  $\CD (P_n)$ module $M$ from
 $\widehat{\CD (P_n)}(\CL_\e )$ such that $\ass_{P_s}(M)=\{ \Gm
 \}$ is an epimorphic image of the $\CD (P_n)$-module $P_t\t (\CD
 (P_s)\t_{P_s}P_s/\Gm )$. Therefore, $M\simeq P_t\t U(\Gm )$
 (Lemma \ref{1l8Jul05}). Now, the results follow from Theorem
 \ref{8Jul05}.  $\Box $

\begin{corollary}\label{c10Jul05}
Keep the notation from Theorem \ref{10Jul05}. If, in addition, the
field $K$ is algebraically closed then
\begin{enumerate}
\item The map $\Max (P_s)=K^s\ra \widehat{\CD (P_n)}(\CL_\e )$,
$\Gm \mapsto \CU (\Gm ):= P_t\t (\CD (P_s)\t_{P_s}P_s/\Gm )$, is a
bijection with inverse $M\mapsto \ass_{P_s}(M)$. In particular,
$\ass_{P_s}( \CU (\Gm ))=\{ \Gm \}$. \item The map $\Max
(\L_\e)\ra \widehat{\CD (P_n)}(\CL_\e )$, $\L (t)_+\t P_s+\L (t)\t
\Gm \mapsto \CU(\Gm )$, is a bijection with inverse $M\mapsto
\ass_{\L_\e}(M)$. In particular, $\ass_{\L_\e}(\CU (\Gm ))=\{ \L
(t)_+\t P_s+\L (t)\t \Gm\}$. \item  $\End_{\CD (P_n)}(\CU(\Gm
))\simeq K$.\item On the simple $\CD (P_n)$-module $\CU (\Gm
)=P_t\t \L (s)\t P_s/\Gm =P_t\t \L (s)\overline{1}$, consider the
filtration of standard type $\{ \CU
 (\Gm )_i:=F_i\overline{1} =\bigoplus_{\alpha \in \mathbb{N}^t, \beta \in \mathbb{N}^s, |\alpha  | +|\beta | \leq
i}Kx^\alpha\der^{[\beta]}\overline{1}\}$. Then
\begin{enumerate}
\item the Poincare series $P_{U(\Gm )}=\frac{1}{(1-\o )^{n+1}}$,
\item the multiplicity $e(U(\Gm ))=1$, \item $\dim_K(U(\Gm )_i)={
i+n\choose n}$ is a  polynomial.
\end{enumerate}
\end{enumerate}
\end{corollary}


\section{Multiplicity of each finitely presented $\CD (P_n)$-module
is a natural number}\label{MEFPNN}

In this section, $K$ is an {\em arbitrary} field of characteristic
$p>0$.

We know already that the multiplicity of a {\em non-holonomic
finitely presented} $\CD (P_n)$-module can be {\em arbitrary
small} (Lemma \ref{a22Jul05}). In this section, we prove that the
multiplicity of a {\em holonomic finitely presented} $\CD
(P_n)$-module is a {\em natural number} (Theorem \ref{t12Jul05}).
This result is a direct consequence of a classification of simple
$\Tk $-modules (Theorem \ref{20Jul05})  and Theorem \ref{23Jun05}.

For each $\bk =(k_1, \ldots , k_n)\in \mathbb{N}^n$, the
subalgebra of $\CD (P_n)=\bigotimes_{i=1}^n\CD (K[x_i])$:
$$ \Tk =T_{\bk, n}:=\bigotimes_{i=1}^n(\L_{i,[p^{k_i}]}\t
K[x_i])=\L_{[\p^{\bk } ]}\t P_n=\bigoplus_{\beta <\p^\bk
}\der^{[\beta ]}\t P_n=\bigoplus_{\beta <\p^\bk}P_n\t \der^{[\beta
]}$$ is a free left and right $P_n$-module of rank $p^{|\bk |}$
where $\L_{[\p^{\bk } ]}:= \bigotimes_{i=1}^n\L_{i, [p^{k_i}]} $,
$|\bk | :=k_1+\cdots +k_n$, and $\beta <\p^\bk$ means
$\beta_i<p^{k_i}$ for all $i$. It is a finitely generated
Noetherian algebra with the centre $\Zk :=K[x_1^{p^{k_1}}, \ldots
,x_n^{p^{k_n}}]$. The algebra $\Tk$ is a free $\Zk$-module of rank
$p^{2|\bk |}$ since $\Tk =\L_{[\p^\bk ]}\t (\oplus_{\alpha <\p^\bk
}Kx^\alpha )\t \Zk$. On the algebra $\Tk$ consider the induced
filtration from the canonical filtration $F=\{ F_i\}$ on the
algebra $\CD (P_n)$: 
\begin{equation}\label{1FCTk}
 \CT_\bk=\{ \CT_{\bk,i}:=\Tk\cap F_i=\bigoplus_{ \beta
<\p^\bk , |\alpha | +|\beta |\leq i} Kx^\alpha
\der^{[\beta]}=\bigoplus_{ \beta <\p^\bk , |\alpha | +|\beta |\leq
i} K \der^{[\beta]}x^\alpha \}.
\end{equation}
The filtration $\CT_\bk$ is the tensor product of the induced
filtrations on each tensor multiple $\L_{i, [p^{k_i}]}\t K[x_i]$
of the algebra $\Tk$. The associated graded algebra $\gr \, \Tk
=\oplus_{i\geq 0}G_{\bk , i}$ is naturally isomorphic (as a graded
algebra) to the {\em tensor product of the commutative algebras}
$\L_{[\p^\bk]}\t P_n$ where
$$ G_{\bk , i}:=\bigoplus_{\beta <\p^\bk , | \alpha |
+|\beta | =i} K \der^{[\beta ]}x^\alpha.$$
 The grading on $\gr \, \Tk$ is the tensor product of natural
 gradings on the tensor multiples. The algebra $\gr\, \Tk$ is an
 {\em affine commutative} algebra with nil-radical $\L_{[\p^\bk ],
 +}\t P_n$ (where $\L_{[\p^\bk ],
 +}:= \oplus_{0\neq \beta <\p^\bk}K\der^{[\beta ]}$) which is a
 prime ideal since
 $$ \gr \, \Tk /(\L_{[\p^\bk ],  +}\t P_n)\simeq (\L_{[\p^\bk ]}/\L_{[\p^\bk ],
 +})\t P_n\simeq K\t P_n\simeq P_n.$$
$T_0:=T_{(0, \ldots , 0)}=P_n$, $\Tk\subseteq  \Tl$ iff $\bk \leq
{\bf l}$ (i.e. $k_1\leq l_1, \ldots , k_n\leq l_n)$. $\CD
(P_n)=\cup_{\bk \in \mathbb{N}^n} \Tk$, $\Tk \Tl \subseteq T_{\max
(\bk , {\bf l} )}$ where $\max (\bk , {\bf l}):= (\max (k_1 ,
l_1), \ldots ,\max (k_n , l_n))$.

\begin{lemma}\label{1Tkpr18Jul05}
\begin{enumerate}
\item  The algebra $\Tk$ is a somewhat commutative algebra with
respect to the finite dimensional filtration $\CT_\bk =\{\CT_{\bk
,i}\}$ having the centre $\Zk=K[x_1^{p^{k_1}}, \ldots
,x_n^{p^{k_n}}]$ and $\GK (\Tk )=n$. In particular, $\Tk $ is a
finitely generated Noetherian algebra. \item The Poincare series
of $\Tk$, $P_{\Tk}=\sum_{i\geq 0}\dim_K(\CT_{\bk
,i})\o^i=\frac{\prod_{i=1}^n(1+\o +\o^2+\cdots +
\o^{p^{k_i}-1})}{(1-\o )^{n+1}}$ and the multiplicity $e(\Tk
)=p^{|\bk |}$. \item The Hilbert function is, in fact, a
polynomial $\dim_K(\CT_{\bk ,i})=\frac{p^{|\bk |}}{n!}i^n+\cdots
$, $i\gg 0$. \item Let $\CZ_\bk =K(x_1^{p^{k_1}}, \ldots
,x_n^{p^{k_n}})$ be the field of fractions of $\Zk$. Then
$\Tk':=\CZ_\bk \t_{\Zk}\Tk\simeq M_{p^{|\bk |}}(\CZ_\bk )$, the
matrix algebra. \item The algebra $\Tk$ is a prime algebra of
uniform dimension $p^{|\bk |}$, and the localization $\CS^{-1}\Tk$
of $\Tk$ at the set $\CS$ of all the non-zero divisors is
isomorphic to the matrix algebra $M_{p^{|\bk |}}(\CZ_\bk )$. \item
The algebra $\Tk$ is preserved by the involution $*$, $\Tk^*=\Tk$,
and so the algebra $\Tk$ is self-dual. \item The algebra $\Tk$ is
faithfully flat over its centre. \item The left and right Krull
dimension of the algebra $\Tk$ is $n$. \item The left and right
global dimension of the algebra $\Tk$ is $n$ but the global
dimension of the associated graded algebra $\gr (\Tk )$  is
$\infty$ if $\bk \neq 0$.
\end{enumerate}
\end{lemma}

{\it Proof}. Repeat the proof of Lemma \ref{Tkpr}. $\Box $

Recall that the algebra $T_\bk=T_{\bk ,n}$ is a somewhat
commutative algebra with respect to the filtration $ \CT_\bk$.

\begin{lemma}\label{l18Jul05}
Let $M$ be a finitely generated $\Tk$-module, $\bk \leq {\bf l}$,
and $M'=\Tl \t_{\Tk }M$ be a $\Tl$-module. Then
\begin{enumerate}
\item $\GK ({}_{\Tl }M')=\GK ({}_{\Tk }M)$. \item $e ({}_{\Tl
}M')=p^{|{\bf l} - \bk |}e ({}_{\Tk }M)$.
\end{enumerate}
\end{lemma}

{\it Proof}. Let $M_0$ be a finite dimensional generating subspace
for the $\Tk$-module $M=\Tk M_0$, $M_i:=\CT_{\bk , i}M_0$, $i\geq
0$. Then $\dim_K(M_i)=\frac{e(M)}{d!}i^d+\cdots $, $i\gg 0$ where
$d=\GK (M)$. $M'=\bigoplus_{0\leq \beta <\p^{{\bf l}-\bk
}}\der^{[\p^\bk \beta ]}\t M$ where $ \der^{[\p^\bk \beta ]}:=
\der_1^{[p^{k_1} \beta_1 ]}\cdots \der_n^{[p^{k_n} \beta_n ]}$,
and
$$\bigoplus_{0\leq \beta <\p^{{\bf l}-\bk }}\der^{[\p^\bk \beta
]}\t M_{i-p^{|{\bf l}|}}\subseteq M_i':=\CT_{{\bf l},
i}M_0\subseteq \bigoplus_{0\leq \beta <\p^{{\bf l}-\bk
}}\der^{[\p^\bk \beta ]}\t M_i, \;\;\; i\gg 0.$$ Therefore,
\begin{eqnarray*}
\frac{p^{|{\bf l}-\bk |}e(M)}{d!}i^d+\cdots &=& p^{|{\bf l}-\bk
|}\dim_K(M_{i-p^{|{\bf l}|}})\leq \dim_K(M_i')
=\frac{e(M')}{d!}i^d+\cdots \\
& \leq & p^{|{\bf l}-\bk |}\dim_K(M_i)=\frac{p^{|{\bf l}-\bk
|}e(M)}{d!}i^d+\cdots ,
\end{eqnarray*}
and so $\GK ({}_{\Tl }M')=\GK ({}_{\Tk }M)$ and $e ({}_{\Tl
}M')=p^{|{\bf l} - \bk |}e ({}_{\Tk }M)$. $\Box $

\begin{theorem}\label{c23Jun05}
Let $M'=\Tk M_0'$ be a nonzero finitely generated $\Tk$-module,
$\dim_K(M_0')<\infty $, $\bk =(k_1, \ldots , k_n)\in
\mathbb{N}^n$, $k=\max (k_1, \ldots , k_n)$, and $M:= \CD
(P_n)\t_{\Tk }M'$.
 Let $\{
M_i':= \CT_{\bk , i}M_0'\}$ be a standard filtration for the
$\Tk$-module $M'$  and $\dim_K(M_i' )=\frac{e(M')}{d!}i^d+\cdots $
for $i\gg 0$ where $d=\GK (M')$. Let $\{ M_i:=F_iM_0'\}$ be the
filtration of standard type on the $\CD (P_n)$-module $M$. Then
\begin{enumerate}
\item $\dim_K(M_i )=\frac{e(M')}{p^{| \bk |}(n+d)!}i^{n+d}+\cdots
$ is an almost polynomial of period $p^k$ with coefficients from
$\frac{1}{p^{k(n+d)}(n+d)!}\mathbb{Z}$,  and
$e(M)=\frac{e(M')}{p^{|\bk |}}\in \frac{1}{p^{|\bk |}}\mathbb{N}$.
\item
 The dimension $\Dim (M)=n+d\geq n$ is equal to $t-1$ where $t$ is  the order of the
 pole of the Poincare series $P_M(\o )=\sum_{i\geq
 0}\dim_K(M_i)\o^i$ at the point $\o =1$, and the multiplicity
 $e(M)=(1-\o )^{\Dim (M)+1}P_M(\o )|_{\o =1}$. The dimension $\Dim (M)$
  of $M$ can be any natural number from the interval $[n,2n]$.
\end{enumerate}
\end{theorem}

{\it Proof}. The subalgebra $\L^{[p^\bk ]}:= \bigotimes_{i=1}^n
\L_i^{[p^{k_i}]}$ of $\CD (P_n)$ has the induced filtration
$$\{\L^{[p^\bk ]}_i:=\L^{[p^\bk ]}\cap
F_i=\bigoplus_{p^{k_1}\beta_1 +\cdots +p^{k_n}\beta_n\leq
i}K\der_1^{[p^{k_1}\beta_1]}\cdots \der_n^{[p^{k_n}\beta_n]}\}.$$
Therefore,
$$ P:= \sum_{i\geq 0}\dim_K(\L^{[p^\bk ]}_i)\o^i=\frac{1}{(1-\o
)\prod_{i=1}^n(1-\o^{p^{k_i}})}\;\;\; {\rm and}\;\;\; e_P:=
(1-\o)^{n+1}P|_{\o =1}=\frac{1}{p^{|\bk |}}.$$ It follows from the
equality $M=\L^{[p^\bk ]}\t M'$ that $M_i=\sum_{s+t\leq
i}\L^{[p^\bk ]}_s\t M_t'$. Therefore, $R:=\sum_{i\geq
0}\dim_K(M_i)\o^i=(1-\o )PQ$ where $Q:=\sum_{i\geq
0}\dim_K(M_i')\o^i$. By Corollary \ref{1izvres},
$e(M)=e_R=e_Pe_Q=\frac{1}{p^{|\bk |}}e(M')$ and $\Dim (M)=n+d\geq
n$, and so
 $\dim_K(M_i
)=\frac{e(M')}{p^{|\bk |}(n+d)!}i^{n+d}+\cdots $, by Lemma
\ref{izvres}. The rest is obvious (Lemma \ref{izvres}). $\Box$

\begin{theorem}\label{21Jul05}
(A classification of simple $T_k$-modules where $T_k=T_{k,1}$).
Let $K$ be a field of characteristic $p>0$ and $k\geq 0$.
\begin{enumerate}
\item The map $\Max (K[x])\ra \widehat{T}_k$, $\Gm \mapsto
[T_k(\Gm )]$ is a bijection with inverse $[M]\mapsto \ass_{K[x]}
(M)$ where $$T_k(\Gm ):=\begin{cases} T_k\t_{\TkGm} \TkGm  / \TkGm
(\Gm , \L_{[p^{k(\Gm )}], +})\simeq T_k\t_{\TkGm} K[x]/\Gm & ,
k\geq k(\Gm
), \\
T_k/T_k (\Gm ,\L_{[p^{k}], +})\simeq K[x]/\Gm &, k <k(\Gm ).
\end{cases}$$
 \item $\dim_K \, T_k(\Gm ):=\begin{cases} p^{k-k(\Gm )} [K[x]/\Gm
 :K]=p^k[(K[x]/\Gm)^{sep}:K]& ,
k\geq k(\Gm
), \\
[K[x]/\Gm
 :K]=p^{k (\Gm )}[(K[x]/\Gm)^{sep}:K]&, k <k(\Gm
),
\end{cases}$

and so

$$\frac{\dim_K \, T_k(\Gm
)}{p^k[(K[x]/\Gm )^{sep}:K]}:=\begin{cases} 1 & , k\geq k(\Gm ),
\\p^{k (\Gm )-k}&, k <k(\Gm ).
\end{cases}$$
 \item ${\rm End}_{T_k}(T_k(\Gm ))\simeq \begin{cases}
 {\rm End}_{T_{k(\Gm )}}(T_{k(\Gm )}(\Gm )) \simeq (K[x]/\Gm
 )^{sep}\simeq K[x^{p^{k(\Gm )}}]/(g)& ,
k\geq k(\Gm
), \\
K[x^{p^k}]/(g)&, k <k(\Gm ),
\end{cases}$

where $\Gm =(g)$ and $g=f(x^{p^{k(\Gm )}})$. ${\rm
End}_{T_k}(T_k(\Gm ))$ is a subfield of $K[x]/\Gm$ that contains
$(K[x]/\Gm )^{sep}$. ${\rm End}_{T_k}(T_k(\Gm ))=(K[x]/\Gm
)^{sep}$ iff $k\geq k(\Gm )$. \item $\CD (K[x])\t_{T_k}T_k (\Gm
)\simeq
\begin{cases}
 U(\Gm )& ,
k\geq k(\Gm
), \\
U(\Gm )^{p^{k(\Gm )-k}}&, k <k(\Gm ),
\end{cases}$

where $U(\Gm )$ is the simple $\CD (K[x])$-module from Lemma
\ref{1Lisim}. \item If $k\leq k(\Gm )$ then the factor algebra
$T_k/T_k\Gm \simeq M_{p^k}(K[x^{p^k}]/(g))$ where $\Gm =(g)$.

\end{enumerate}
\end{theorem}

{\it Proof}. Let $\CD =\CD (K[x])$.

$4$. Let $T_k (\Gm )$ be as in the second part of statement 1.  If
$k\geq k(\Gm )$ then $\CD\t_{T_k}T_k(\Gm )\simeq \CD
\t_{T_k}T_k\t_{T_{k(\Gm )}} T_{k(\Gm )}/T_{k(\Gm )}(\Gm ,
\L_{[p^{k(\Gm )}], +})\simeq U(\Gm )$.

If $k<k(\Gm)$ then the $\CD $-module $M:=\CD \t_{T_k}T_k(\Gm )$ is
an epimorphic image of the $\CD $-module
$\CD\t_{K[x]}K[x]/\Gm\simeq U(\Gm )^s$ for some $s\geq 1$ (Lemma
\ref{1l8Jul05}). Therefore, $M\simeq U(\Gm )^t$ for some $t$. By
Theorem \ref{23Jun05},
$$e(M)=p^{-k}\dim_K(T_k(\Gm ))=p^{k(\Gm
)-k}[(K[x]/\Gm)^{sep}:K],$$
and by Theorem \ref{8Jul05}, $e(U(\Gm )^t)=t[(K[x]/\Gm)^{sep}:K]$.
Therefore, $t=p^{k(\Gm )-k}$.

$1$. If $k\geq k(\Gm )$ then the $T_k$-module $T_k(\Gm )$ is
simple since $\CD_{T_k}$ is faithfully flat and the induced $\CD
$-module $\CD\t_{T_k}T_k(\Gm )$ is simple (by statement $4$).

If $k<k(\Gm )$ then the $K[x]$-module $T_k(\Gm )$ is simple, and
so the $T_k$-module $T_k(\Gm )$ is simple. Now, statement $1$
follows from statement $4$ and Theorem \ref{8Jul05}.

$2$ and $3$. These statements are obvious.

$5$. It follows from the decomposition 
\begin{equation}\label{TkTkg}
T_k/T_kg\simeq \L_{[p^k]}\t K[x]/(g)=\bigoplus_{0\leq
i<p^k}\der^{[i]}K[x]/(g)=\bigoplus_{0\leq i,
j<p^k}\der^{[i]}x^jK[x^{p^k}]/(g)
\end{equation}
that the algebra $T_k/T_kg$ is a simple algebra with the centre
$K[x^{p^k}]/(g)$ (use $\ad \, x$ and the fact that the $K[x]/(g)$
is  a field), and $\dim_K(T_k/T_kg)=p^k[K[x]/\Gm :K]$. By
(\ref{TkTkg}),  the $T_k/T_kg$-module
$$ U':=T_k/T_k(g, \L_{[p^k], +})\simeq K[x]/\Gm $$ is simple,
$\dim_K(U')=[K[x]/\Gm :K]=p^k[K[x^{p^k}]/(g) :K]$, and ${\rm
End}_{T_k/T_kg}(U')\simeq K[x^{p^k}]/(g)$. This implies that
$T_k/T_kg\simeq M_{p^k}(K[x^{p^k}]/(g))$ (this also proves
statement 1, the case  $k<k(\Gm )$).
 $\Box $

Let $\bk =(k_1, \ldots , k_n)\in \mathbb{N}^n$. We are going to
classify simple $\Tk$-modules (Theorem \ref{20Jul05}). The algebra
$\Tk$ is a somewhat commutative algebra which is a finitely
generated module over its centre. By Quillen's Lemma, {\em every
simple $\Tk$-module has finite dimension over $K$}.  Given a
finite dimensional $T_k$-module $M$. Then $M=\bigoplus_{\Gm \in
\Max (P_n)}M^\Gm $ is a direct sum of its submodules $M^\Gm
:=\cup_{i\geq 1}\ann_M(\Gm^i)$. If, in addition, the $\Tk$-module
$M$ is simple then $M=M^\Gm$ for a uniquely determined maximal
ideal $\Gm $ of $P_n$ and $M$ is an epimorphic image of the finite
dimensional $\Tk$-module $\Tk / \Tk \Gm \simeq \Tk
\t_{P_n}P_n/\Gm\simeq \L_{[p^\bk ]}\t P_n/\Gm$, $\dim_K(\Tk
\t_{P_n}P_n) = p^{|\bk |}[P_n/\Gm :K]$.

Suppose that $\bk \leq k(\Gm ):=(k_1', \ldots , k_n')$ (i.e.
$k_1\leq k_1',\ldots , k_n\leq k_n'$). Let $g(\Gm )=(g_1, \ldots ,
g_n)$ where $g_i(x_i)=f_i(x_i^{p^{k_i'}})$. We keep the notation
as in (\ref{GmPnms}). Consider natural maps
$$ \prod_{\nu =1}^n \G_\nu \simeq \bigotimes_{i=1}^n
K[x_i^{p^{k_i'}}]/(g_i)\ra \bigotimes_{i=1}^n
K[x_i^{p^{k_i}}]/(g_i)\stackrel{\phi}{\longrightarrow}\bigotimes_{i=1}^n
K[x_i]/(g_i)\stackrel{\pi}{\longrightarrow}P_n/\Gm.$$ By
(\ref{GmPnms}), we have the inclusions of fields: 
\begin{equation}\label{Gkbm}
(P_n/\Gm )^{sep}=\G (\Gm )\subseteq \G (\bk , \Gm ):= \im (\pi
\circ \phi )\subseteq P_n/\Gm .
\end{equation}
Consider the factor algebra (Theorem \ref{21Jul05})
$$ \Tk / \Tk g(\Gm )\simeq \bigotimes_{i=1}^n
T_{k_i}/T_{k_i}g_i\simeq \bigotimes_{i=1}^n
M_{p^{k_i}}(K[x^{p^{k_i}}]/(g_i))\simeq M_{p^{|\bk |}}
(\bigotimes_{i=1}^n K[x^{p^{k_i}}]/(g_i)).$$ The $\Tk$-module $\Tk
/ \Tk \Gm$ is, in fact, a $\Tk / \Tk g(\Gm )$-module, or even,
$M_{p^{|\bk |}}(\G (\bk , \Gm ))$-module (since $e_\mu \ra 0$ if
$\mu \neq \nu$, see (\ref{GmPnms})). Let 
\begin{equation}\label{Vbkm}
V(\bk , \Gm ):=\G (\bk , \Gm )^{p^{|\bk |}}
\end{equation}
be the only simple module of the matrix algebra $M(\bk , \Gm
):=M_{p^{|\bk |}}(\G (\bk , \Gm ))$. Then, $\dim_K \, V(\bk , \Gm
)=p^{|\bk |}[\G (\bk , \Gm ):K]$, and 
\begin{equation}\label{EMkm}
{\rm End}_{M(\bk , \Gm )}(V(\bk , \Gm ))\simeq \G (\bk , \Gm ).
\end{equation}
It follows that $\Tk / \Tk \Gm \simeq V(\bk , \Gm )^{p^\nu}$ where
$$ p^\nu =\frac{\dim_K(\Tk / \Tk \Gm )}{\dim_K(V(\bk , \Gm ))}=
\frac{p^{|\bk |}[P_n/\Gm :K]}{p^{|\bk |}[\G (\bk , \Gm
):K]}=\frac{[P_n/\Gm :K]}{[\G (\bk , \Gm ):K]}$$ by (\ref{Gkbm}).
Therefore, $V(\bk , \Gm )$ {\em is the only simple $\Tk$-module
which is annihilated by a power of the maximal ideal $\Gm $
(provided $\bk\leq k(\Gm  )$)}.

For $\alpha =(\alpha_1, \ldots , \alpha_n), \beta =(\beta_1,
\ldots , \beta_n)\in \mathbb{N}^n$, let $\min (\alpha , \beta )=
(\min (\alpha_1 , \beta_1 ), \ldots , \min (\alpha_n , \beta_n ))$
and $\max (\alpha , \beta )= (\max (\alpha_1 , \beta_1 ), \ldots ,
\max (\alpha_n , \beta_n ))$.

\begin{theorem}\label{20Jul05}
(A classification of simple $\Tk$-modules). Let $K$ be a field of
characteristic $p>0$ and $\bk  =(k_1, \ldots , k_n)\in
\mathbb{N}^n$.
\begin{enumerate}
\item The map $\Max (P_n)\ra \widehat{T}_\bk$, $\Gm \mapsto [\Tk
(\Gm )]$ is a bijection with inverse $[M]\mapsto \ass_{P_n} (M)$
where
$$\Tk (\Gm ):=\begin{cases} \Tk\t_{\TkGm} V(\Gm ) & , k\geq
k(\Gm
), \\
V(\bk , \Gm )&, k \leq k(\Gm ),\\ \Tk \t_{T_{\min (\bk , k(\Gm
))}} V( \min (\bk , k(\Gm )), \Gm )& , \text{otherwise}.
\end{cases}$$
 \item $\dim_K \, T_{\bk}(\Gm ):=\begin{cases} p^{|\bk |} [(P_n/\Gm)^{sep}:K]& ,
k\geq k(\Gm
), \\
p^{|\bk |}[\G (\bk , \Gm ):K]&, k \leq k(\Gm ),\\
p^{|\bk |}[\G (\min (\bk , k(\Gm )),  \Gm ):K]&, \text{otherwise},
\end{cases}$

$$r(\bk , \Gm ):=\frac{\dim_K \, \Tk (\Gm
)}{p^{|\bk |}[(P_n/\Gm )^{sep}:K]}:=\begin{cases} 1 & , k\geq
k(\Gm ),\\
[\G (\bk , \Gm ): (P_n/\Gm )^{sep}]&, k \leq k(\Gm ),\\
[\G (\min (\bk , k(\Gm )), \Gm ):(P_n/\Gm)^{sep}]&,
\text{otherwise},
\end{cases}$$
and $r(\bk , \Gm )=p^s$ for some $s=s(\bk , \Gm )\in \mathbb{N}$.
 \item ${\rm End}_{\Tk }(\Tk (\Gm ))\simeq \begin{cases}
 (P_n/\Gm )^{sep}& ,
k\geq k(\Gm
), \\
\G (\bk , \Gm )&, k \leq k(\Gm ),\\
\G (\min (\bk , k(\Gm )), \Gm )& ,\text{otherwise}.
\end{cases}$

 ${\rm End}_{\Tk }(\Tk (\Gm ))$ is a subfield of $P_n/\Gm$ that contains
$(P_n/\Gm )^{sep}$. \item $\dim_K \, \Tk (\Gm )=p^{|\bk |}
\dim_K\, {\rm End}_{\Tk }(\Tk (\Gm ))$.

\item $\CD (P_n)\t_{\Tk }\Tk (\Gm ))\simeq U(\Gm )^{r(\bk , \Gm
)}$.
\end{enumerate}
\end{theorem}

{\it Proof}. $1$. Let  $\Gm \in \Max (P_n)$. If $\bk \geq k(\Gm )$
then the $\CD (P_n)$-module $\CD (P_n)\t_{\TkGm} V(\Gm ) \simeq
\CD (P_n)\t_{\Tk } (\Tk \t_{\TkGm} V(\Gm ))$ is simple (Theorem
\ref{8Jul05}). Therefore, $\Tk \t_{\TkGm} V(\Gm )$ must be a
simple $\Tk$-module.

If $k\leq k(\Gm )$ then $V(\bk , \Gm )$ is a simple $\Tk$-module.

In the remaining case, one can prove that any nonzero $\Tk
$-submodule of $M:= \Tk \t_{\Tl }V({\bf l} , \Gm )$, ${\bf l}:=
\min (\bk , k(\Gm ))$, has a nonzero intersection with the simple
$\Tl$-submodule $V({\bf l}, \Gm )$ of $M$. Therefore, $M$ is a
simple $\Tk$-module. The rest of statement $1$ is obvious (see
Theorem \ref{8Jul05} and the arguments preceding Theorem
\ref{20Jul05}).

$2$. If $\bk \geq k(\Gm )$ then $$\dim_K (\Tk (\Gm ))=p^{|\bk
-k(\Gm )|}\dim_K\, V(\Gm )=p^{|\bk -k(\Gm )|}p^{|k(\Gm )|}[\G (\Gm
):K]=p^{|\bk |}[(P_n/\Gm )^{sep}:K].$$

If $\bk \leq k(\Gm )$ then the result follows from (\ref{Vbkm}).
In the third case, let $l=\min (\bk , k(\Gm ))$. Then $$\dim_K
(\Tk (\Gm ))=p^{|\bk -l|}\dim_K\, V(l, \Gm )=p^{|\bk
-l|}p^{|l|}[\G (l, \Gm ):K]=p^{|\bk |}[\G (l, \Gm ):K].$$ The rest
of statement $2$ is obvious.

$3$. Evident.

$4$. This follows from statement $2$.

$5$. By Lemma \ref{1l8Jul05}, the $\CD (P_n)$-module $N:=\CD
(P_n)\t_{\Tk }\Tk (\Gm )$ is isomorphic to $U(\Gm )^r$ for some
$r\in \mathbb{N}$. By Theorem \ref{c23Jun05}, the multiplicity of
the $\CD (P_n)$-module $N$ is equal to $ e(N)=p^{-|\bk |}\dim_K\,
\Tk (\Gm )$. By Theorem \ref{8Jul05}, $e(U(\Gm )^r)=r [(P_n/\Gm
)^{sep}:K]$, hence $$r=\frac{\dim_K \, \Tk (\Gm )}{p^{|\bk
|}[(P_n/\Gm )^{sep}:K] }=r(\bk , \Gm ).\;\;\; \Box $$

\begin{corollary}\label{ac20Jul05}
$p^{|\bk |}|\dim_K(M)$ for all finite dimensional $\Tk$-modules
$M$.
\end{corollary}

\begin{theorem}\label{t12Jul05}
Let $M$ be a nonzero holonomic finitely presented $\CD
(P_n)$-module. Then its multiplicity is a natural number.
\end{theorem}

{\it Proof}. This follows directly from Corollary \ref{ac20Jul05},
(\ref{Mfpr}),  and Theorem \ref{c23Jun05}. $\Box$


\section{Holonomic sets of subalgebras with multiplicity, every holonomic
 $\CD (P_n)$-module has finite length}\label{HSSA}
In this section, $K$ is an {\em arbitrary} field of characteristic
$p>0$ if it is not stated otherwise.

In this section, the  concept of {\em holonomic set of subalgebras
with multiplicity} is introduced which is a crucial one in the
proof of the analogue of the inequality of Bernstein for the
algebra $\CD (P_n)$ (Theorem \ref{c27Jn05}) and in the proof of
the fact that each holonomic $\CD (P_n)$-module has finite length
and the length does not exceed the multiplicity (Theorem
\ref{hDPnfl}). It is proved that $n\leq \Dim (L)\leq 2n$ for each
nonzero finitely generated $\CD (P_n)$-module $L$, and,  for {\em
each real} number $d\in [n, 2n]$, there exists a cyclic $\CD
(P_n)$-module $M$ with $\Dim (M)=d$ (Theorem \ref{3Jl05}), and
there exists a cyclic non-holonomic $\CD (P_n)$-module $N$ with
$\Dim (N)=n$ (Proposition \ref{Mbknh}).

{\bf Holonomic sets of subalgebras}. Let $A$ be an algebra over an
arbitrary field $K$ with a {\em finite dimensional} filtration
$\CA =\{ A_i\}_{i\geq 0}$  such that $\Dim (A):=\g (\dim_K \,
A_i)<\infty$. Any subalgebra $B$ of the algebra $A$ has the
induced finite dimensional filtration $\CB =\{ B_i:=B\cap A_i\}$
and $\Dim (B):=\g (\dim_K \, B_i)\leq \Dim (A)<\infty$.

{\it Definition}. A set $\CC =\{ C_\nu \}_{\nu \in \CN }$ of
subalgebras of the algebra $A$ is called a {\bf sub-holonomic} set
if there exists a real positive number $h_\CC$ such that for each
nonzero $A$-module $M$ there exists $\nu \in \CN$ and a finitely
generated $C_\nu$-submodule $M_\nu$ of $M$ such that $\Dim
({}_{C_\nu}M_\nu)\geq h_\CC$ or, equivalently, there exists a
nonzero finite dimensional vector subspace $V$ of $M$ such that
$\g (\dim_K(C_{\nu , i}V))\geq h_\CC$ for some $\nu$ where $\{
C_{\nu , i}:=C_\nu \cap A_i\}$ is the induced filtration on the
algebra $C_\nu$.

Surprisingly, the following simple observation yields an idea of
another proof of the inequality of Bernstein for the ring of
differential operators in positive characteristic, and, more
importantly, it produces an analogue of multiplicity.

\begin{lemma}\label{cf2}
Let $A$, $\CC =\{ C_\nu \}_{\nu \in \CN}$, and $h_\CC$ be as
above. Then $\Dim (M)\geq h_\CC$ for all nonzero finitely
generated $A$-modules $M$.
\end{lemma}

{\it Proof}. For a nonzero finitely generated $A$-module $M$, we
have $\g (\dim_K \, C_{\nu , i}V)\geq h_\CC$ for some nonzero
finite dimensional $K$-subspace $V$ of $M$. Let $M_0$ be a finite
dimensional generating subspace for the $A$-module $M$ that
contains $V$. Then
$$ \Dim (M)=\g (\dim_K\, A_iM_0)\geq \g (\dim_K\, C_{\nu ,
i}V)\geq h_\CC. \;\;\; \Box $$

{\it Definition}. A set $\CC =\{ C_\nu \}_{\nu \in\CN}$ of
subalgebras of the algebra $A$ is called a {\bf sub-holonomic} set
of {\bf degree} $n$ and with {\bf leading coefficient} $l$ where
$n$ and $l$ are positive real numbers if for each nonzero
$A$-module $M$ there exists a nonzero finite dimensional
$K$-vector subspace $V\subseteq M$ and an algebra $C_\nu$ such
that $\dim_K(C_{\nu , i}V)\geq li^n+\cdots $ (where the three dots
mean a function which is negligible comparing to $i^n$, i.e.
$o(i^n)$). If $n$ is a {\em natural} number then $e:=n!l$ is
called the {\bf multiplicity} for $\CC$. If, in addition, $n=h_A$
then the set $\CC$ is called a {\bf holonomic} set of subalgebras
with {\bf leading coefficient} $l$  (or {\bf multiplicity} $e$)
for the algebra $A$ where $h_A:= \inf \{ \Dim (M)\, | \, M$ is a
nonzero finitely generated $A$-module$\}$ is the {\em holonomic
number} for the algebra $A$ with respect to the filtration $\CA $.

\begin{theorem}\label{ChBI1}
If there exists a holonomic set $\CC =\{ C_\nu \}_{\nu \in \CN}$
of subalgebras with the leading coefficient $l_\CC$ for the
algebra $A$ then every holonomic $A$-module has finite length.
Moreover, if $\{ M_i\}$ is a filtration of standard type on a
holonomic $A$-module $M$ then the length of the $A$-module $M$ is
$\leq \frac{l(M)}{l_\CC}$ where $l_\CC$ is the leading
coefficients for $\CC$, $\dim_K(M_i)\leq l(M)i^n+\cdots $, $i\gg
0$, and the three dots mean $o(i^n)$.
\end{theorem}

{\it Proof}.  It suffices to prove the last statement. Suppose to
the contrary that there exists a holonomic $A$-module $M$ of
length $> \frac{l(M)}{l_\CC}$, we seek a contradiction. Then one
can choose a strictly ascending chain of submodules in $M$:
$0=M_0'\subset M_1'\subset \cdots \subset M_t'\subseteq M$ with $t
> \frac{l(M)}{l_\CC}$. For each factor module $M_j'/M_{j-1}'$, fix a
nonzero finite dimensional subspace $\overline{V}_j\subseteq
M_j'/M_{j-1}'$ such that $\dim_K(C_{\nu (j), i}\overline{V}_j)\geq
l_\CC i^n+\cdots $, $i\gg 0$,  for some $\nu (j)$. Let $V_j$ be a
finite dimensional subspace of $M_j'$ such that
$\overline{V}_j=V_j+M_{j-1}'$. Fix $s\geq 1$ such that $V_1+\cdots
+V_t\subseteq M_s$. Then for $i\gg 0$,
\begin{eqnarray*}
 tl_\CC i^n+\cdots &\leq &
\sum_{j=1}^t \dim (C_{\nu (j), i}\overline{V}_j)\leq \dim
(\sum_{j=1}^tC_{\nu (j), i}V_j)\leq \dim \, M_{i+s}\\
&\leq &l(M)(i+s)^n+\cdots = l(M)i^n+\cdots,
\end{eqnarray*}
 and so $tl_\CC \leq l(M)$, a contradiction.
$\Box $

{\it Definition}. We say that a subalgebra of $\CD (P_n)$ is of
type $P_s\t \L (n-s)$ (resp. of type $P_s\t
\bigotimes_{i=1}^{n-s}\L_{-1}^{[p^{k_i}]})$ if after changing, if
necessary, the order of the tensor multiples in $\CD (P_n)=\CD
(P_1)\t \cdots \t \CD (P_1)$ the algebra is equal to $P_1^{\t s}\t
\L_{-1}^{\t (n-s)}$ (resp. $P_s\t
\bigotimes_{i=1}^{n-s}\L_{-1}^{[p^{k_i}]})$).

For $\e \in \{ \pm 1\}^n$, $|\e |$ denotes the number of {\em
negative} coordinates (eg, $|(-1, \ldots , -1)|=n$ and $|(1,
\ldots , 1)|=0$).

\begin{theorem}\label{27Jn05}
Let $K$ be an arbitrary field of characteristic $p>0$. For any
nonzero $\CD (P_n)$-module $M$ there exists a subalgebra $\L$ of
the type $P_s\t \bigotimes_{i=1}^{n-s}\L_{-1}^{[p^{k_i}]}$
 of $\CD (P_n)$ for some $k_i\geq 0$ and a finite dimensional
$K$-subspace $V$ of $M$ such that $\dim_K (V)\geq p^{k_1+\cdots
+k_{n-s}}$ and the natural map $\L \t V\ra \L V$, $\l \t v\mapsto
lv$ (in $M$), is an isomorphism of $\L$-modules.
\end{theorem}

{\it Proof}. The polynomial algebra $P_n$ is a commutative
Noetherian domain, so any {\em maximal} (with respect to
inclusion) element of the set of annihilators $\{ \ann_{P_n}(v)\,
| \, 0\neq v\in M\}$ is a {\em prime} ideal. Fix  such a prime
ideal, say $\Gp =\ann_{P_n}(v)$ for some $0\neq v\in M$. Without
loss of generality one can assume that $M=\CD (P_n)v$. Then the
$\CD (P_n)$-module $M$ is an epimorphic image of the $\CD
(P_n)$-module $\CD (P_n)/\CD (P_n)\Gp \simeq \CD (P_n)\t_{P_n}
P_n/\Gp = \cup_{i\geq 1}\ann (\Gp^i)$. So, any element of $M$ is
annihilated by a power of the ideal of $\Gp$. To prove the theorem
we use induction on $n$.

The case $n=1$. There are two cases: either $\Gp =0$ or otherwise
$\Gp$ is a maximal ideal of the polynomial algebra $P_1:=K[x]$. If
$\Gp =0$ then $K[x]v \simeq K[x]$ is an isomorphism of
$K[x]$-modules, and so it suffices to take $s=1$ and $V=Kv$. If
$\Gp \neq 0$ then the ideal $\Gp$ is generated by an irreducible
polynomial of $K[x]$. Then the result follows from Lemma
\ref{1Lisim}.

Suppose that $n>1$ and the theorem is true for all $n'<n$. Now, we
use a second downward induction on the Krull dimension $d=\Kdim
(P_n/\Gp )$ of the algebra $P_n/\Gp$ starting with $d=n$, i.e.
$\Gp =0$. In this case, it suffices to take $\e =(1, \ldots , 1)$
and $V=Kv$, since $P_nv\simeq P_n$ is an isomorphism of
$P_n$-modules.

Suppose now that $d<n$ and the result is true for all $d'$ such
that $d<d'\leq n$. The field of fractions $Q={\rm Frac} (P_n/\Gp
)$ of the domain $P_n/\Gp $ has transcendence degree $d$ over the
field $K$, and it is generated by the elements $\bx_i =x_i+\Gp$,
$i=1, \ldots , n$. Up to order of the elements $\bx_i$, we can
assume that the elements $\bx_1 , \ldots , \bx_d$ are
algebraically independent over $K$, and so $Q$ is the  finite
 field extension of its subfield $Q_d:=K(\bx_1 , \ldots , \bx_d)$
of rational functions. Then $P_n=P_d\t P_{n-d}$ where
$P_d=K[x_1,\ldots , x_d]$ and $P_{n-d}=K[x_{d+1}, \ldots , x_n]$.
Correspondingly, $\CD (P_n)=\CD (P_d)\t \CD  (P_{n-d})$ and the
localization $Q_d\t_{P_d}\CD (P_n)=Q_d\t_{P_d}\CD (P_d)\t \CD
(P_{n-d})\simeq \CD (Q_d)\t \CD (P_{n-d})$ of the algebra $\CD
(P_n)$ at $P_d\backslash \{ 0\}$ contains the subalgebra $Q_d\t
\CD (P_{n-d})\simeq \CD_{Q_d}(Q_d [x_{d+1}, \ldots , x_n])$ which
is the ring of differential operators over the field $Q_d$ of the
polynomial algebra $Q_d[x_{d+1}, \ldots , x_n]$ in $n-d$ variables
over the field $Q_d$. By the choice of the prime ideal $\Gp$ and
the elements $x_1, \ldots , x_d$, the $\CD (P_n)$-module $M$ is a
{\em submodule} of its localization $Q_d\t_{P_d}M$ (use the fact
that $\Gp \cap P_d=0$ and $M=\cup_{i\geq 1}\ann (\Gp^i))$ which is
a $Q_d\t_{P_d}\CD (P_n)$-module, and, by restriction, it is a
$\CD_{Q_d}(Q_d[x_{d+1}, \ldots , x_n])$-module. Since $n-d<n$, by
induction, one can find a subalgebra $\L'=P_s
(Q_d)\t_{Q_d}\bigotimes_{i=1}^{n-d-s}\L_{-1}^{[p^{k_i}]}$ for some
$k_i\geq 0$ and
 a finite dimensional $Q_d$-submodule of $\CD_{Q_d}(Q_d[x_{d+1},
  \ldots , x_n])$, say $\CV =Q_d\t V$, of $Q_d\t_{P_d}M$ (where
 $V$ is a finite dimensional $K$-submodule of $M$) such that
 $\dim_{Q_d}(\CV )=\dim_K(V)\geq p^{k_1+\cdots +k_{n-d-s}}$ and the natural  map
 $\L'\t_{Q_d}\CV \ra \L'\CV$ is an isomorphism of $\L'$-modules.
 Let $\L =P_d\t P_s\t \bigotimes_{i=1}^{n-d-s}\L_{-1}^{[p^{k_i}]}$ (a subalgebra of $\CD (P_n)$).
 Then the natural map $\L \t V\ra \L V$ is an isomorphism. By
 induction, the proof now is complete.
$\Box $

There is another proof of the inequality of Bernstein in prime
characteristic.

\begin{theorem}\label{c27Jn05}
Let $K$ be an arbitrary field of characteristic $p>0$. Then $\Dim
(M)\geq n$ for each nonzero finitely generated $\CD (P_n)$-module
$M$.
\end{theorem}

{\it Proof}. By Theorem \ref{27Jn05}, $\Dim (M)\geq \Dim
(P_s\t\bigotimes_{i=1}^{n-s}\L_{-1}^{[p^{k_i}]})=n$ for some $s$
and $k_i\geq 0$. $\Box $

The next theorem gives {\em explicitly} examples of sets of
holonomic subalgebras with multiplicity $1$ for the algebra $\CD
(P_n)$.
\begin{theorem}\label{1c27Jn05}
Let $K$ be an arbitrary field of characteristic $p>0$, $\CC =\{$a
subalgebra of $\CD (P_n)$ of type
$P_s\t\bigotimes_{i=1}^{n-s}\L_{-1}^{[p^{k_i}]}$, $0\leq s\leq n$,
$k_i\geq 0\}$ and $\CC'=\{ \L_\e \, | \, \e\in \{\pm 1\}^n \}$.
Then the sets $\CC $ and $\CC'$ are holonomic sets of subalgebras
with multiplicity $1$ for the ring of differential operators $\CD
(P_n)$ (equipped with the canonical filtration).
\end{theorem}

{\it Proof}. $\Dim (A)=n$ for all algebras $A$ from $\CC \cup
\CC'$. By Theorem \ref{27Jn05}, $\CC$ is a holonomic set of
subalgebras with multiplicity $1$. Each algebra from the set $\CC$
is a subalgebra of one of the algebras from the set $\CC'$, and
$\CC'\subseteq \CC$.  Therefore, $\CC'$ is a holonomic set of
subalgebras with multiplicity $1$ for the algebra $\CD (P_n)$.
$\Box$

\begin{theorem}\label{hDPnfl}
Let $K$ be an arbitrary field of characteristic $p>0$. Then each
holonomic  $\CD (P_n)$-module has finite length and its length
does not exceed the multiplicity (i.e. the length of $M$ $ \leq
\frac{l(M)}{n!}$, see Theorem \ref{ChBI1}).
\end{theorem}

{\it Proof}. This follows from  Theorems \ref{ChBI1} and
\ref{1c27Jn05}. $\Box$

\begin{theorem}\label{2Sep05}
Each holonomic $\CD (P_n)$-module is cyclic.
\end{theorem}

{\it Proof}. Repeat the characteristic zero proof which uses only
that facts that each holonomic module has finite length and the
ring of differential operators is simple and it is not an artinian
module over itself. $\Box $

{\bf An example of a cyclic non-holonomic $\CDPn$-module $M$ with
$\Dim (M)=n$}. Consider the subalgebra  $\L = \L_{-1}$ in $\CD
(P_1)$. Given an {\em infinite} sequence of natural numbers $\bk$:
$0<k_1<k_2<\cdots $. Consider the cyclic $\L$-module
$$M(\bk )=\L /\L (\der^{[j]}\, | \, j\in [1, p^{k_1}-1]\cup [p^{k_1}+1, p^{k_2}-1]
\cup [p^{k_2}+1, p^{k_3}-1]\cup \ldots )\simeq K
\overline{1}\oplus \oplus_{s\geq 1}K\der^{[p^{k_s}]}\overline{1}$$
where $\overline{1}$ is the canonical generator for the
$\L$-module $M(\bk )$. Consider the filtration of standard type on
$M(\bk )$ (with respect to the canonical filtration on $\CD
(P_1))$ $ \{ M(\bk )_i:=\L_i K\overline{1}=K\overline{1}\oplus
K\der^{[p^{k_1}]}\overline{1}\oplus \cdots \oplus
K\der^{[p^{k_s}]}\overline{1} \}$ where $s=s(i)$ satisfies
$p^{k_s}\leq i<p^{k_{s+1}}$, and so $\dim_K(M(\bk )_i)=s(i)+1$.
For each $t\geq 1$, $\oplus_{t\geq
s}K\der^{[p^{k_s}]}\overline{1}$ is a submodule of $M(\bk )$, the
corresponding factor module is denoted by $M(k_1, \ldots ,
k_{s-1})=K\overline{1}\oplus K\der^{[p^{k_1}]} \overline{1}\oplus
\cdots \oplus K\der^{[p^{k_{s-1}}]} \overline{1} $. In particular,
$M(\emptyset )=K$.

The next lemma shows that the growth of the module $M(\bk )$ can
be arbitrary slow.

\begin{lemma}\label{Laslgr}
For any non-decreasing function $f:\mathbb{N}\ra \mathbb{N}$ that
takes infinitely many values and $f(0)=1$, there exists a module
$M(\bk )$ such that $\dim_K(M(\bk )_i)\leq f(i)$ for all $i\geq 0$
(for an arbitrary non-decreasing  function $f$ with $f(0)=1$ there
exists a $\L$-module $M(k_1,\ldots , k_s)$ such that
$\dim_K(M(k_1,\ldots ,k_s)_i)\leq f(i)$ for all $i\geq 0$).
\end{lemma}

{\it Proof}. One can easily find an infinite sequence of natural
numbers $0<k_1<k_2<\cdots $ satisfying the property that $\#
 \{ j\,
| \, p^{k_j}<f(i)\} \leq f(i)$ for all $i\geq 0$. $\Box $

\begin{proposition}\label{Mbknh}
There exists a cyclic non-holonomic non-Noetherian  $\CD
(P_n)$-module $M$ such that $\Dim (M)=n$.
\end{proposition}

{\it Proof}. Fix a $\L$-module $M(\bk )$ from Lemma \ref{Laslgr}
which has zero growth, i.e. $\g (d_i)=0$ where $d_i=\dim_K(M(\bk
)_i)$. The $\L$-module $M(\bk )$ is not a Noetherian module, hence
the induced $\CD (P_1)$-module $\CD (P_1)\t_{\L}M (\bk )=P_1\t
M(\bk )$ is a cyclic non-Noetherian $\CD (P_1)$-module. Since $\CD
(P_n)=\CD (P_{n-1})\t \CD (P_1)$, the $\CD (P_n)$-module
$$ \CM (\bk ):= P_{n-1}\t (\CD (P_1)\t_{\L}M (\bk ))\simeq P_n\t
M(\bk )$$ is a cyclic non-Noetherian $\CD (P_n)$-module. Let $\{
\CM_i\}$ be the filtration of standard type associated with the
generating space $K\overline{1}$ for the $\CD (P_n)$-module $\CM
(\bk )$ and the canonical filtration on $\CD (P_n)$. Then
$$ \dim_K(\CM_i)={i+n\choose n}+{i+n-p^{k_1}\choose n}+\cdots +{i+n-p^{k_{d_i-1}}\choose
n}\leq d_i {i+n\choose n}.$$ It follows that $\Dim (\CM (\bk ))=\g
(\dim_K\, \CM_i)\leq \g (d_i{i+n\choose n})=\g (d_i)+n=n$.

Fix an arbitrary natural number $l$, then for all $i\gg 0$,
\begin{eqnarray*}
\dim_K\, \CM_i & > &{i+n\choose n}+{i+n-p^{k_1}\choose n}+\cdots
+{i+n-p^{k_l}\choose n} \\
&\geq & (l+1){i+n-p^{k_l}\choose n}=\frac{(l+1)}{n!}i^n+\cdots .
\end{eqnarray*}
Therefore, $\Dim (\CM (\bk ))=n$ and $\CM (\bk )$ is not a
holonomic $\CD (P_n)$-module.  $\Box $

{\bf An example of a cyclic $\CD (P_n)$-module $M$ with $\Dim
(M)=d$ for each $d\in [n,2n]$}.  Given an ascending sequence $b$
of positive real numbers $b_0=0<b_1<b_2<\cdots $ with $\lim_{i\ra
\infty} b_i=\infty$ and a sequence $s$ of positive real numbers
$s_1, s_2, \ldots $. Consider a {\em continuous piecewise linear}
function $f=f_{b,s}:\mathbb{R}_+\ra \mathbb{R}_+:=\{ r\in
\mathbb{R}\, | \, r\geq 0\}$ such that $f(0)=1$ and on  each
interval $[b_{i-1}, b_i]$ it is a linear function with slope
$s_i$. Then $b$ and $s$ are called the sequence of {\em breaking
points} and {\em slopes} for $f_{b,s}$ respectively.

Let us explain an idea of the proof of Lemma \ref{l3Jl05} which is
an essential step in proving Theorem \ref{3Jl05}. For any $r\in
\mathbb{R}$ such that $0<r<1$, each linear function $ax+b$ with
$a>0$ grows faster then the function $y=x^r+1$. The function
$y=x^r+1$ can be approximated by a function $f_{b,s}$ such that
both functions have the {\em same} growth $r$ and the graph of the
function $f_{b,s}$ lies {\em below} the graph of the function
$y=x^r+1$. When the slopes tend to zero  sufficiently fast then
the restriction of the function $f_{b,s}$ to the set of natural
numbers has the same growth. If we alter such a restriction at any
subset of natural numbers such that the values at infinitely many
breaking points remain unchanged, the new function from
$\mathbb{N}$ to $\mathbb{R}_+$  is increasing, and its graph lies
below the graph of $f_{b,s}$, then the altered function has growth
$r$. Such an altered function will be the function that defines
the  growth of the $\L$-module $M_r$ from Lemma \ref{l3Jl05}.

\begin{lemma}\label{l3Jl05}
Let $\L =K[\der^{[1]}, \der^{[2]}, \ldots , ]$ and $r\in
\mathbb{R}$ , $0<r<1$.
\begin{enumerate}
\item There exists a cyclic $\L$-module $M_r$ such that $\Dim
(M_r)=r$. \item The $\CD (P_1)$-module $\CM_r :=\CD
(P_1)\t_{\L}M_r$ has dimension $\Dim (\CM_r)=1+r$.
\end{enumerate}
\end{lemma}

{\it Proof}. $1$.  In this proof all functions are from
$\mathbb{N}$ to  $\mathbb{R}_+$. We are going to find an
approximation of the function $y=x^r+1$ by a function of the type
$f=f_{b,s}$ where $s: s_1, p^{-k_1}, s_2, p^{-k_2}, s_3, p^{-k_3},
\ldots $ where $0<k_1<k_2<\cdots $ and $b:
b_0=0<p^{k_1}<b_1<p^{k_2}<b_2<\cdots $. Fix  a sufficiently big
natural number, say $k_1$. Then $s_1$ is the slope of the linear
function passing through the points $(0,1)$ and $(p^{k_1}, 2)$,
and so $f(p^{k_1})=2$. Let $b_1$ be the largest {\em natural}
number of the form $i_1p^{k_1}$ such that $i_1\in \mathbb{N}$ and
$f(b_1)<y(b_1)$. Then fix a sufficiently large natural number, say
$k_2$, such that $b_1<p^{k_2}$. Then $s_2$ is the slope of the
linear function passing through the points $(b_1, f(b_1))$ and
$(p^{k_2}, f(p^{k_2}):=f(b_1)+1)$. Let $b_2$ be the largest
natural number of the type $i_2p^{k_2}$ such that $ i_2\in
\mathbb{N}$ and $p^{k_2}\leq b_2$ and $f(b_2)<y(b_2)$. We continue
in a similar fashion. The graph of the function $f$ lies below the
graph of the function $y$. `Sufficiently big'  in the choices
above means that $\lim_{i\ra
\infty}\frac{y(b_i)-f(b_i)}{f(b_i)}=0$ (this can be easily
achieved if the sequence $0<k_1<k_2<\cdots $ grows sufficiently
fast, this condition guarantees that the values of the function
$f$ at the breaking points $b_i$ are getting `closer and closer'
to the values of the function $y=x^r+1$). Then $\g (f)=\g (y)=r$.
For each $n\geq 1$, let $I_n=\{ jp^{k_n}, 1\leq j \leq i_n\}$,
$I:= \cup_{n\geq 1}I_n\cup \{ 0\}$, and $I':= \mathbb{N}\backslash
I$. Consider the $\L$-module $M_r:= \L /\L (\der^{[i]}\, | \, i\in
I')$ and its filtration of  standard type $\{ M_{r,i}\}$ induced
from the canonical filtration on the algebra $\L$. Then
$\dim_K(M_{r,j})\leq f(j)$ for all $j\geq 0$, and
$\dim_K(M_{r,b_\nu })=f(b_\nu )$ for all $\nu\geq 1$. Therefore,
$\Dim (M):=\g (\dim_K(M_{r,j}))=\g (f)=r$.

$2$. It follows from $\CM_r =P_1\t M_r$ that $\Dim (\CM_r)=1+r$
since $\dim_K(P_{1, i})=i+1$ is a polynomial.

$\Box $

\begin{theorem}\label{3Jl05}
Let $K$ be a field of characteristic $p>0$. Then $n\leq \Dim
(M)\leq 2n$ for each nonzero finitely generated $\CD (P_n)$-module
$M$, and for each real number $d$ such that $n\leq d\leq 2n$ there
exists a cyclic $\CD (P_n)$-module $M$ such that $\Dim (M)=d$.
\end{theorem}

{\it Proof}. Let $M$ be a nonzero finitely generated $\CD
(P_n)$-module. Then $\Dim (M)\geq n$ by Theorem \ref{c27Jn05}, and
$\Dim (M)\leq \Dim (\CD (P_n))=2n$.  Therefore, $n\leq \Dim
(M)\leq 2n$.

Given a real number $d$ such that $n\leq d\leq 2n$. Then $d=n+s+r$
for some $s\in \mathbb{N}$ and $0\leq r<1$. If $r=0$ then $\Dim
(\CD (P_s)\t P_{n-s})=2s+n-s=d$. If $r\neq 0$ then $\Dim (\CD
(P_s)\t P_{n-s-1}\t \CM_r)=2s+n-s-1+1+r=n+s+r=d$ where the $\CD
(P_1)$-module $\CM_r$ is from Lemma \ref{l3Jl05}. Obviously, the
$\CD (P_n)$-modules $\CD (P_s)\t P_{n-s}$ and $\CD (P_s)\t
P_{n-s-1}\t \CM_r$ are cyclic.  $\Box $

Department of Pure Mathematics

University of  Sheffield

Hicks Building

Sheffield S3~7RH

UK

email: v.bavula@sheffield.ac.uk


\begin{thebibliography}{99}

\bibitem{BavIzv} V. Bavula, Identification of the Hilbert function
and the Poicar\'{e} series, and the dimension of modules over
filtered rings, {\it Russian Acad. Sci. Izv. Math.} {\bf 44}
(1995), 225--246.

\bibitem{Bavcafd} V. Bavula, Filter dimension of algebras and modules, a simplicity
criterion for generalized Weyl algebras. {\it Comm. Algebra} {\bf
24} (1996), no. 6, 1971--1992.

\bibitem{Bavjafd} V. Bavula, Krull, Gelfand-Kirillov, and filter dimensions of simple affine
   algebras. {\it J. Algebra} {\bf  206} (1998), no. 1, 33--39.

\bibitem{bie98} V. Bavula, Krull, Gelfand-Kirillov, filter, faithful and Schur dimensions.
   {\it Infinite length modules} (Bielefeld, 1998), 149--166, Trends Math., Birkhäuser,
   Basel, 2000.

\bibitem{Bav-holmodII} V. Bavula, Dimension, multiplicity, holonomic modules, and an analogue of the inequality of Bernstein
for  rings of differential operators in prime characteristic, II.


\bibitem{BogvadJA95} R. Bogvad, Some results on $\CD $-modules on
Borel varieties in characteristic $p>0$. {\it J. of Algebra} {\bf
173} (1995), 638--667.

\bibitem{Ber72} I.N. Bernstein, Modules over a ring of differential operators.
 An investigation of the fundamental solutions of equations
with constant coefficients. {\it Funkcional. Anal. i Prilozen.}
{\bf 5} (1971),  no. 2, 1--16.

\bibitem{Haastert87} B. Haastert, \"{U}ber Differentialoperatoren und $D$-Moduln in positiver
Charakteristik. {\it Manuscripta Math.} {\bf 58} (1987), no. 4,
385--415.


\bibitem{KL} G. Krause and T.  Lenagan, {\em Growth of algebras and Gelfand-Kirillov
 dimension}. Revised edition. Graduate Studies in Mathematics, 22.
American Mathematical Society, Providence, RI, 2000.

\bibitem{Ma} H. Matsumura, {\em Commutative ring theory.}
  Cambridge University Press,
Cambridge, 1989.



\bibitem{MR} J. C. McConnell and J. C. Robson, Noncommutative Noetherian rings,
Wiley, Chichester,  1987.

\bibitem{Meb-Nar-MacLNM90} Z. Mebkhout and L. Narvaez-Macarro, Sur
les coefficients de Rham-Grothendieck des vari\'et\'es
alg\'ebriques, {\it Lecture Notes in Mathematics}, vol. 1454,
Springer-Verlag, Berlin/New York, 1990.

\bibitem{Pierceb}   R. Pierce, Associative algebras. Springer-Verlag, New York-Berlin, 1982.



\bibitem{KSm-vdB} K. E. Smith and M. van den Bergh, Simplicity
of rings of differential operators in prime characteristic. {\it
Proc. London Math. Soc.} {\bf 75} (1997), no. 1, 32--62.

\bibitem{Smith85LNM} S. P.  Smith, Differential operators on
commutative algebras, in {\em Ring Theory, LNM 1197}, (1985),
164-177.

\bibitem{Smith85LNM} S. P.  Smith, Differential operators on
 the affine and projective lines, in {\em Ring Theory, LNM 1220}, (1985), 157--177.



\end{thebibliography}
\end{document}